\documentclass{article}[12pt]
\usepackage{amsmath,amsfonts,amsthm}
\usepackage{amssymb}

\usepackage{epsfig}

\usepackage{graphics}
\usepackage{graphicx}

\usepackage{latexsym}
\usepackage{graphics}
\usepackage{amsmath}
\usepackage{amsthm}
\usepackage{xspace}
\usepackage{amssymb}
\usepackage{color}
\usepackage{cite}
\usepackage{latexsym}
\usepackage{graphics}
\usepackage{amsmath}
\usepackage{amsthm}
\usepackage{xspace}
\usepackage{amssymb}
\usepackage{epsfig}

\usepackage{enumitem}

\DeclareMathOperator*{\argmax}{arg\,max}
\DeclareMathOperator*{\argmin}{arg\,min}

\newcommand{\beql}[1]{\begin{equation}\label{#1}}
\newcommand{\eeql}{\end{equation}}
\newcommand{\eqn}[1]{(\ref{#1})}

\newcommand{\sasha}[1]{#1}

\newcommand{\R}{\mathbb{R}}
\newcommand{\pr}{\mathbb{P}}
\newcommand{\E}{\mathbb{E}}

\def\P{{\mathbb P}}

\newcommand{\cx}{{\cal X}}

\newcommand{\cj}{{\cal J}}

\newcommand{\bI}{{\bf I}}

\newtheorem{thm}{Theorem}[section]
\newtheorem{lem}{Lemma}[section]

\newtheorem{assumption}{Assumption}[section]
\newtheorem{definition}{Definition}[section]
\newtheorem{condition}[assumption]{Condition}

\begin{document}

\title{Parallel server systems with cancel-on-completion redundancy 
}

\author
{
Alexander L. Stolyar \\
University of Illinois at Urbana-Champaign\\
Urbana, IL 61801, USA \\
\texttt{stolyar@illinois.edu}
}

\date{\today}

\maketitle

\begin{abstract}

We consider a parallel server system with so-called cancel-on-completion redundancy.
There are $n$ servers and
multiple job classes $j$. An arriving class $j$ job consists of $d_j$ components, 
placed on a randomly selected subset of servers;
the job service is complete as soon as $k_j$ components out of $d_j$ (with $k_j \le d_j$) complete their service,
at which point the unfinished service of all remaining $d_j-k_j$ components is canceled. 
The system is in general non-work-conserving, \sasha{in the sense that
the average amount of new workload added to the system by an arriving class $j$ job is not defined a priori -- it depends
on the system state at the time of arrival.}
This poses the main challenge for the system analysis. 

\sasha{For the system with a fixed number of servers $n$ our main results include: the stability properties; 
the property that the stationary distributions of the relative server workloads remain tight, uniformly in the system load.}

\sasha{We also consider the mean-field asymptotic regime when $n\to\infty$ while each job class arrival rate per server remains constant. The main question we address here is: under which conditions the steady-state asymptotic independence (SSAI) of server workloads holds, and in particular when the SSAI for the full range of loads (SSAI-FRL) holds. 
(Informally, SSAI-FRL means that SSAI holds for any system load less than $1$.)
We obtain sufficient conditions for SSAI and SSAI-FRL. In particular, we prove that SSAI-FRL holds in the important special case when job components of each class $j$ are i.i.d. with an increasing-hazard-rate distribution.}

\end{abstract}




\section{Introduction}

We consider a parallel server system with so-called cancel-on-completion (c.o.c.) redundancy (cf. \cite{Vulimiri-2013, shah2015redundant, Harcol-Balter-2017, SnSt2020}), which has been introduced as a means of improving reliability and/or reducing delays in data storage/processing systems. When c.o.c. redundancy is employed, 
an arriving job consists of multiple components, placed on a randomly selected subset of servers.
A job service is complete as soon as a certain number of components out of total complete their service,
at which point the unfinished service of all remaining components 
is immediately canceled. 
The key property of this model, which poses main challenges for its analysis, is that it is in general non-work-conserving, in the sense that
the average amount of new workload added to the system by an arriving job is not defined a priori; therefore, for example, the system 
load is not defined a priori. The results of this paper concern both the system with fixed number of servers and the mean-field asymptotic regime when the number of servers goes to infinity while the job arrival rate per server remains constant. 
The main question we address for the asymptotic regime is whether the steady-state asymptotic independence of server workloads holds.

A more specific definition of the model is as follows. There are $n$ identical servers, processing work at unit rate. New jobs arrive as a Poisson process of rate $\lambda^n n$. There is a finite set of job classes, and an arriving job is of class $j$ with probability $\pi_j^n$.
A class $j$ job consists of $d_j$ components, which are placed on $d_j$ servers selected uniformly at random;
the components' sizes (workloads) $(\xi^{(j)}_1,\ldots, \xi^{(j)}_{d_j})$ are drawn according to an exchangeable distribution $F_j$.
Each server processes its work (components of different jobs) in the First-Come-First-Served order. 
A class $j$ job service is complete as soon as $k_j$ components (with $k_j \le d_j$) of the job complete their service, at which point the unfinished service of the remaining $d_j -k_j$ components of the job is ``canceled;'' 
we will call this $(d_j,k_j)$-c.o.c. redundancy.
We study both the system with a fixed $n$ and the mean-field asymptotic regime when $n\to\infty$, $\lambda^n \to \lambda$,
and $\pi^n_j \to \pi_j$ for all $j$. 

\sasha{Due to cancellations, unless $k_j = d_j$, the actual amount of workload that a job component adds to the server it is placed on, may be smaller than the component size. As a result, the system is in general {\em non-work-conserving} in the following sense: both the distribution and {\em the expectation} of
the random total workload that a class $j$ job actually
brings to the system {\em depend of the system state at the time of arrival.} 
Due to the non-work-conservation property of c.o.c redundancy, even the system stochastic stability/instability in general is not known a priori.
Furthermore, the non-work-conservation renders some very intuitive qualitative properties 
non-obvious for a c.o.c. redundancy system.}
For example, consider a system with fixed finite $n$. Let $\bar \lambda^n$ be the supremum of those $\lambda^n$ for which the system 
is stable. Then it is natural to expect that, as $\lambda^n \uparrow \bar \lambda^n$, the system steady-state load $\rho^n(\lambda^n)$
(i.e. the probability that a server is busy) increases to $1$. 
This property, which we refer to as {\bf stability for the full range of loads (stability-FRL)}, 
is {\em not} automatic for our system. But, we prove that it does, in fact, hold for our model.
\sasha{(In Section~\ref{sec-work-discuss} below we discuss the non-work-conservation and related challenges in more detail.)}

Another property that we prove for the systems with fixed $n$ is the
{\bf uniform (in $n$ and $\lambda^n < \bar \lambda^n$) tightness of stationary distributions of the relative server workloads.} 
Here the 'relative server workloads' refers to the empirical distribution of server workloads, centered by their average. 
Moreover, we show that, essentially, this uniform tightness extends to ``free'' systems as well, in which the ``server workloads''
are {\em not} lower bounded by (regulated at) $0$; we also discuss a connection between free systems and stability of our original systems.
These fixed-$n$ results are both of independent interest and they serve as important ``ingredients'' of the asymptotic analysis. 

Most of this paper is devoted to the mean-field asymptotic regime, with the main property of interest for a given $\lambda$ being
{\bf the steady-state asymptotic independence (SSAI)} of server workloads, namely the property that, as $n\to\infty$
(with $\lambda^n \to \lambda$, $\pi^n_j \to \pi_j, ~\forall j$), the steady-state workloads of servers within a fixed finite set become independent.
Suppose the SSAI holds for some values of $\lambda$
(we will prove that it always holds 
for sufficiently small positive $\lambda$) , 
and denote by $\bar \lambda$ the supremum of those values.
Denote by $\rho(\lambda)= \lim_n \rho^n(\lambda^n)$ the limiting system load for $\lambda < \bar \lambda$.
It is natural to expect that as $\lambda \uparrow \bar \lambda$, the system limiting steady-state load $\rho(\lambda) \uparrow 1$.
This property, which will be referred to as {\bf the steady-state asymptotic independence for the full range of loads (SSAI-FRL)},
is {\em not} automatic for our system. We prove that SSAI-FRL does hold under a certain uniqueness condition on the fixed points of the process mean-field limits. We further identify some sufficient conditions for that uniqueness to hold. In particular,
we prove that SSAI-FRL holds in the important special case
when each job class has i.i.d. component sizes with an increasing-hazard-rate (IHR) distribution. 
\sasha{(Whether or not there exist cases when SSAI-FRL does {\em not} hold for the model in this paper, is an open problem and is a subject of future research.)}

The analysis of mean-field limits and their fixed points is a key part of the proofs of our asymptotic results. 
This part, as well as our general approach, are quite generic, relying almost exclusively on {\em two} fundamental properties of the model: {\bf (a) monotonicity} and {\bf (b) the property that, on average, ``the servers with lower workload receive larger expected additional workload of arriving jobs.''}
This is in contrast to paper \cite{SnSt2020}, which also studies a model with multi-component jobs, satisfying properties (a) and (b) above, but additionally satisfying the work-conservation; the analysis in \cite{SnSt2020} relies on work-conservation in the crucial way. 

The summary of our main results is as follows.

\begin{itemize} 

\item Finite $n$ model.

\begin{itemize} 

\item We prove the stability-FRL property (Theorem~\ref{thm-crit-load-finite}).

\item We prove uniform tightness of stationary distributions of relative server workloads (Theorem~\ref{th-closeness}). Moreover, we show that this uniform tightness extends to ``free'' systems as well (Theorem~\ref{th-closeness-free}),
and discuss a connection between free systems and stability.

\end{itemize}

\item Asymptotic regime: $n\to\infty$, $\lambda^n \to \lambda$, $\pi^n_j \to \pi_j, ~\forall j$. 

\begin{itemize} 

\item We derive general properties of mean-field limits and their fixed points, under general assumptions on the components' sizes distributions (these are the results in Section~\ref{sec-ml}). In particular, we prove that fixed points are completely characterized as solutions of a certain functional differential equation (Lemma~\ref{lem-fp-iff-de}).

\item For the special, ``truncated''  (``finite-frame'') model, where the workloads are ``clipped'' at some finite level, 
we prove that SSAI-FRL always holds (Theorem~\ref{thm-main-finite}).

\item We prove SSAI-FRL under a certain uniqueness assumption on the mean-field limit fixed points (Theorem~\ref{thm-main}). We identify sufficient conditions for that uniqueness to hold
 (Theorem~\ref{thm-main-D-monotonicity}). This, in particular, allows us to prove 
SSAI-FRL in the important special case (Theorem~\ref{thm-ihr})
when each job class has i.i.d. IHR component sizes (or, more generally, the joint component size distribution is a mixture of such distributions) .

\item We prove some other sufficient conditions for SSAI to hold for a given $\lambda$
(Theorems~\ref{thm-majorization}, \ref{thm-inh-subcrit}, \ref{thm-main-majorize}). 
In particular, SSAI holds for any $\lambda$ such that the system is inherently subcritically loaded (Theorem~\ref{thm-inh-subcrit}), and we provide non-trivial conditions to verify inherent subcriticality.  

\end{itemize}

\end{itemize}

Prior work specifically on c.o.c. redundancy models includes \cite{Vulimiri-2013, shah2015redundant, Harcol-Balter-2017}. 
Paper \cite{hellemans2019performance} considers redundancy models more general than c.o.c. 
Cancel-on-completion is not the only form of redundancy considered in prior literature -- another form being 
cancel-on-start (c.o.s.) redundancy, when the remaining $d_j-k_j$ job components are canceled as soon as $k_j$ (out of $d_j$)
components {\em start} their service; under the c.o.s. redundancy the system is work-conserving. 
(The extensively studied power-of-$d$-choices load balancing scheme \cite{VDK96,BLP2012-jsq-asymp-indep} is a special case of c.o.s. redundancy, with 
a single job class with $k_j=1 < d_j$.)
For a more extensive overview of the literature on redundancy models, which includes papers
\cite{Vulimiri-2013,shah2015redundant,gardner2015reducing,Harcol-Balter-2017,adan2018fcfs,ayesta2018unifying,ayesta2019redundancy,VDK96,BLP2012-jsq-asymp-indep},
we refer reader to \cite{SnSt2020} and references therein. (We note that, in most of the prior work, what we call job components are called job ``replicas.'')

Paper \cite{SnSt2020} proves SSAI, under any subcritical load,
for a {\em work-conserving} system with multi-component jobs,
such that each job class may be of the two different types -- water-filling or least-load. 
(Given that the system is work-conserving, this automatically implies SSAI-FRL.)
The least-load job classes cover the c.o.s. redundancy. The water-filling job classes cover, in particular, the c.o.c. redundancy, {\em but only in the special case of i.i.d. exponentially distributed component sizes,}
because then a c.o.c. job can be equivalently viewed as a water-filling job. 

\sasha{Besides the above (work-conserving) special case covered in \cite{SnSt2020}, the author is not aware of any prior work for c.o.c. redundancy models, proving the SSAI-FRL, or even identifying and proving non-trivial sufficient conditions for the SSAI. This is the main distinction of the present paper from the prior work.}

While there are essentially no prior works {\em proving} SSAI for c.o.c. models, 
the {\em SSAI conjecture} is often employed to obtain estimates of the system performance metrics when the number of servers is large (cf. \cite{Vulimiri-2013,Harcol-Balter-2017,hellemans2019performance}). 

\subsection{\sasha{Discussion of non-work-conservation and its implications}}
\label{sec-work-discuss} 

In queueing literature work-conservation usually refers to some form of the property that work 
``cannot be discarded'' by a system. One way to formally define this property is as follows: the distribution of the random total amount of work brought to the system by a job of given class $j$ is independent of the system state at the time of arrival, and the job leaves the system only when the entire amount of its associated work is completed by the system. 
Let us refer to this as ``work-conservation in the strong sense.'' 
In this paper, we understand work-conservation as the following weaker property: the 
{\em expected} total amount of work $s_j$ brought to the system by a job of given class $j$ does not depend on the system state at the time of arrival, and the job leaves the system only when the entire amount of its associated work is completed by the system. 

As we already noted, due to cancellations, the system under c.o.c. redundancy is in general non-work-conserving (even in the specified weaker sense). This is contrast to, for example, a system with cancel-on-start (c.o.s.) redundancy \cite{SnSt2020}, 
described above, 
where a class $j$ job brings exactly $s_j=k_j \E \xi^{(j)}_1$ expected new workload, thus
making a c.o.s. system work-conserving. (In the degenerate case when $k_j = d_j$ for {\em all} job classes $j$, there is no difference between c.o.c. and c.o.s. redundancies, because there are no cancellations.) 

There are cases when a system, which is ostensibly non-work-conserving, in fact is work-conserving. An example is the system with c.o.c. redundancy in the special case when, for each job class $j$, the component sizes are i.i.d. exponential. 
In this case it is easy to observe  that the total amount of work brought by a class $j$ job
happens to be equal in distribution to $\sum_1^{k_j} \xi^{(j)}_i$, regardless of the system state at the time of arrival;
this is because, in this case, a c.o.c. job can be equivalently viewed as a water-filling job (see \cite{SnSt2020}).
Therefore, $s_j = k_j \E \xi^{(j)}_1$,  and the system is indeed work-conserving. (In fact, this system is work-conserving even in the strong sense.) 

Consider stability properties of a multi-component-job system with a given $n$. 
Assume that if the system is stable for $\lambda^n$, then it is also stable for any smaller arrival rate parameter. (As we will see, this fact is straightforward for the model in this paper, as well as for many other models,
including that in \cite{SnSt2020}.)
Then the property which we called {\em the stability for the full range of loads} (stability-FRL) can be thought of as a generalization of the property that `the system is stable as long as its nominal load is less than $1$'
for a work-conserving system. Indeed, 
for a work-conserving system, there is a well-defined nominal system load, namely the steady-state probability that a server is busy, {\em assuming stability holds}. Specifically, the nominal load is simply 
$\rho^n(\lambda^n) = [\sum_j \pi_j^n s_j] \lambda^n$. Then the stability-FRL is equivalent to the `stability 
as long as $\rho^n(\lambda^n) < 1$, or $\lambda^n < \bar \lambda^n = 1/[\sum_j \pi_j^n s_j]$.'
For many work-conserving systems, including a quite general system in \cite{SnSt2020}, 
the stability for $\lambda^n < \bar \lambda^n$, and then the stability-FRL,
is easy to show. Let us now consider a non-work-conserving system, specifically the c.o.c. redundancy 
system in the present paper. In this case, there is no a priori defined nominal load $\rho^n(\lambda^n)$. Instead, as we increase
$\lambda^n$, and as long as the system remains stable, actual load $\rho^n(\lambda^n)$ is some, generally non-linear, unknown a priori, function of $\lambda^n$. As $\lambda^n$ approaches $\bar \lambda^n$, which is 
the supremum of those $\lambda^n$ for which the system 
is stable, the property that $\rho^n(\lambda^n)\uparrow 1$, i.e. the stability-FRL, is not automatic.
We do prove this property for our system (and, in fact, give a characterization of $\bar \lambda^n$, although not explicit).

Let us now discuss the asymptotic regime, $\lambda^n \to \lambda$, $\pi_j^n \to \pi_j, ~\forall j$, and the properties which we called {\em the steady-state asymptotic independence} (SSAI) and {\em the steady-state asymptotic independence for the full range of loads} (SSAI-FRL). The SSAI is equivalent to the property that, as $n\to\infty$, the steady-state (random) empirical distribution of the server workloads in the system converges to the single (deterministic) point $x^{*,\lambda}  = (x_w^{*,\lambda}, ~w\ge 0),$ where $x_w^{*,\lambda}$ gives the (limiting) fraction of servers with workload greater than $w$. 
In the terminology that we will introduce later, $x^{*,\lambda}$ is a {\em fixed point of the mean-field limit} (ML-FP). 
Consider a work-conserving system. Then, SSAI-FRL is equivalent to the property that `SSAI holds for all $\lambda$ such that the (limiting) nominal load $\rho(\lambda) = [\sum_j \pi_j s_j] \lambda < 1$, or $\lambda < \bar \lambda = 1/[\sum_j \pi_j s_j]$.'  The proof of this property for a (work-conserving) model in \cite{SnSt2020} relies in essential way on the
the facts that an ML-FP $x^{*,\lambda}$ is unique and is such that the (limiting) fraction of busy servers is exactly $\rho$, i.e.
$x^{*,\lambda}_0 = \rho$;
these facts are relatively easy consequences of work-conservation. Let us now consider the non-work-conserving model in this paper. Proving SSAI-FRL reduces to proving the following for any fixed $\rho<1$. 
Consider a converging sequence $\lambda^n \to \lambda$ 
such that, for all $n$, $\rho^n(\lambda^n) = \rho$. (Such a sequence exists, because for each $n$ we have the stability-FRL.)
Then the steady-state empirical distribution of the server workloads converges to a single (deterministic) point $x^{*,\lambda}$,
which is an ML-FP with $x^{*,\lambda}_0 = \rho$. The key difficulty posed by the non-work-conservation 
is that, while the existence of an ML-FP, with $x^{*,\lambda}_0 \le \rho$, is not hard to obtain, it is not clear that an ML-FP 
for a given $\lambda$ is unique. In this paper we prove that, in essence, {\em if the uniqueness of an ML-FP for each $\lambda$
holds, then the SSAI-FRL holds as well.} (This fact -- which is the essence of Theorem~\ref{thm-main} -- is far from straightforward. Large part of this paper is devoted to and leads to the proof of Theorem~\ref{thm-main}; an informal discussion
of the proof key ideas is given at the beginning of Section~\ref{sec-proof-main-th}.
Also, Theorem~\ref{thm-main} is quite generic, not relying much on the specific structure of the c.o.c. redundancy model, beyond
the two fundamental properties (a) and (b) mentioned earlier.) We further show that the ML-FP uniqueness condition, and then the SSAI-FRL, hold in the important special case of the i.i.d. IHR components. (This part relies on the c.o.c. redundancy model
structure to a much larger degree.)

Finally, we note that the challenge of establishing SSAI-FRL is not limited to specifically multi-component-job models or specifically non-work-conservation. As a example, a join-the-idle-queue (JIQ) model with general job size distribution, studied 
in \cite{FS2016}, is work-conserving and is not a multi-component-job model. 
The stability-FRL for this model is straightforward. However, the SSAI for it has only been established for the loads less than $1/2$, as opposed to less than $1$. Therefore, SSAI-FRL for this system is still an open problem, to the best of our knowledge.
We note that the JIQ model in \cite{FS2016} does not have certain monotonicity properties, which the model 
in our paper does have.


\subsection{Paper layout}
\label{sec-layout} 

The rest of the paper is organized as follows. Basic notation and terminology, used throughout the paper, are introduced in Section~\ref{sec-basic-notation}. Section~\ref{sec-formal-model} formally defines our model and the asymptotic regime. 
Sections~\ref{sec-results-finite} and \ref{sec-results-infinite} state our main results for the truncated (finite-frame) and the original (infinite-frame) models, respectively. In Sections~\ref{sec-monotonicity-etc} and \ref{sec-equiv-view} we establish some monotonicity -- and related to it -- properties of the model,
which serve as important tools of our analysis. 
Sections~\ref{sec-proof444} and \ref{sec-closeness-proof} contain the proofs of 
Theorems~\ref{thm-crit-load-finite} and \ref{th-closeness}, respectively, which are our main results for a system with finite $n$.
In Section~\ref{sec-ml} we define and study the mean-field limits, which serve as a key tool of our asymptotic analysis.
Sections~\ref{sec-proof-main-th}-\ref{sec-proof-ihr} present proofs of all asymptotic results stated in Sections~\ref{sec-results-finite} and \ref{sec-results-infinite}. 
In Section~\ref{sec-further-conditions} we state and prove some additional results on SSAI.
Finally, in Section~\ref{sec-free} we state and prove the tightness result for the ``free'' system.

\section{Basic notation and terminology}
\label{sec-basic-notation}

We denote by $\R$ and $\R_+$ the sets of real and real non-negative numbers, respectively, and by $\R^n$ and $\R_+^n$ the corresponding $n$-dimensional product sets. By $\bar \R \doteq \R \cup \{\infty\} \cup \{-\infty\}$ we denote the two-point compactification of $\R$, where 
$\infty$ and $-\infty$ are the points at infinity and minus infinity, with the natural topology and a consistent with it metric,
so that $\bar \R$ [and $\R$] is complete and separable.
Analogously,
$\bar \R_+ \doteq \R_+ \cup \{\infty\}$ is the one-point compactification of $\R_+$. 
For a vector $q=(q_i) \in \R^n$,  $|q| \doteq \sum_i |q_i|$.
Inequalities applied to vectors [resp.  functions] are understood component-wise [resp. for every value of the function argument]. 
We say that a function is RCLL [resp. LCRL] 
if it is {\em right-continuous with left-limits} [resp. {\em left-continuous with right-limits}]. 
A scalar function $f(w)$ is called $c$-Lipschitz if it is Lipschitz continuous with constant $c$.
The sup-norm of a scalar function $f(w)$ is denoted $\|f(\cdot)\| \doteq \sup_w |f(w)|$; the corresponding convergence is denoted by $\stackrel{\|\cdot\|}{\longrightarrow}$. {\em U.o.c.} convergence means {\em uniform on compact sets} convergence, and is denoted by $\stackrel{u.o.c.}{\longrightarrow}$. We use notation: $a\vee b \doteq \max\{a,b\}$, 
$a\wedge b \doteq \min\{a,b\}$.
Abbreviation {\em w.r.t.} means {\em with respect to};
{\em a.e.} means {\em almost everywhere w.r.t. Lebesgue measure};
WLOG means {\em without loss of generality}; RHS and LHS means {\em right-hand side} and {\em left-hand side}, 
respectively.

We denote by $\tilde \cx$ [resp., $\cx$] the set of non-increasing RCLL functions $x= (x_w, ~w\in \R)$, [resp., $x=(x_w, ~w\in \R_+)$,]
taking values in $[0,1]$. An element $x \in \tilde \cx$ [resp., $x\in \cx$] is naturally interpreted as complementary distribution function
on $\bar \R$ [resp., $\bar \R_+$], with $1-x_w$ being the measure of $[-\infty,w]$ [resp., $[0,w]$]. An element $x \in \tilde \cx$ is {\em proper} if $x_\infty \doteq \lim_{w\to \infty} x_w =0$ and $x_{-\infty} \doteq \lim_{w\to -\infty} x_w =1$, and {\em improper} otherwise; in any case, $x_\infty$ and $1 - x_{-\infty}$ give the measures of $\{+\infty\}$ and $\{-\infty\}$, respectively. Analogously, an element 
$x \in \cx$ is proper if $x_\infty =0$ and  and improper otherwise. 
(With some abuse of standard terminology, in this paper we sometimes refer to elements of  $\tilde \cx$ and $\cx$ as distributions.)
Denote by $\tilde \cx^{pr}$ [resp., $\cx^{pr}$] the subset of proper elements of $\tilde \cx$ [resp., $\cx$].

We will equip the space $\tilde \cx$ with the topology of weak convergence of the corresponding distributions on $\bar \R$;
equivalently, $y \to x$ if and only if $y_w \to x_w$ for each $w \in (-\infty,\infty)$ where $x$ is continuous.
\sasha{Furthermore, on $\tilde \cx$ we consider the following metric $L$, consistent with the weak convergence topology.
We map an element $x \in \tilde \cx$ into the (proper or non-proper) distribution function $\phi=\phi(x)$ on $\R$, where
$\phi_s = 1- x_{-\log (1-s)}$ for $0 \le s < 1$, $\phi_s = 1- x_{\log (1+s)}$  for $-1 \le s < 0$, $\phi_s = 0$ for $s<-1$, $\phi_s = 1$ for $s\ge 1$.
Then $L(x,y) \doteq \hat L(\phi(x),\phi(y))$, where $\hat L$ is the Levy-Prohorov metric (cf.  \cite{Ethier_Kurtz}). 
It is easy to see that the convergence in metric $L$ is indeed equivalent to the weak convergence.}
Clearly, $\tilde \cx$ is compact.
Space $\cx$ inherits the topology and metric of space $\tilde \cx$ -- for these purposes we view any $x \in \cx$ as 
an element of $\tilde \cx$ via the convention that $x_w = 1$ for $w<0$. Clearly, $\cx$ is also compact.

For a function $x \in \tilde \cx$ [resp., $x\in \cx$] we denote by $x^{-1}=(x_u^{-1}, ~u \in [0,1]),$
its inverse, defined as
$$
x_u^{-1} = \sup \{w~|~x_w > u\},
$$
with the convention that $x_u^{-1} = -\infty$ [resp., $x_u^{-1} = 0$]
if the set in $\sup$ is empty. 

Unless explicitly specified otherwise, we use the following conventions regarding random elements and random processes.
A measurable space is considered equipped with a Borel $\sigma$-algebra, induced by the 
metric which is clear from the context. A random process $Y(t), ~t\ge 0,$ always takes values in a complete separable metric space (clear
from the context), and has RCLL sample paths; the sample paths are elements of the Skorohod $J_1$-space with the corresponding metric.

For a random process $Y(t), ~t\ge 0,$ we denote by $Y(\infty)$ the random value of $Y(t)$ in a stationary regime (which will be clear from the context). Symbol $\Rightarrow$ signifies convergence of random elements in distribution; $\stackrel{P}{\longrightarrow}$ means convergence in probability.
 {\em W.p.1} means {\em with probability one.}
{\em I.i.d.} means {\em independent identically distributed.}
Indicator of event or condition $B$ is denoted by $\bI(B)$. If $X,Y$ are random elements taking values in a set on which a partial order $\le$ is defined, then the stochastic order
$X \le_{st} Y$ means that $X$ and $Y$ can be coupled (constructed on a common probability space) so that $X \le Y$ w.p.1.

 \subsection{Fluid limits}
 
 In several places we rely on the {\em fluid limit technique} \cite{RS92, Dai95, St95, Bramson-book}
 to establish stability (positive recurrence) of Markov processes. We now fix the corresponding definitions and terminology,
 which are not the most general, but suffice for the purposes of this paper.
 
Consider a continuous time Markov process $Q(t)=(Q_1(t), \ldots, Q_n(t)), ~t\ge 0,$ in $\R^n$. 
Consider a sequence of processes $Q^{(r)}(\cdot)$, indexed by by $r\uparrow \infty$, with initial states such that 
$|Q^{(r)}(0)|=r$ and $(1/r) Q^{(r)}(0) \to q(0)$, for some fixed $q(0) \in \R^n$, $|q(0)| = 1$. 
We say that {\em fluid limits $q(\cdot)$ of $Q(\cdot)$ satisfy a set of properties $A$} if trajectories $q(t), ~t\ge 0$, of any process
being a distributional limit of the sequence of processes $(1/r) Q^{(r)}(rt), ~t\ge 0,$ along a subsequence of $r \uparrow \infty$,
satisfy properties $A$ w.p.1.
If any such distributional limit is concentrated on the unique trajectory $q(\cdot)$ (for a fixed $q(0)$), we say that the fluid limit 
is unique and is equal to this $q(\cdot)$. If for some fixed $T>0$, fluid limits satisfy property
$$
|q(t)| = 0, ~~\forall t \ge T,
$$
we say that fluid limits are stable. The fluid limit technique 
allows one to claim that 
``typically'' (under mild additional conditions) the stability of fluid limits implies stability (positive recurrence) of the Markov process itself.

\section{Model and asymptotic regime}
\label{sec-formal-model}

\subsection{Model. Cancel-on-completion redundancy}

There are $n$ identical servers, each processing its work at rate $1$.
The {\em workload} of a server at a given time is its unfinished work,
i.e. the time duration until the server becomes idle assuming no new job arrivals. 
There is a finite set $\cj$
of job classes $j$. 
Jobs  arrive as Poisson process of rate $\lambda^n n$, for some $\lambda^n \ge 0$.
The probability that an arriving job is of class $j$ is $\pi_j^n \ge 0$, $\sum_j \pi_j^n = 1$.
(Our asymptotic regime, which will be introduced shortly, will be such that $n \to \infty$, $\lambda^n \to \lambda$
and each $\pi_j^n \to \pi_j$.
That's why we have notations $\lambda^n$ and $\pi_j^n$ indicating the dependence on $n$.)
Each job class $j$ has three parameters: integers $k_j$ and $d_j$
such that $1\le k_j \le d_j$, and the exchangeable probability distribution $F_j$ on $\R_+^{d_j}$. 
(Exchangeability of $F_j$ means that it is invariant w.r.t. permutations of components.)
When a class $j$ job arrives, $d_j$ servers are selected uniformly at random (without replacement); these servers form the {\em selection set} of the job. The job places $d_j$ {\em components} on the selected servers; the component {\em sizes} (workloads) 
$(\xi^{(j)}_1,\ldots, \xi^{(j)}_{d_j})$ are drawn according to distribution $F_j$, independently of the process history up to the job arrival time.
Each server processes its work (components of different jobs) in the First-Come-First-Served (FCFS) order. 
A job of class $j$, to be completed, requires $k_j$ (out of $d_j$) components to be processed, and as soon as $k_j$ components of the job {\em complete} their service (i.e. receive the amounts of service equal to their sizes), the unfinished service of the remaining $d_j -k_j$ components of the job is ``canceled'' and 
the job immediately leaves the system.
(Hence the name cancel-on-completion.) We will call this {\em $(d_j,k_j)$-c.o.c. redundancy}.
Throughout the paper we use notation $\bar d \doteq \max_j d_j$.

\sasha{Consider a server at time $t$. Let $t+w$ be the time at which all the components (of different jobs) that are currently (at time $t$) in its queue
leave the system, either due to their service completions or service cancellations. This $w$ we will call the {\em server workload} at time $t$.
In other words, the server workload $w$ at time $t$ is the actual amount of its current unfinished work.
Now consider a class $j$ job arrival at time $t$.  Note that the workloads added 
to the selected servers by this arriving job 
depend not only on the realization of the component sizes,
but also on the current workloads of the selected servers or, more specifically, on their {\em relative} workloads (i.e. workload differences).
For a simple illustration, 
suppose $d_j=3$, $k_j=2$.
Suppose a class $j$ job arrives at time $t$ with the realization of the components' sizes being $(7,4,3)$, and these components are placed on the selected severs $(S_1, S_2, S_3)$
with workloads at time $t$ being $(5,3,6)$. Then, after this job arrival, the new workloads of the selected servers will be
$(9,7,9)$, because the service of the component placed on server $S_1$ will be cancelled at time $t+9$ upon the completion 
of the service of the component placed on server $S_3$. 
(Without cancellation, the new workloads would be $(12,7,9)$.)
Therefore, the workloads added to the selected servers
by this job arrival will be $(4,4,3)$, as opposed to the component sizes  $(7,4,3)$. 
We see that, indeed, 
unless $k_j = d_j$, the realization 
of the added workloads 
depends on the relative workloads of selected servers, as well as the 
realization of the component sizes. Consequently, in general, the joint distribution of the random 
added workloads that the job brings to the selected servers, depends on the relative workloads of those servers.
Moreover, both the distribution and the expectation of the total amount of added workload depend on the relative workloads
of the selected servers. In this sense, the c.o.c. redundancy model is {\em non-work-conserving}. 
In particular, the system load -- the average rate per server at which the new work arrives -- is not known a priori.}

In this paper we will often use ``particle'' language to describe the system dynamics. Namely, we will identify each server with a ``particle,'' the server workload with ``particle location,''  and the workload evolution with ``particle movement.''
The basic dynamics of a particle (server workload) is as follows. If/when a job arrival adds workload $\kappa \ge 0$ to the server,
the particle jumps ``right'' by $\kappa$. Between job arrivals, the particle (workload) moves ``left'' at the constant speed $-1$ until and unless it ``hits'' the ``regulation boundary'' at $0$; if the particle does hit the boundary $0$, it stays at it until and unless it jumps right due to a job arrival.  

Throughout this paper we assume that the component sizes have finite mean (Asssumption~\ref{cond-finite-mean}) 
and that a non-triviality Assumption~\ref{cond-non-triv1} holds.

\begin{assumption}
\label{cond-finite-mean}
For all $j$, 
\beql{eq-finite-mean}
\E \xi^{(j)}_1 < \infty ~~ \mbox{(and then $\E \sum_{i=1}^{d_j} \xi^{(j)}_i = d_j \E \xi^{(j)}_1 < \infty$),}
\eeql
and WLOG, $\E \xi^{(j)}_1 >0$.
\end{assumption}

\begin{assumption}
\label{cond-non-triv1}
There exists class $j$ with $F_j$ such that
\beql{eq-non-triv1}
\pr \{ \sum_i \bI(\xi^{(j)}_i =0) \ge k_j \} < 1.
\eeql
\end{assumption}

Assumption~ \ref{cond-non-triv1} implies that, uniformly on the system state at the time of a class $j$ job arrival, the expected amount of 
workload that the job brings to the system is at least some $\gamma>0$. If Assumption~ \eqn{cond-non-triv1} does {\em not} hold, the system behavior is degenerate in that when the system reaches the state with
all workloads equal zero, it can never leave this state.

Note that in the special case when for each class $j$ either $k_j=d_j$ or the component sizes are i.i.d. exponentially distributed,
the model is within the framework of the model in \cite{SnSt2020}, with job classes being of either the least-load or water-filling type
(in the terminology of \cite{SnSt2020}).

In addition to the basic model defined above in this subsection, we will also consider its truncated version, 
defined as follows. Let a constant $c \ge 0$ be fixed. The model is the same as the basic one, except after an arriving job adds workloads to its selected servers, the server workloads which happen to exceed $c$ are immediately reset (truncated) to $c$. Note that the basic model can be viewed as a special case of the truncated one, with $c=\infty$. From now on we refer to parameter $c$ as the {\em frame size}; 
to the basic (non-truncated) model as the {\em infinite-frame model}, with frame $[0,\infty)$; and to the truncated model with $c<\infty$ as the {\em finite-frame} model, with frame $[0,c]$. The purpose of considering the finite-frame model is two-fold: it is of independent interest and, more importantly, it serves as a tool for the analysis of the infinite-frame model.
 
 \subsection{Asymptotic regime. Mean-field scaled process} 
 \label{sec-asymp-regime}

 We consider the sequence of systems with $n\to\infty$, $\lambda^n \to \lambda \ge 0$, 
 and $\pi_j^n \to \pi_j$ for all $j$,
 while all other system parameters remain fixed. 
 
From now on, the upper index $n$ of a variable/quantity will indicate that it pertains to the system with $n$ servers,
or $n$-th system.
Let $W_i^n(t)$ denote the workload
of server $i$ at time $t$ in the $n$-th system. (When $W_i^n(t)=0$ we say that server $i$ at time $t$ is empty.)
Clearly, for each $n$ the process $W^n(t) = (W^n_1(t), \ldots, W^n_n(t)), ~t\ge 0,$ is Markov with state space $\R_+^n$.
We adopt the convention that its sample paths are RCLL. 
The process state $(0,\ldots,0)$, i.e. such that all servers are empty, we will call the empty state.
\sasha{The process $W^n(\cdot)$ is renewal, with the empty state serving as a renewal atom, reachable from any other state;
the renewals occur when the process enters the empty state.
The process is irreducible, which follows from the reachability of the atom and the fact that the job arrival process is Poisson.}

Consider also the following {\em mean-field} scaled quantities:
\beql{eq-x-def2}
x^n_w(t) \doteq (1/n) \sum_i \bI\{W_i^n(t)> w\}, ~~ w \ge 0.
\end{equation}
That is, $x^n_w(t)$ is the fraction of servers $i$ with $W_i^n(t)> w$.
Then
 $x^n(t)=(x^n_w(t), ~w\ge 0)$ is a projection of $W^n(t)$, and $x^n(t), ~t\ge 0,$ is also a Markov process.
Note that $x^n_0(t)$ 
is the fraction of busy servers (the instantaneous system load).

Clearly, for any $n$ and $t$, $x^n(t) \in \cx^{pr}$. Denote by $\cx^{(n)} \subset \cx^{pr}$
the state space $\cx^{(n)}$ of the Markov process $x^n(t), ~t\ge 0.$ 
Therefore,
for any $n$, we can and will view
$x^n(t), ~t\ge 0$ as a Markov process with (common) state space 
$\cx^{pr} \subset \cx$,
and with sample paths being RCLL functions of $t\ge 0$ (taking values in $\cx^{pr}$).
It is a renewal irreducible process, with renewal atom being the empty state $x^\emptyset$, defined by 
\beql{eq-empty-def}
x^\emptyset_0=0.
\eeql

We say that the process $W^n(\cdot)$ [resp. $x^n(\cdot)$] is {\em stable} if it is {\em positive Harris recurrent}, which for this process
simply means 
that the empty state is reachable from any other state w.p.1 and the 
expected time to return to the empty state (after leaving it) is finite. 
\sasha{(For a general definition of positive Harris recurrence, cf. \cite{Dai95, Bramson-book}.)}
Obviously, $W^n(\cdot)$ is stable if and only if $x^n(\cdot)$ is.
Note that, trivially, a finite-frame system is always stable. For the infinite-frame system,
due to the fact that the system is non-work-conserving, the stability condition is not automatic and is non-trivial.

Stability of the process $W^n(\cdot)$ [resp. $x^n(\cdot)$] implies that it has unique stationary distribution. If the process is stable,
let $W^n(\infty)$ [resp. $x^n(\infty)$] be a random element with values in $\R^n$ [resp., in $\cx^{pr} \subset \cx$],
whose distribution is the stationary distribution of the process; in other words,
it is a random process state in the stationary regime. 
\sasha{(Due to the process renewal structure, with absolutely continuous distribution of a renewal cycle duration,
the stability also implies the convergence $W^n(t) \Rightarrow W^n(\infty)$ [or $x^n(t) \Rightarrow x^n(\infty)$]
starting any proper initial state.)}
We also adopt a convention that $W^{n}(\infty) = (\infty,\ldots,\infty)$
for an unstable process. 

Later in the paper we will need some additional notation, associated with the space $\cx$.
The infinite state 
$x^{**} \in \cx$ is defined by
\beql{eq-inf-def}
x^{**}_\infty=1.
\eeql
For an $x \in \cx$, let us define its {\em frame size} $c$ as follows:
 $c = \infty$ if $ x_w >0$ for all $w < \infty$, and $c = \min \{w~|~ x_w=0\}$, otherwise.
 We say that $x$ has infinite frame $[0,\infty)$ or finite frame $[0,c]$, when $c=\infty$ and $c<\infty$, respectively.
 Note that $x^\emptyset$ is the only element with zero frame size.
 Let $x^{**,c}$ denote the maximum element with frame size $c \ge 0$: $x^{**,c}_w = 1$ for $w<c$, and $x^{**,c}_c = 0$.

\section{Main result for the finite-frame (truncated) system}
\label{sec-results-finite}

\begin{thm}
\label{thm-main-finite}

Consider a finite-frame (truncated)  system, with frame $[0,c]$, $c<\infty$.

(i) For each $\lambda \in [0,\infty)$, there exists $\rho(\lambda) \in [0,1)$ and 
a unique element $x^{*,\lambda} \in \cx$ with finite frame $[0,c]$, with $x^{*,\lambda}_0 = \rho(\lambda)$, such that,
as $n\to\infty$ with $\lambda^n \to \lambda \ge 0$ and $\pi_j^n \to \pi_j$ for all $j$
(while all other system parameters remain fixed), 
\beql{eq-conv-main-finite}
x^n(\infty) \Rightarrow x^{*,\lambda}. 
\eeql

(ii) Function $\rho(\lambda)$ is a strictly increasing continuous function, mapping $[0,\infty)$ onto $[0,1)$. 
The dependence of element $x^{*,\lambda}$ on $\lambda$ is strictly increasing continuous in $[0,\infty)$.
(Strictly increasing here means that $x^{*,\lambda}_w$ is strictly increasing for each $w \ge 0$.)
Furthermore, $x^{*,\lambda} \to x^{**,c}$ as $\lambda \uparrow \infty$. 

As a corollary of (i), we also have the following {\em steady-state asymptotic independence} property, for any $\lambda \ge 0$ 
and any fixed integer $m \ge 1$:
\beql{eq-indep-main-finite}
[W^n_1(\infty), \ldots, W^n_m(\infty)] \Rightarrow [W^*_1, \ldots, W^*_m], ~~ n\to\infty,
\eeql
where random variables $W^*_i$ are i.i.d., with
$\P\{W^*_i > w\} = x^{*,\lambda}_w, ~w\ge 0$. (This is a corollary, because, due to the symmetry between servers,
in steady-state, the joint distribution of workloads of servers $1,\ldots,m$ is equal to that for a set of $m$ servers chosen 
uniformly at random.)

\end{thm}

The proof of Theorem~\ref{thm-main-finite} is in Section~\ref{sec-proof-main-th-finite}. 

\section{Main results for the infinite-frame (non-truncated) system}
\label{sec-results-infinite}

\subsection{Stability properties for a fixed $n$}

Our first main result for the infinite-frame system, Theorem~\ref{thm-crit-load-finite}, 
concerns its stability properties for each fixed $n$.
For the infinite-frame system with a fixed $n$, consider process $W^{n} (\cdot)$ with different values of $\lambda^n$,
and let us use notation $W^{n,\lambda^n} (\cdot)$ to indicate the dependence on $\lambda^n$.
(Here the value of $\lambda^n$ varies for a fixed $n$, while the probabilities $\pi_j^n$ remain fixed.)

Recall that the system is non-work-conserving.
As a result, given the arrival rate $\lambda^n$, the question of the process stability is non-trivial.
Moreover, even some ``very intuitive'' properties are not automatic. Let us denote by
\beql{eq-rho-full-range}
\rho^n(\lambda^n) \doteq \pr\{W_1^{n,\lambda^n}(\infty) > 0\}
\eeql
the system steady-state load, assuming the system is stable for a given $\lambda^n$.
Then it is natural to expect that as we continuously increase $\lambda^n$ from $0$, the steady-state load $\rho^n(\lambda^n)$
continuously increases from $0$ and $\rho^n(\lambda^n)\uparrow 1$ as $\lambda^n$ approaches some critical value 
$\bar \lambda^n$. 
This property, however, is not automatic. We prove it in the following  Theorem~\ref{thm-crit-load-finite}.
(In addition to being of interest in itself, we will use this property in the proofs of our asymptotic results.) 

Denote 
\beql{eq-crit-load-finite-def}
\bar \lambda^n \doteq \sup \{\lambda^n ~|~ W^{n,\lambda^n}(\cdot) ~\mbox{is stable}\}.
\eeql

\begin{thm}
\label{thm-crit-load-finite}
For a fixed $n$ consider the infinite-frame system and the corresponding process $W^{n,\lambda^n} (\cdot)$ with 
$\lambda^n$ being a parameter, and $\bar \lambda^n$ defined in \eqn{eq-crit-load-finite-def}.
Then, the following holds.\\
(i) $0 < C_1 < \bar \lambda^n < C_2 < \infty$, with $C_1,C_2$ independent of $n$. Process $W^{n,\lambda^n}(\cdot)$ is stable for any $\lambda^n < \bar \lambda^n$.\\
(ii) Function $\rho^n(\lambda^n)$ is a strictly increasing continuous (one-to-one) mapping of $[0,\bar \lambda^n)$ onto $[0,1)$.
\end{thm}

The proof of Theorem~\ref{thm-crit-load-finite} is in Section~\ref{sec-proof444}.

\subsection{Tightness of stationary distributions of centered states}
\label{sec-tightness-centered}

For some (not all) of our results, we will need the following
additional assumptions, namely finite second moments of the component sizes (Assumption~\ref{cond-second-moment}) and 
further (non-restrictive) non-triviality conditions (Assumption~\ref{cond-ntriv-add}).

\begin{assumption}
\label{cond-second-moment}
For any $j$, $\E[\xi_1^{(j)}]^2 < \infty$ (and then $\E[\sum_i \xi_i^{(j)}]^2 < \infty$).
\end{assumption}

\begin{assumption}
\label{cond-ntriv-add}
There exists a class $j$, with $k_j < d_j$ and the joint distribution $F_j$ of the component sizes such that:

(i)  $\P\{ \max_i \xi^{(j)}_i - \xi^{(j)}_{i_1} >0\} > 0$, where $\xi^{(j)}_{i_1}$ is the $k_j$-th smallest among the component sizes $\xi^{(j)}_i$.

(ii) $\pr \{ \sum_i \bI(\xi^{(j)}_i =0) \ge k_j \} < 1$.
\end{assumption}

Assumption~\ref{cond-ntriv-add}(i) holds in most cases of interest. For example, it automatically holds when there 
is a class $j$ with i.i.d. component sizes, with a component size distribution {\em not} concentrated on a single point.
 Assumption~\ref{cond-ntriv-add}(ii) is a slightly stronger version of Assumption~\ref{cond-non-triv1} 
(adding that $k_j < d_j$ holds for the class $j$).  
This assumption automatically holds, for example, when there is a class $j$ with $k_j < d_j$ and i.i.d. component sizes.

For an element $y\in \tilde \cx$, denote by $\bar y$ the mean of the corresponding distribution,
$$
\bar y = \int_{-\infty}^\infty u d[-y_u] = \int_0^\infty y_u du - \int_{-\infty}^0 (1-y_u) du,
$$
with the usual convention that the mean is finite when both integrals in the RHS are finite.
When $\bar y$ is finite, denote by $\mathring y = (\mathring y_u,~u\in \R)$ the centered version of $y$, namely
$$
\mathring y_u = y_{u+\bar y}, ~u\in \R,
$$
and denote
$$
\Phi_\ell(y) = \int_{-\infty}^\infty |u|^\ell d[- \mathring y_u] = \int_{-\infty}^\infty |u-\bar y|^\ell d[-y_u], ~~\ell \ge 1,
$$
$$
\tilde \cx_\ell = \{y\in \tilde \cx ~|~ \mbox{$\bar y$ is finite}, ~ \Phi_\ell(y) < \infty\}, ~~\ell \ge 1.
$$
Denote by $\mathring{\cx}_1$ the subset of those $y \in \tilde \cx_1$ with mean $\bar y =0$.

Recall that $\cx^{(n)} \subset \cx^{pr}$ is the state space of the process $x^n(\cdot)$. 
For any $n$ 
denote by $\mathring{\cx}^{(n)} \subset \tilde \cx^{pr}$ the state space of the process $\mathring x^n(\cdot)$.  
(Process $\mathring x^n(\cdot)$ is not Markov -- it is a projection of Markov process $\mathring x^n(\cdot)$.
Of course, $\mathring{\cx}^{(n)}$ is a projection of $\cx^{(n)}$.)
 Obviously, for any $n$ and any $\ell \ge 1$,
$\mathring{\cx}^{(n)} \subset \mathring{\cx}_1 \cap \tilde \cx_\ell$.

\begin{thm}
\label{th-closeness}
\sasha{Consider the infinite-frame system.} 
Suppose the additional Assumptions~\ref{cond-second-moment} and \ref{cond-ntriv-add}  hold.
Then there exist $\bar C>0$ and $\bar n$ such that, uniformly in $n \ge \bar n$ and $\lambda^n$
such that the process $x^n(\cdot)$ is stable, we have
\beql{eq-tight-uniform}
\E \Phi_1(\mathring x^n(\infty)) \le \bar C.
\eeql
\end{thm}

The proof of Theorem~\ref{th-closeness} is given in Section~\ref{sec-closeness-proof}. 
\sasha{Note that the theorem implies the tightness of the family of distributions of $\mathring x^n(\infty)$,
for $n \ge \bar n$ and those $\lambda^n$ for which stability holds,
 in the space of distributions on $\tilde \cx^{pr}$. 
(Hence the title of this subsection.)
Indeed, for any $\delta>0$ we can choose a constant C large enough, so that $\P \{\Phi_1(\mathring x^n(\infty)) \le C\} \ge 1-\delta$,
uniformly within the specified family; and set 
$\{\Phi_1(x) \le C\}$ is compact in $\tilde \cx^{pr}$.}

\sasha{The proof of Theorem~\ref{th-closeness} uses a quadratic Lyapunov function $\Phi_2(\mathring x^n(t))$, 
which is a function of the set of particle locations centered by their average. This provides the intuition for why the bound \eqn{eq-tight-uniform} -- and the tightness of the family of distributions of $\mathring x^n(\infty)$ -- holds uniformly in $\lambda^n$. 
Informally speaking, when $\Phi_2(\mathring x^n(t))$ is large, it will have a negative drift regardless of how close 
the particle locations are to the boundary at $0$ -- only particle {\em relative} locations matter.
(This, in turn, is due to the fundamental model property (b) referred to in the Introduction: ``the servers with lower workload receive larger expected additional workload of arriving jobs.'')} 
In fact, a version of Theorem~\ref{th-closeness} holds for an artificial ``free'' system, where particle locations (server workloads)
are {\em not} regulated at $0$; namely, they evolve in $\R$ rather than in $\R_+$, with each particle keeping moving left at constant rate $-1$
unless/until it jumps right due to a job arrival.  The corresponding result, Theorem~\ref{th-closeness-free}, is given in Section~\ref{sec-free},
where we also discuss its connection to Theorem~\ref{thm-crit-load-finite}.

\sasha{Assumption~\ref{cond-second-moment} -- the finiteness of second moments of component sizes -- is employed in the proof of Theorem~\ref{th-closeness} in the estimates of the steady-state drift of the quadratic Lyapunov function $\Phi_2(\mathring x^n(t))$. Informally speaking, it guarantees that the expected increment of $\Phi_2(\mathring x^n(t))$
due to a job arrival is equal, up to a bounded additive term, to the first-order approximation of the expected increment;
see \eqn{eq-generator-bound} and its proof in Section~\ref{sec-generator-bound}.}

\subsection{Steady-state asymptotic independence results}

\sasha{The following Theorem~\ref{thm-main} assumes a certain uniqueness condition, namely Condition~\ref{cond-unique}, given later in Section~\ref{sec-cond-unique}. This uniqueness condition is in terms of fixed points of the system mean-field limits -- these notions are defined later in the paper, which requires a fair amount of preliminaries and analysis. That is why the formal statement of Condition~\ref{cond-unique} is postponed until the point (Section~\ref{sec-cond-unique}) where we will be in position to do so.
}

\begin{thm}
\label{thm-main}

Consider an infinite-frame (non-truncated)  system, i.e with frame $[0,\infty)$. 
Suppose the additional Assumptions~\ref{cond-second-moment} and \ref{cond-ntriv-add} hold.
Suppose also that Condition~\ref{cond-unique} (given later in Section~\ref{sec-cond-unique}) holds for all $\lambda\ge 0$.

Then there exists $\bar \lambda \in (0,\infty)$ such that the following holds.

(i) For each $\lambda \in [0,\bar \lambda)$, there exists $\rho(\lambda) \in [0,1)$ and 
a unique proper element $x^{*,\lambda} \in \cx$, with $x^{*,\lambda}_0 = \rho(\lambda)$,
such that,
as $n\to\infty$, with $\lambda^n \to \lambda \ge 0$ and $\pi_j^n \to \pi_j$ for all $j$
(while all other system parameters remain fixed), 
\beql{eq-conv-main}
x^n(\infty) \Rightarrow x^{*,\lambda}.   
\eeql

(ii) Function $\rho(\lambda)$ is strictly increasing continuous, mapping $[0,\bar \lambda)$ into $[0,1)$. Function $x^{*,\lambda}$ (of $\lambda$) is strictly increasing continuous in $[0,\bar \lambda)$, and such that $x^{*,\lambda} \to x^{**}$ as $\lambda \uparrow \bar \lambda$.

(iii) If $\lambda \ge \bar \lambda$, 
\beql{eq-conv-main2}
x^n(\infty) \Rightarrow x^{**}, ~~n\to\infty.
\eeql

As corollaries of (ii) and (iii), we also have the following {\em steady-state asymptotic independence} properties,  
for any fixed integer $m \ge 1$. 

(ii') If $\lambda < \bar \lambda$,
\beql{eq-indep-main}
[W^n_1(\infty), \ldots, W^n_m(\infty)] \Rightarrow [W^*_1, \ldots, W^*_m], ~~ n\to\infty,
\eeql
where random variables $W^*_i$ are i.i.d., with
$\P\{W^*_i > w\} = x^{*,\lambda}_w, ~w\ge 0$.

(iii') If $\lambda \ge \bar \lambda$,
\beql{eq-indep-main2}
[W^n_1(\infty), \ldots, W^n_m(\infty)] \Rightarrow [\infty, \ldots, \infty], ~~ n\to\infty.
\eeql

(Property (ii') follows from (ii), because, due to the symmetry between servers,
in steady-state, the joint distribution of workloads of servers $1,\ldots,m$ is equal to that for a set of $m$ servers chosen 
uniformly at random. Property (iii') follows from (iii) by an analogous argument, which includes the possibility of pre-limit systems being either stable or unstable.)

\end{thm}

The proof of Theorem~\ref{thm-main} is in Section~\ref{sec-proof-main-th}.

\sasha{The next main result on the steady-state asymptotic independence is Theorem~\ref{thm-ihr} below.
It shows that the conclusions of Theorem~\ref{thm-main} hold for an important wide class of systems,
in which for each class $j$ the component sizes are i.i.d. with a distribution $H_j$ having increasing hazard rate
(IHR, see Definition~\ref{def-ihr}). Theorem~\ref{thm-ihr} is proved by showing that, in essence, such systems 
satisfy conditions of Theorem~\ref{thm-main}. Most importantly, i.i.d. IHR condition on component sizes implies Condition~\ref{cond-unique},
via an ``intermediate'' condition which we call $D$-monotonicity -- see Definition~\ref{def-d-monotone} and Theorem~\ref{thm-main-D-monotonicity} in Section~\ref{sec-cond-sufficient}.}

\begin{definition}[Increasing hazard rate (IHR)]
\label{def-ihr}
A distribution $H(\cdot)$ on $\R_+$ has increasing hazard rate (IHR), if
\beql{eq-ihr-cond}
\frac{H(y_1+\Delta) - H(y_1)}{1-H(y_1)} \le \frac{H(y_2+\Delta) - H(y_2)}{1-H(y_2)}, ~~\forall ~ 0 \le y_1 \le y_2 ~\mbox{such that $H(y_2)<1$}, ~\forall \Delta \ge 0.
\eeql
\end{definition}

An IHR distribution $H(\cdot)$ necessarily has the following structure. Let $y_{max} = \sup\{y \ge 0 ~|~ H(y) < 1\}
\le \infty$ (with $y_{max} = 0$ if $H(0)=1$).
Then there exists a non-negative non-decreasing hazard rate function $h(y), ~y\in [0,y_{max})$, such that:
$$
1-H(y) = \exp\{-\int_0^y h(u)du\}, ~~ y<y_{max}.
$$
Thus, $H(\cdot)$ has density $[1-H(y)] h(y)$ in $[0,y_{max})$; when $y_{max}<\infty$, the distribution may have an atom at $y_{max}$ with mass
$1-H(y_{max}-) =  \exp\{-\int_0^{y_{max}} h(u)du\}$.

Examples of IHR distributions include: the exponential distribution (which, in fact, has constant hazard rate); deterministic distribution 
concentrated on a single point $a \ge 0$; uniform distribution on a finite interval $[0,a]$. Also, if a random variable $A$ has an IHR distribution
and $a\ge 0$ is a constant, then the truncated random variable $A \wedge a$ also has an IHR distribution.

We will say that a distribution $H$ is a mixture of distributions $H_\chi$, parameterized by $\chi$,
if $H$ is obtained by averaging $H_\chi$ w.r.t. some probability measure 
on the values of $\chi$.

\begin{thm}
\label{thm-ihr}
Consider an infinite-frame (non-truncated)  system, i.e. with frame $[0,\infty)$. 
Suppose the distribution $F_j$ for each class $j$ is such that the component sizes are i.i.d. with an IHR distribution $H_j(\cdot)$ 
(in which case Assumption~\ref{cond-second-moment} holds automatically); 
or, more generally, $F_j$ may be a mixture of such distributions, 
\sasha{in which case we make Assumption~\ref{cond-second-moment}.} 
(In this theorem we do {\em not} a priori assume Condition~\ref{cond-unique} or make 
additional Assumption \ref{cond-ntriv-add}.)
Then the conclusions of Theorem~\ref{thm-main} hold.
\end{thm}

The proof of Theorem~\ref{thm-ihr} is in Section~\ref{sec-proof-ihr}.

\sasha{Some additional results on the steady-state asymptotic independence 
for infinite-frame systems are presented in Section~\ref{sec-further-conditions}.}

\subsection{Conditions sufficient for Condition~\ref{cond-unique} to hold.}
\label{sec-cond-sufficient}

For each job class class $j$, with its parameters $k_j, d_j, F_j$, and each integer $m \le d_j$, let us define a customer {\em sub-class} $j^m$ as follows: it is an (artificial) class with parameters $k_j^m, d_j^m, F_j^m$, where $k_j^m = k_j \wedge m$, $d_j^m = d_j \wedge m$, 
and the exchangeable $m$-dimensional distribution $F_j^m$ is the projection of distribution $F_j$ on first $m$ coordinates. In words, sub-class $j^m$ is the same as class $j$, except that only $m$ components out of $d_j$ are considered arrived, while the remaining $d_j-m$ components are ``ignored;'' the cancellation of the residual job service occurs when any $k_j \wedge m$ job components complete service. Note that when $m=d_j$ the sub-class $j^m$ is equal to the class $j$ itself; when $m=0$, $j^m$ is the ``empty sub-class'' -- nothing happens in the system when such job arrives.

Consider an infinite frame system. Consider a fixed job class $j$ (or sub-class $j^m$).
Let $(Z_1, \ldots, Z_{d_j})$ be the ordered workload vector for a selection set (of a class $j$ job); namely, $Z_i$ are the workloads of selected servers,
ordered to form a non-decreasing sequence, $Z_{i} \le Z_{i+1}$. Denote by $D=(D_1, \ldots, D_{d_j})$
the workload-differential vector, where $D_1 =0$, $D_i = Z_{i} - Z_{i-1}, ~i=2,\ldots,d_j$. 
Denote by $\eta_j(D)$ the (random) total workload brought by 
a class $j$ job, given the selection set workload-differential vector is $D$.
It is easy to see that 
\beql{eq-eta-cont}
\mbox{$\E \eta_j(D)$ is continuous in $D$.}
\eeql

\begin{definition}[$D$-monotonicity]
\label{def-d-monotone}
Consider an infinite frame system.
A class $j$ (or sub-class $j^m$) is called $D$-monotone decreasing, if its distribution $F_j$ of component sizes is such that $\E \eta_j(D)$ is non-increasing in $D$, i.e.
$D \ge D'$ implies $\E \eta_j(D) \le \E \eta_j(D')$.
\end{definition}

\begin{thm}
\label{thm-main-D-monotonicity}
(i) If all classes $j$ and all their sub-classes $j^m$ are D-monotone decreasing,
  then Condition~\ref{cond-unique} holds for all $\lambda \ge 0$.
  
(ii) Suppose the distribution $F_j$ for a class $j$ is such that the component sizes are i.i.d. with an IHR distribution $H_j(\cdot)$. 
(Or, more generally, $F_j$ may be a mixture of such distributions.)
Then class $j$ and all its sub-classes $j^m$ are
$D$-monotone decreasing.
\end{thm}

The proof of Theorem~\ref{thm-main-D-monotonicity} is in Section~\ref{sec-cond-unique}.

Further discussion of Condition~\ref{cond-unique} and related additional results on the steady-state asymptotic independence 
for infinite-frame systems are presented in Section~\ref{sec-further-conditions}.

\section{Monotonicity and continuity properties}
\label{sec-monotonicity-etc}

In this section we consider a system with fixed finite $n$, and
describe basic monotonicity and continuity properties of our model. For that, as in e.g. \cite{St2014_pull,SnSt2020},
it will be convenient to adopt a more general view of the model, which allows some of the server workloads to be infinite. Obviously, this generalization is only relevant for infinite-frame systems (because in finite-frame systems the workloads are uniformly upper bounded by definition). Specifically, if server $i$ workload is initially infinite, $W_i^n(0)=\infty$, then, by convention, it remains infinite at all times, $W_i^n(t)=\infty$, $t\ge 0$. The same workload placement rules apply, with the convention that an infinite workload remains infinite when ``more'' workload is added to it. Under this more general view of the model, $W^n(\cdot)$ 
is still a well defined Markov process,
but its state space is $\bar \R_+^n$, not $\R_+^n$. Similarly, $x^n(\cdot)$ is
still a well defined Markov process,
but its state space is a subset of $\cx$, not a subset of $\cx^{pr}$. 

Note that the fraction of infinite-workload servers remains constant at all times: $x^n_\infty(t)=x^n_\infty(0)$, $t\ge 0$.
This means that if the initial state $x^n(0)$ is proper, i.e. $x^n_\infty(0)=0$, 
the process ``lives'' on proper states thereafter. The version of the process confined to proper states (in $\cx^{pr}$),
is the process as we defined it originally, and we sometimes refer to it as {\em standard}. 
The version of the process allowing improper states in $\cx$ (that is, allowing infinite server workloads) will
be referred as {\em generalized}. In the rest of the paper, unless explicitly stated otherwise, we refer to the standard version
of the process. 

All continuity and monotonicity properties that we need stem from the following simple fact. (Its proof is very straightforward, so we skip it.)

\begin{lem}
\label{lem-monotone-cont-basic}
Consider a fixed set of servers, selected by a class-$j$ job, labeled WLOG by $1, 2, \ldots, d_j$, with (deterministic) workloads given by vector
$W=(W_1, \ldots, W_{d_j})$. Let $\xi^{(j)} = (\xi^{(j)}_1, \ldots, \xi^{(j)}_{d_j})$ be the (deterministic) vector of the component sizes of the job. Denote 
by $\hat W=(\hat W_1, \ldots, \hat W_{d_j})$ the vector of the server workloads after the job is added. Then $\hat W$ is monotone non-decreasing and continuous in both $W$ and $\xi^{(j)}$.
\end{lem}

As a corollary of Lemma~\ref{lem-monotone-cont-basic}, we obtain

\begin{lem}
\label{lem-monotone1}
For a given $n$ consider two versions of the generalized process. The first version $W^{n}(\cdot)$ has the frame size $c^n$,
the arrival rate parameter $\lambda^n$ and the job class probability distribution $\pi^n = (\pi^n_j)$, as defined in our model.
The second version $W^{n,(\ell)}(\cdot)$ has the frame size $c^{n,(\ell)}$,
the arrival rate parameter $\lambda^{n,(\ell)}$, and a possibly different job class probability distribution $\pi^{n,(\ell)} = (\pi^{n,(\ell)}_j)$.
Suppose $c^n \le c^{n,(\ell)}$ and $\lambda^n \pi^n \le \lambda^{n,(\ell)} \pi^{n,(\ell)}$.
Suppose, $W^{n}(0) \le W^{n,(\ell)}(0)$. Then processes $W^{n}(\cdot)$ and $W^{n,(\ell)}(\cdot)$ can be coupled so that, w.p.1,
\beql{eq-key-monotone}
W^{n}(t) \le W^{n,(\ell)}(t), ~\forall t\ge 0.
\eeql
\end{lem}

{\em Proof.} We couple the job arrival processes in the natural way, so that when there is a job arrival in the former system there is also the same class job arrival in the latter system, with the same selection set and same component sizes. 
 It is easy to see that if \eqn{eq-key-monotone} holds just before (any) job arrival, it holds right after the job arrival as well, and also 
 in the time interval until the next arrival. By induction on the job arrival times, we see that \eqn{eq-key-monotone} prevails at all times.
 $\Box$

Another corollary of Lemma~\ref{lem-monotone-cont-basic} is the following Lemma~\ref{lem-monotone2}, which, in essence, a more general version of Lemma~\ref{lem-monotone1}. It is formulated in terms of the processes' sample path, for convenience.
\sasha{In Lemma~\ref{lem-monotone2} we take a still more general view of the system. Namely, we allow that the ``left regulation boundary (point),'' which particles cannot cross, is not necessarily $0$ and not necessarily fixed. Instead, this regulation point,
denoted by $B^{n}(t)$ is a function of time $t\ge 0$.}

\begin{lem}
\label{lem-monotone2}
For a given $n$ consider two versions of the generalized process realizations, $W^{n}(\cdot)$ and $W^{n,(\ell)}(\cdot)$,
corresponding to two systems with infinite frame sizes. The realizations are such that
both systems receive the same sequence of common job arrivals (same job arrival times, job classes, selections sets, and component sizes).
In addition, the realization of $W^{n,(\ell)}(\cdot)$ is such that additional amounts of workload may be added to any servers at any times.
The initial states (particle locations in the two systems) are such that $W^{n}(0) \le W^{n,(\ell)}(0)$, and the initial regulation boundaries 
(points)
$B^{n}(0)$ and $B^{n,(\ell)}(0)$ are such that
$B^{n}(0) \le B^{n,(\ell)}(0)$. Assume that one or both regulation points may move left in a way such that their instantaneous speeds are either $-1$ or $0$;
however, $B^{n}(t) \le B^{n,(\ell)}(t)$ holds at all times. Then, $W^{n}(t) \le W^{n,(\ell)}(t)$ holds for all $t\ge 0$.
\end{lem}

 Yet another corollary of Lemma~\ref{lem-monotone-cont-basic} is the following
 
\begin{lem}
\label{lem-cont1}
For a given $n$ consider a generalized process $W^{n}(\cdot)$ 
corresponding to the system with frame size $c^n$ and arrival rate $\lambda^n$.
Consider also a {\em sequence} of generalized processes $W^{n,(\ell)}(\cdot)$, indexed by $\ell \to \infty$, with the
frame sizes $c^{n,(\ell)}$ and arrival rates $\lambda^{n,(\ell)}$.
(For all these processes the probabilities $\pi_j^n$ are same.)
Suppose, as $\ell \to \infty$, $W^{n,(\ell)}(0) \to W^{n}(0)$, $c^{n,(\ell)} \to c^n$ and arrival rates $\lambda^{n,(\ell)} \to \lambda^n$.
Then these processes can be coupled so that, w.p.1,
\beql{eq-key-cont}
W^{n,(\ell)}(t) \to W^{n}(t), ~\forall t \ge 0, ~~\mbox{and}~~ W^{n,(\ell)}(\cdot) \stackrel{J_1}{\longrightarrow} W^{n}(\cdot).
\eeql
\end{lem}

 {\em Proof.} We use natural coupling of process $W^{n}(\cdot)$ and all processes $W^{n,(\ell)}(\cdot)$
 so that, w.p.1, on any finite time interval $[0,T]$, for all sufficiently large $\ell$,
 the number of job arrivals, the job arrival times and classes, as well as component sizes, of the $W^{n,(\ell)}(\cdot)$ process
 are equal to those of the $W^{n}(\cdot)$ process. Then, w.p.1,  
 $W^{n,(\ell)}(t) \to W^{n}(t)$ for any $t \ge 0$
 and, moreover, on any finite time interval $[0,T]$, for all sufficiently large $\ell$, trajectory $W^{n,(\ell)}(\cdot)$
 has jumps at the same times as $W^{n,(\ell)}(\cdot)$. This implies \eqn{eq-key-cont}.
 $\Box$
 
It follows from Lemma~\ref{lem-monotone1} that the following properties hold for a fixed $n$. 
For any $\lambda^n$ the generalized process $W^{n,\emptyset}(\cdot)$ starting 
from the empty initial state is stochastically increasing in time, i.e.
\beql{eq-monotine-from-empty}
W^{n,\emptyset}(t_1) \le_{st} W^{n,\emptyset}(t_2), ~~t_1 \le t_2.
\eeql
\sasha{Indeed, $W^{n,\emptyset}(0) = (0,\ldots,0) \le_{st} W^{n,\emptyset}(t_2-t_1)$, and \eqn{eq-monotine-from-empty}
follows by applying Lemma~\ref{lem-monotone1} with $t=t_1$. The monotonicity property \eqn{eq-monotine-from-empty}
will be used extensively in our analysis.} 
In particular, it implies that 
$$
W^{n,\emptyset}(t) \Rightarrow W^{n}(\infty), ~~t\to \infty,
$$
where the distribution of the random element $W^{n}(\infty)$ in $\bar \R_+^n$ is the {\em lower invariant measure} of the generalized process 
$W^{n}(\cdot)$, i.e. it is the (uniquely defined) stationary distribution of $W^{n}(\cdot)$ which is stochastically dominated by any other stationary distribution. Note that, as already discussed earlier, the process $W^{n}(\cdot)$ for a finite-frame system is trivially stable,
and therefore the lower invariant measure is its unique stationary distribution. For the infinite-frame system, clearly, (standard) process 
$W^{n}(\cdot)$ is stable if and only if the lower invariant measure is proper, meaning that $W^n_i(\infty)$ is finite w.p.1. for all $i$;
if so, the lower invariant measure is the unique stationary distribution of the standard process. 
Also note that,
formally speaking, whether the standard process $W^{n}(\cdot)$ is
stable or not, the distribution concentrated on the point $(\infty, \ldots, \infty)$
is always a stationary distribution of the generalized process.
Finally, note that, clearly, if the (standard) process $W^{n}(\cdot)$ is stable for arrival rate $\lambda^n$, it is also stable for any smaller arrival rate.

\begin{lem}
\label{lem-stabil-solidarity}
(i) The (standard) process $W^{n}(\cdot)$ [resp. $x^n(\cdot)$] is stable if and only if 
the lower invariant measure is proper.

(ii) The (standard) process $W^{n}(\cdot)$ [resp. $x^n(\cdot)$] is unstable if and only if 
the lower invariant measure is concentrated on
the point at infinity, $(\infty, \ldots, \infty)$. 
In other words, the instability is equivalent to the property that
for any initial state $W^{n}(0)$ we have
\beql{eq-unstable}
W^{n}(t) \Rightarrow (\infty, \ldots, \infty), ~t\to\infty.
\eeql
\end{lem}

{\em Proof.} (i) This is already observed in the discussion preceding the lemma.

(ii) The sufficiency is obvious -- \eqn{eq-unstable} directly contradicts the existence of a proper stationary distribution.
Let us prove the necessity. 
Suppose, the process $W^{n}(\cdot)$ is unstable. We need to prove \eqn{eq-unstable} for the process
starting from the empty initial state $W^n(0)$.
Recall that $W^{n}(t)$ is stochastically monotone increasing in $t$, converging in distribution to $W^{n}(\infty)$, whose
distribution is the lower invariant measure. By instability, the lower invariant measure is {\em not} proper, that is 
$\P \{ W_i^{n}(\infty) =\infty\} =\delta >0$ for some $i$. In fact, by symmetry the
constant $\delta \in (0,1]$ is same for all $i$. Proceeding exactly as in the proof of Lemma 4 in \cite{St2014_pull}, we show that this
$\delta$ must satisfy $\delta(1-\delta) \le 0$, and therefore $\delta=1$. Thus, $W^{n}(\infty) = (\infty,\ldots,\infty)$.
$\Box$

Lemma~\ref{lem-stabil-solidarity} justifies our earlier convention in Section~\ref{sec-asymp-regime} by which 
$W^{n}(\infty) = (\infty,\ldots,\infty)$
for an unstable process. 

From Lemmas~\ref{lem-monotone1}, \ref{lem-cont1} and \ref{lem-stabil-solidarity} we obtain the following

\begin{lem}
\label{lem-stationary-monotone}
Recall the definition of $\bar \lambda^n$ in \eqn{eq-crit-load-finite-def}.
For a fixed $n$ and $\lambda^n < \bar \lambda^n$, the (standard) process $W^{n, \lambda^n}(\cdot)$ is stable;
denote by $W^{n, \lambda^n}(\infty)$ a random element whose distribution is the stationary distribution 
(equal to the lower invariant measure) of the process.
Then $W^{n, \lambda^n}(\infty)$ is continuous (in the sense of convergence in distribution) and stochastically
non-decreasing for $\lambda^n \in [0,\bar \lambda^n)$.
Furthermore, $\rho^n(\lambda^n) \doteq \P\{W_1^{n, \lambda^n}(\infty) > 0\}$ is a non-decreasing
continuous function of $\lambda^n \in [0,\bar \lambda^n)$.
\end{lem}

{\em Proof.} The fact that $W^{n, \lambda^n}(\infty)$ is stochastically non-decreasing follows from the process 
monotonicity in $\lambda^n$, and the distribution of $W^{n, \lambda^n}(\infty)$ being the lower invariant measure.
If we have a sequence $\lambda^{n,\ell} \to \lambda^n \in [0,\bar \lambda^n)$, the corresponding processes can be
coupled so that w.p.1 we have the convergence of their busy periods and idle periods;
this implies the convergence $W^{n, \lambda^{n,\ell}}(\infty) \Rightarrow W^{n, \lambda^n}(\infty)$.
$\Box$

Finally, we will need the following fact, which easily follows from Lemma~\ref{lem-monotone1} and 
the stochastic monotonicity \eqn{eq-monotine-from-empty} of the process starting from the empty state.

\begin{lem}
\label{lem-load-monotone}
For a fixed $n$, denote by  $W^{n, \lambda^n, \emptyset}(t), ~t\ge 0,$ the process with parameter $\lambda^n$, starting from the 
empty state. Denote by $\rho^n(\lambda^n,t) \doteq \P\{W_1^{n, \lambda^n, \emptyset}(t) > 0\}$ the system expected load at time $t$.
Then, $\rho^n(\lambda^n,t)$ is continuous, strictly increasing in both $\lambda^n$ and $t$. We also have: $\rho^n(\lambda^n,0)=0$
for any $\lambda^n$; $\rho^n(0,t)=0$ and $\rho^n(\lambda^n,t) \uparrow 1$ as $\lambda^n \to\infty$, for any $t$.
If the process is stable for $\lambda^n$, then $\rho^n(\lambda^n,t) \uparrow \rho^n(\lambda^n)$ as $t\to\infty$.
\end{lem}

\section{An equivalent view of a generalized process}
\label{sec-equiv-view}

In this section we
observe that a generalized process (for a 
system with some servers having infinite workload) can be equivalently viewed as a standard process for 
the corresponding ``reduced'' system, which includes only those servers with finite workloads and in which the set of job classes is expanded to include all sub-classes of each class $j$.
We will use this equivalent view in the proofs of our main results. 

Suppose our system with index $n$ is such that initially (and then at any time) $\phi^n n$ servers have infinite workloads; the remaining 
$(1-\phi^n) n$ servers -- let us call them ``standard'' -- have initially (and then at any finite time) finite workloads. Clearly, {\em as far as behavior of the standard servers is concerned,} this system is equivalent to the following one, which we will call a $\phi^n$-reduced system. The system contains $(1-\phi^n) n$ servers, all of which are standard; 
the total job arrival rate is $\lambda^n n$, same as in the original system, which means the rate per server is 
$\lambda^n/(1-\phi^n)$; the set of job classes includes all subclasses $j^m$ of all original
classes $j$; the probability that an arriving job is of sub-class $j^m$ is $\pi_{j^m}^n = \pi_j^n \gamma_{j^m}^n$, where
$\gamma_{j^m}^n$ is the probability that if $d_j$ servers are selected uniformly at random from $n$ servers,
exactly $m$ of them will be within the subset of $(1-\phi^n) n$ regular servers. Clearly, the process describing the behavior of the reduced system is standard.

Now, suppose the asymptotic regime that we consider is such that $\phi^n \to \phi \in [0,1)$. Then, for each sub-class $j^m$,
we have
$$
\pi_{j^m}^n \to \pi_{j^m} \doteq \pi_j \frac{d_j !}{m! (d_j-m)!} (1-\phi)^m \phi^{d_j-m}.
$$
Clearly, {\em as far as behavior of the standard servers is concerned,} 
considering such sequence of systems 
is equivalent to considering the corresponding sequence of $\phi^n$-reduced systems. The sequence of reduced systems fits into the same asymptotic regime as our original sequence of systems, but has: different number of servers ($(1-\phi^n) n$ instead of $n$); 
different per server job arrival rate $\lambda^n/(1-\phi^n) \to \lambda/(1-\phi)$; different set of classes (which includes all sub-classes $j^m$ of all classes $j$); and different probabilities $\pi_{j^m}^n$ of class occurrence, 
such that $\pi_{j^m}^n \to \pi_{j^m}$. 

Note that if all sub-classes of the original system are $D$-monotone, the same is obviously true for any reduced system.

Finally, when we consider a reduced system, its frame size $c$ may be finite as well as infinite. Recall that, if frame size $c<\infty$, the workloads of (standard) servers are truncated at level $c$. 

\section{Proof of Theorem~\ref{thm-crit-load-finite}}
\label{sec-proof444}

\sasha{This proof relies on the basic monotonicity properties (in Section~\ref{sec-monotonicity-etc}) and on the fluid limit technique
for establishing stability. (The proof does {\em not} use the material in Section~\ref{sec-equiv-view}.)}

Recall that we consider a system with a fixed $n$ and the process $W^{n,\lambda^n} (\cdot)$, with the superscript 
$\lambda^n$ indicating the dependence on the arrival rate. Also recall the notation $W^{n,\lambda^n,\emptyset}(\cdot)$ 
for the process starting specifically from the empty initial state $W^{n,\lambda^n,\emptyset} (0) = (0,\ldots,0)$.

(i) The uniform upper and lower bounds on $\bar \lambda^n$
 follow from the fact that the expected 
total amount of new workload brought by an arriving job is both upper and lower bounded uniformly in $n$ and uniformly in system states at the arrival time. The stability for any $\lambda^n < \bar \lambda^n$ follows from the basic monotonicity properties in 
Section~\ref{sec-monotonicity-etc}; specifically, the stability for a given $\lambda^n$ implies stability for any smaller arrival rate as well.

(ii) Recall notation $\rho^n(\lambda^n) = \pr\{W_1^{n,\lambda^n}(\infty) > 0\}$. 
By Lemma~\ref{lem-stationary-monotone},
as $\lambda^n$ increases in the interval $[0,\bar \lambda^n)$,  $\rho^n(\lambda^n)$ strictly increases, continuously, 
from $0$ to some value $\nu \le 1$. It remains to prove that $\nu = 1$. The proof is by contradiction. Suppose, $\nu < 1$. 

As the first step, we prove that (if $\nu<1$ were to hold) the process $W^{n,\bar \lambda^n}(\cdot)$ is stable, and 
\beql{eq-load-nu}
\rho^n(\bar \lambda^n)=\nu.
\eeql
Consider a sequence of times $t_\ell \to \infty$. For each $t_\ell$, by Lemma~\ref{lem-load-monotone},
there exists unique $\beta^{(\ell)}$
such that $\pr\{W_1^{n,\beta^{(\ell)},\emptyset}(t_\ell) > 0\}=\nu$. Again by Lemma~\ref{lem-load-monotone}, $\beta^{(\ell)}$ is non-increasing
in $t_\ell$.
We must have $\lim_\ell \beta^{(\ell)} \ge \bar \lambda^n$. Indeed, if some $\beta^{(\ell)} < \bar \lambda^n$,
then for a small $\epsilon>0$, $\pr\{W_1^{n,\bar \lambda^n-\epsilon,\emptyset}(t_\ell)> 0\}  > \pr\{W_1^{n,\beta^{(\ell)},\emptyset}(t_\ell) > 0\} = \nu$;
this is, however, not possible because $W^{n,\bar \lambda^n-\epsilon}(\cdot)$ is stable, and then
$\pr\{W_1^{n,\bar \lambda^n-\epsilon,\emptyset}(t_\ell)> 0\}  \le \pr\{W_1^{n,\bar \lambda^n-\epsilon}(\infty)> 0\} < \nu$. Now, since 
$\beta^{(\ell)} \ge \bar \lambda^n$ for all $t_\ell$, we have 
$\pr\{W_1^{n,\bar \lambda^n,\emptyset}(t_\ell)> 0\} \le \pr\{W_1^{n,\beta^{(\ell)},\emptyset}(t_\ell)> 0\} = \nu$ for all $t_\ell$.
Recall that $t_\ell \to \infty$, which implies that $W^{n,\bar \lambda^n}(\cdot)$ is stable and
$\rho^n(\bar \lambda^n) \le \nu$, where the strict inequality is impossible by the definition of $\bar \lambda^n$.
Thus we obtain \eqn{eq-load-nu}.

\sasha{The second step is to show that \eqn{eq-load-nu}, with $\nu<1$, leads to a contradiction. 
This is the key step. The basic intuition for it is as follows. If the system with $\lambda^n$ exactly equal to $\bar \lambda^n$ 
is stable, with the load equal to $\nu<1$, the process ``hits'' the empty state (when all workloads are at the regulation boundary at $0$) at a positive rate. Using monotonicity arguments, this implies that maximum of the workloads has a strictly negative drift. (More precisely, the process fluid limit is such that the maximum of workloads, when away from $0$, has negative, bounded away from $0$ derivative.) Then we observe that a system with $\lambda^n=\mu$ being slightly larger than $\bar \lambda^n$, 
is equivalent, {\em up time rescaling}, to the system with $\bar \lambda^n = \bar \lambda^n$ and the speed at which workloads move to the left being slightly larger (i.e., smaller in absolute value) than $-1$. (Note that a time rescaling does not affect
a system stability/instability.) Using this, we can show that, when $\lambda^n=\mu$ is sufficiently close to $\bar \lambda^n$,
the maximum of workloads in the equivalent system also must have a strictly negative drift. This would apply the stability
of the equivalent system, and then also the stability of the system with $\lambda^n=\mu > \bar \lambda^n$ -- a contradiction with the definition of $\bar \lambda^n$.}

We now proceed with the formal argument leading to a contradiction. Suppose, \eqn{eq-load-nu} holds with $\nu<1$.
We start with constructing a process that dominates $W^{n,\bar \lambda^n}(\cdot)$, and easily seen to be also stable.
Consider the process, let us call it $\hat W^{n,\bar \lambda^n}(\cdot)$, such that, initially, all $\hat W_i^{n,\bar \lambda^n}(0)$ 
are set to be equal to $\max_i W_i^{n,\bar \lambda^n}(0)$,
and the left regulation boundary is reset from $0$ to that maximum.
Using Lemma~\ref{lem-monotone2}, the proceses $\hat W^{n,\bar \lambda^n}(\cdot)$ and $W^{n,\bar \lambda^n}(\cdot)$
can be coupled so that, w.p.1, $W^{n,\bar \lambda^n}(t) \le \hat W^{n,\bar \lambda^n}(t)$ for all $t$; thus 
$\hat W^{n,\bar \lambda^n}(\cdot)$ is a stochastic upper bound of $W^{n,\bar \lambda^n}(\cdot)$.
Since $\hat W^{n,\bar \lambda^n}(\cdot)$ is just a (shifted to the right) version of the stable process 
$W^{n,\bar \lambda^n}(\cdot)$,
it is also stable, and then the steady-state probability of all workloads being at the regulation boundary is some $\delta>0$. 
Now consider yet another version of process $W^{n,\bar \lambda^n}(\cdot)$, let us call it 
$\tilde W^{n,\bar \lambda^n}(\cdot)$. This version is same as $\hat W^{n,\bar \lambda^n}(\cdot)$, with the same 
regulation, except the following: in the time intervals when {\em all} workloads 
are at the regulation boundary, the regulation boundary itself moves 
``left''  at the constant rate $-1$ until and unless it hits $0$ -- starting that time the regulation boundary stays at $0$.
By Lemma~\ref{lem-monotone2}, the process $\tilde W^{n,\bar \lambda^n}(\cdot)$ is also a stochastic upper bound of $W^{n,\bar \lambda^n}(\cdot)$.  Comparing $\tilde W^{n,\bar \lambda^n}(\cdot)$ to $\hat W^{n,\bar \lambda^n}(\cdot)$, we observe that the regulation boundary of the former 
has the average drift $-\delta <0$, until the boundary hits $0$. Indeed, if the points in time when all workloads are at the regulation boundary are the process ``renewals,'' the process $\tilde W^{n,\bar \lambda^n}(\cdot)$ has the ``same'' renewal cycles as $\hat W^{n,\bar \lambda^n}(\cdot)$, except while it stays in the renewal state and the regulation boundary is away from $0$, the regulation boundary itself -- along with all workloads -- moves left at rate $-1$. Note that the expected duration of a renewal cycle is finite (by 
the stability of $\hat W^{n,\bar \lambda^n}(\cdot)$). Thus, he average drift of 
of the regulation boundary of 
$\tilde W^{n,\bar \lambda^n}(\cdot)$ is $-\delta$.

Consider now the version of process $W^{n,\mu}(\cdot)$, with the arrival rate parameter $\mu > \bar \lambda^n + \delta_1$.
Specifically consider $\mu$ being close to $\bar \lambda^n$ -- we will specify later how close.
Consider the process $W^{n,\mu}(\cdot)$, slowed down by the factor $\gamma = \bar \lambda^n/\mu < 1$, namely the
process
$$
W^{n,\mu, (slow)}(t) = W^{n,\bar \lambda^n}(\gamma t), ~t\ge 0.
$$
Obviously, $W^{n,\mu, (slow)}(\cdot)$ is stable if and only if $W^{n,\mu}(\cdot)$ is.
The process $W^{n,\mu, (slow)}(\cdot)$ is the same as $W^{n,\bar \lambda^n}(\cdot)$, with the same job arrival rates
and the same regulation boundary at $0$, but with workloads decreasing at rate $-\gamma$ (instead of $-1$) between the job arrivals.
Consider the following process $\tilde W^{n,\mu, (slow)}(\cdot)$. (Roughly speaking, 
it has the same relation to $W^{n,\mu, (slow)}(\cdot)$ as 
 $\tilde W^{n,\bar \lambda^n}(\cdot)$ does to $W^{n,\bar \lambda^n}(\cdot)$.)
 Initially, $\tilde W_i^{n,\mu, (slow)}(0)= \max_s W_s^{n,\mu, (slow)}(0)$, and the regulation boundary is at $0$.
 Unless all workloads are at the regulation boundary, the boundary it moves {\em right} at the constant positive speed $1-\gamma$;
 note that {\em with respect to the boundary} the workloads that are not on it, move at the rate $-1$.
 When, and only when, all workloads are at the regulation boundary, the boundary (along with all workload) moves left at 
 speed $-\gamma$; the regulation boundary can never go to the left of $0$. 
 Process $\tilde W^{n,\mu, (slow)}(\cdot)$ is a stochastic upper bound of $W^{n,\mu, (slow)}(\cdot)$.
 Note that the renewal cycles of
 $\tilde W^{n,\mu, (slow)}(\cdot)$ ``with respect to the regulation boundary,'' again, have the same structure as those of $\hat W^{n,\bar \lambda^n}(\cdot)$; in particular, the steady-state probability of all workloads being at the boundary is $\delta$. 
What is different is the movement of the regulation boundary itself. Combining all these observations,
we conclude that the regulation boundary of the
process $\tilde W^{n,\mu, (slow)}(\cdot)$ 
 has the average drift
$$
-\delta_1=\delta (-\gamma) + (1-\delta) (1-\gamma) = -\delta + (1-\gamma),
$$
as long as the boundary is away from $0$. We now choose $\mu$ sufficiently close to $\bar \lambda^n$, 
i.e. $\gamma$ sufficiently close to $1$, so that $-\delta_1<0$.
This, along with the finite expectation of a renewal cycle duration,
implies the stability of $\tilde W^{n,\mu, (slow)}(\cdot)$. For example,   
the {\em fluid limit technique} for stability \cite{RS92, Dai95, St95, Bramson-book}
can be applied; indeed, in this case
the fluid limit 
$q(\cdot)=(q_i(\cdot))$ of the process $\tilde W^{n,\mu, (slow)}(\cdot)$ is easily seen to be unique and such that 
all $q_i(t)$ are equal at all times and $q'_i(t) = - \delta_1 <0$ when $q_i(t)>0$;
this implies the stability of the fluid limit, and then the stability of $\tilde W^{n,\mu, (slow)}(\cdot)$.
 Then, by the stochastic dominance,  $W^{n,\mu, (slow)}(\cdot)$
is also stable. Finally, this implies the stability of $W^{n,\mu}(\cdot)$, which, since $\mu > \bar \lambda^n$, contradicts
the definition of $\bar \lambda^n$.
 This contradiction, obtained from the assumption $\nu<1$, proves that $\nu=1$, which completes the proof.
$\Box$

\section{Proof of Theorem~\ref{th-closeness}}
\label{sec-closeness-proof}

\sasha{This proof is mostly self-contained. It does use the basic continuity properties in Section~\ref{sec-monotonicity-etc}. (The proof does {\em not} use the material in Sections~\ref{sec-equiv-view} and \ref{sec-proof444}.) The basic intuition for this proof is given in the discussion of Theorem~\ref{th-closeness} following its statement in Section~\ref{sec-tightness-centered}.}

If $x^n(t)$ is the state at $t$, it can be equivalently described as
$p^n(t)=(w_1(t), w_2(t), \ldots, w_n(t); z(t))$, where $z(t) = - \bar x^n(t) \le 0$ is the location of the origin w.r.t. to the mean $\bar x^n(t)$,
and 
$w_1, w_2, \ldots, w_n \ge z$ are the locations of the $n$ particles (not necessarily ordered in any way)
w.r.t. the mean $\bar x^n(t)$. 
Slightly abusing notation by defining $\Phi_\ell(p^n(t))$ as $\Phi_\ell(\mathring x^n(t))$ for the corresponding
$\mathring x^n(t)$, that is
$$
\Phi_\ell(p^n(t)) = \frac{1}{n} \sum_{i=1}^n |w_i(t)|^\ell.
$$
The evolution of $p^n(t)$ is as follows. Note that here the particle locations $w_i(t)$ may be negative as well as non-negative.
Point $z(t)$ serves as a ``regulation boundary,'' which evolves in time; $z(t)\le 0$ at all times.
Denote by $I_z(t)$ the (possibly empty) subset of particles, 
which are located exactly at the ``regulation boundary''  $z(t)$ at time $t$.
Between the times of the particle jumps (job arrivals), 
the  boundary $z(t)$ moves ``right'' at the constant non-negative speed $(n-|I_z(t)|)/n$; 
the particles in $I_z(t)$ (those located exactly at the boundary $z(t)$) 
stay at  the boundary and therefore move with it right at constant non-negative speed $(n-|I_z(t)|)/n$; 
each particle that is not at $z(t)$ (i.e. not within subset $I_z(t)$) moves ``left'' at the constant non-positive speed $-|I_z(t)|/n$;
when a particle ``hits'' boundary $z(t)$, it joins the set $I_z(t)$. (The average of the particle locations stays at $0$, as it should by
the $p(t)$ definition.)
It is easy to see that, for any $\ell\ge 1$, $\Phi_\ell(p^n(t))$ is absolutely continuous with the derivative (existing almost everywhere) 
$(d/dt)\Phi_\ell(p^n(t)) \le 0$, and moreover $(d/dt)\Phi_\ell(p(t)) < 0$ as long as  $0< |I_z(t)| < n$. At a time when one or more particles jump (upon a job arrival), the following occurs. Let $\kappa_i(t)$ be the amount of new workload added to server $i$ by the
job arriving at time $t$ (which may be non-zero only for the servers within the selection set of the job). Then, particle $i$ jump size at $t$ 
is $\zeta_i=\kappa_i(t) - \sum_s \kappa_s(t)/n$ (which may be positive or negative); the point $z(t)$ jumps (left) by 
$ - \sum_s \kappa_s(t)/n$.

Consider the following function $G^{(n)}(\mathring x^n), ~\mathring x^n \in \mathring\cx^{(n)},$ 
$$
G^{(n)} (\mathring x^n) = G^{(n)} (p^n)
= 2 \lambda^n n \frac{1}{n} \sum_i w_i \E  \zeta_i^{(\mathring x^n)},
$$
where $\zeta_i^{(\mathring x^n)}$ is the random jump size (which can have any sign) of particle $i$ upon a job arrival when the state is 
$\mathring x^n$. (The sizes $\zeta_i^{(\mathring x^n)}$ are dependent across $i$, of course.) 
As we will see, function $G^{(n)} (p^n)$ can be thought of as the ``first-order approximation of the generator of process $p^n(\cdot)$, applied to function $\Phi_2(p^n)$;'' but we do not even claim that $\Phi_2(p^n)$ is within the generator domain.

For future reference, note that, for each $n$, $G^{(n)} (\mathring x^n)$ is continuous in $\mathring x^n \in \mathring\cx^{(n)}$.
Also note that each 
$\E  \zeta_i^{(\mathring x^n)}$ is the quantity of the order $O(1/n)$, which motivates the definition
\beql{eq-bar-zeta-def}
\bar \zeta_i^{(\mathring x^n)} \doteq n \E \zeta_i^{(\mathring x^n)}.
\eeql
We can then write
$$
G^{(n)} (\mathring x^n) = G^{(n)} (p)
= 2 \lambda^n \frac{1}{n} \sum_i w_i \bar \zeta_i^{(\mathring x^n)}.
$$
We define the function $\bar\zeta^{(\mathring x^n)}(w), ~w\in \R,$ as follows: 
$\bar \zeta^{(\mathring x^n)}(w) = \bar \zeta^{(\mathring x^n)}_i$, where $i$ is the particle whose location $w_i$ is the closest to $w$ on the left; we also adopt a conventions that, if $w_i$ is the location 
of the left-most particle, then $\bar\zeta^{(\mathring x^n)}(w) = \bar \zeta_i^{(\mathring x^n)}$ for all $w< w_i$.
Clearly, function $\bar\zeta^{(\mathring x^n)}(w)$ is a piece-wise constant non-increasing function, and we can write
\beql{eq-gen1}
G^{(n)} (\mathring x^n) = G^{(n)} (p^n)
= 2 \lambda^n \int_{-\infty}^{\infty} \bar\zeta^{(\mathring x^n)}(w) w d[-\mathring x^n_w] \le 0,
\eeql
where the inequality in \eqn{eq-gen1} follows because
$$
\int_{-\infty}^{\infty} \bar\zeta^{(\mathring x^n)}(w) w d[-\mathring x^n_w] 
= \int_{-\infty}^{\infty} [\bar\zeta^{(\mathring x^n)}(w) - \bar\zeta^{(\mathring x^n)}(0)] w d[-\mathring x^n_w] 
+ \bar\zeta^{(\mathring x^n)}(0) \int_{-\infty}^{\infty}  w d[-\mathring x^n_w],
$$
$\bar\zeta^{(\mathring x^n)}(w) - \bar\zeta^{(\mathring x^n)}(0)$ is non-increasing with value $0$ at $0$,
and $\mathring x^n \in \mathring{\cx}_1$
(and then $\int_{-\infty}^{\infty} w d[-\mathring x^n_w] = 0$).

Next, we claim the following property: there exists a sufficiently large $C>0$ and some $\epsilon > 0$, such that,
uniformly in all sufficiently large $n$ and all $\mathring x^n \in \mathring\cx^{(n)}$ with $\Phi_1(\mathring x^n) \ge C$,
\beql{eq-drift}
G^{(n)} (\mathring x^n) = G^{(n)} (p^n) \le -\epsilon \Phi_1(\mathring x^n).
\eeql
The proof of \eqn{eq-drift} is given in Section~\ref{sec-eq-drift}.

From \eqn{eq-drift} and \eqn{eq-gen1} we obtain that, uniformly in all sufficiently large $n$,
\beql{eq-drift2}
G^{(n)}(\mathring x^n) = G^{(n)}(p^n) \le -\epsilon \Phi_1(\mathring x^n) +\epsilon C.
\eeql

Denote by $\Phi_2^{(C_2)}(p^n) = \Phi_2(p^n) \wedge C_2$ the function $\Phi_2$ truncated at level $C_2$.
Given that this is a continuous bounded function, and using (Lemma~\ref{lem-monotone-cont-basic}) the continuity  in $p^n$
of the workload amounts added by a job, 
it is not hard to see that it is within the domain 
of the generator $A^{(n)}$ of process $p^n(\cdot)$. 

Next, we claim the following fact:
there exists $C_1>0$ such that for any fixed $C_2>0$, 
uniformly in all large $n$ and $p^n$ such that $\Phi_2(p^n) \le C_2$, we have
\beql{eq-generator-bound}
A^{(n)} \Phi_2^{(C_2)}(p^n) \le G^{(n)}(p^n) + C_1,
\eeql
and then
\beql{eq-generator-bound5}
A^{(n)} \Phi_2^{(C_2)}(p^n) \le  -\epsilon \Phi_1(p^n) + C_3,
\eeql
with $C_3 = \epsilon C + C_1$.
The proof of \eqn{eq-generator-bound} is given in Section~\ref{sec-generator-bound}.

Bound \eqn{eq-generator-bound5} in turn implies that
 for any fixed $C_2>0$, 
\beql{eq-generator-bound2}
A^{(n)} \Phi_2^{(C_2)}(p^n)  \le [-\epsilon \Phi_1(p^n) + C_3] \bI \{\Phi_2(p^n) \le C_2\}.
\eeql

Recalling that $p^n(\infty)$ is the random value of $p^n(t)$ in the stationary regime, 
we have for all large $n$:
$$
0 = \E A^{(n)}  \Phi_2^{(C_2)}(p^n(\infty)) \le 
\E \left[ (-\epsilon \Phi_1(\mathring x^n(\infty)) + C_3) \bI \{\Phi_2(\mathring x^n(\infty)) \le C_2\} \right]
$$
$$
\le -\epsilon \E [\Phi_1(\mathring x^n(\infty))  \bI \{\Phi_2(\mathring x^n(\infty)) \le C_2\}] + C_3,
$$
and then
$$
\E [\Phi_1(\mathring x^n(\infty))  \bI \{\Phi_2(\mathring x^n(\infty))  \le C_2\}] \le C_3/\epsilon.
$$
Letting $C_2 \uparrow \infty$, we finally obtain that
$$
\E \Phi_1(\mathring x^n(\infty)) \le C_3/\epsilon
$$
for all sufficiently large $n$, and then 
$$
\E \Phi_1(\mathring x^n(\infty)) \le \bar C
$$
holds for all $n$ for some large $\bar C>0$.
$\Box$

\subsection{Proof of \eqn{eq-drift}.}
\label{sec-eq-drift}

The definition of $\bar \zeta_i=\bar \zeta_i^{(\mathring x^n)}$ in \eqn{eq-bar-zeta-def} can be interpreted as follows:  
$\bar \zeta_i$ is the expected amount of new workload 
$\E [\kappa_i | i\in S]$ the server $i$ receives, {\em conditioned on it being selected by a job}, multiplied by constant $[n\P\{i\in S\}]$
(which is close to the constant $\sum_j \pi_j d_j$ for all large $n$), and then centered by the total expected amount of new 
workload $\E \sum_s \kappa_s$ brought by a job. 

The proof is by contradiction. Suppose property \eqn{eq-drift} does not hold. Then, we can and do choose a subsequence of $n\to\infty$, and corresponding $\mathring x^n$, so that 
along this subsequence $c_n =\Phi_1(\mathring x^n) \uparrow \infty$ and 
\beql{eq-contr}
G^{(n)} (\mathring x^n) / c_n \to 0.
\eeql
Using the continuity of $G^{(n)} (\mathring x^n)$ in $\mathring x^n \in \mathring\cx^{(n)}$, we also can and
 do choose this subsequence so that all particle locations $w_i$ are distinct (i.e., $\mathring x^n$ has exactly $n$ jumps,
at points $w_i$, each jump is by $1/n$). 
For a fixed $\nu \in (0,1/4)$ and each $n$, denote 
$$
a^+_{\nu,n}= a^+_{\nu,n}(\mathring x^n)= \max \left\{w^+ \ge 0 ~|~ \int_{w^+}^\infty w d[-\mathring x^n_{w}]  \ge \nu c_n \right\},
$$
$$
a^-_{\nu,n}= a^-_{\nu,n}(\mathring x^n)= \min \left\{w^- \le 0 ~|~ \int_{-\infty}^{w^-} |w| d[-\mathring x^n_{w}]  \ge \nu c_n \right\}.
$$
It is easy to see that for a sufficiently small fixed $\delta>0$, uniformly in $n$ and $\nu \in (0,1/4)$,
$$
a^+_{\nu,n} \ge \delta c_n, ~~ a^-_{\nu,n} \le -\delta c_n.
$$
(Otherwise, $\int_{-\infty}^{0} |w| d[-\mathring x^n_{w}] = \int_0^{\infty} w d[-\mathring x^n_{w}] = c_n/2$ could not hold.)
Denote 
$$
q^+ = \mathring x^n_{a^+_{\nu,n}}, ~~ q^- = 1-\mathring x^n_{a^-_{\nu,n}}.
$$
It is easy to see that for any $\delta_1>0$ we can fix a sufficiently small $\nu>0$ such that, uniformly in all large $n$, 
$$
q^+ \le \delta_1, ~~ q^- \le \delta_1.
$$
(If not, since both $a^+_{\nu,n}/c_n$ and $a^-_{\nu,n}/c_n$ are bounded away from $0$,
$\nu$ could not be small.) Let us fix any $\delta_1 \in (0,1/2)$, and then choose a corresponding $\nu$.

Observe that the following holds for each large $n$. The distance between the particles located
at $a^-_{\nu,n}$ and $a^+_{\nu,n}$, let us label them $i^-$ and $i^+$, 
 is at least $2\delta c_n$, and there is at least $(1-2\delta_1)n/2$ (i.e. a non-zero fraction) of particles, located strictly between particles
 $i^-$ and $i^+$. 
 
 Using coupling and Assumption~\ref{cond-ntriv-add} (this is the only place where this assumption is used),
 we obtain that there exists $\delta_2>0$ such that
 \beql{eq-drift-diff}
 \bar\zeta^{(\mathring x^n)}_{i^-} - \bar\zeta^{(\mathring x^n)}_{i^+} \ge \delta_2.
 \eeql
 Indeed, when a job arrival selects both $i^-$ and $i^+$, then, conditioned on this event, 
 by monotonicity, $i^-$ receives at least as much average new workload as $i^+$.
 Also, clearly, the probability that both $i^-$ and $i^+$ are selected by an arriving job, under the condition that $i^-$
 (respectively, $i^+$) is selected vanishes as $n\to\infty$.
 Thus, to show \eqn{eq-drift-diff}, is suffices to consider two probability spaces, $E^-$ and $E^+$, 
 obtained by conditioning on the two non-intersecting events, when, respectively, $i^-$ (but not $i^+$) and $i^-$ (but not $i^+$) is selected
 by a job.
We can couple the constructions of the spaces $E^-$ and $E^+$, so that, w.p.1, 
$i^-$ (on $E^-$) and $i^+$ (on $E^+$) receive equal component sizes, while 
all other selected servers and their component sizes are equal.
Consider an arrival of a job of class $j$, satisfying Assumption~\ref{cond-ntriv-add}.
 At least one of the following two conditions holds for $\mathring x^n$: (a) at least $[(1-2\delta_1)n/2]/2$
 particles are located in $[0,a^+_{\nu,n})$ or (b) at least $[(1-2\delta_1)n/2]/2$
 particles are located in $(a^-_{\nu,n},0]$. In the case (a), with probability bounded away from $0$ 
 all other selected particles are in $[0,a^+_{\nu,n})$; then, Assumption~\ref{cond-ntriv-add}(i) easily implies \eqn{eq-drift-diff},
 because, with bounded away from $0$ probability, $i^-$ will receive new workload equal to its component size, while $i^+$ will receive 
 a strictly smaller (by a bounded away from $0$ quantity) new workload. 
 In the case (b), with probability bounded away from $0$ 
 all other selected particles are in $(a^-_{\nu,n},0]$; then, Assumption~\ref{cond-ntriv-add}(ii) easily implies \eqn{eq-drift-diff},
 because, with bounded away from $0$ probability, $i^-$ will receive some positive (bounded away from $0$) new workload,
 while $i^+$ will receive zero new workload. This completes the proof of \eqn{eq-drift-diff}.
 
 From \eqn{eq-drift-diff} we see that 
 $$
 \mbox{either}~ \bar\zeta^{(\mathring x^n)}(a^-_{\nu,n}) - \bar\zeta^{(\mathring x^n)}(0) \ge \delta_2/2~~\mbox{or}~~
 \bar\zeta^{(\mathring x^n)}(a^+_{\nu,n}) - \bar\zeta^{(\mathring x^n)}(0) \le -\delta_2/2.
 $$
 In either case, $G^{(n)} (\mathring x^n) \ge (\delta_2/2) \nu c_n$, which contradicts \eqn{eq-contr}.
$\Box$

\subsection{Proof of \eqn{eq-generator-bound}.}
\label{sec-generator-bound}

First, consider a fixed state $p^n$ and consider the expected increment $\Delta$ of $\Phi_2(p^n)$ upon a job arrival in this state. 
Denote by $\zeta_i=\zeta_i^{(\mathring x^n)}$ the (random) displacement of $w_i$ due to a job arrival.
Then, 
$$
\Delta= \E \frac{1}{n} \sum_i (w_i+\zeta_i)^2 - \frac{1}{n} \sum_i w_i^2 = 
\E \frac{1}{n} \sum_i [2w_i \zeta_i + \zeta_i^2],
$$
where expectation $\E$ is with respect to the distribution of the job class, selection set and the component sizes.
For $\zeta_i$ we have:
$$
\zeta_i = \kappa_i - \frac{1}{n} \sum_s \kappa_s, 
$$
where $\kappa_i=\kappa_i^{(\mathring x^n)}$ the (random) amount of workload added to server $i$ due to a job arrival.
Note that $\E \zeta_i$ is the quantity of the order $O(1/n)$, since $\sum_s \kappa_s$ is of order $O(1)$
and $\E \kappa_i $ is of order $O(1/n)$ (because $\P\{i\in S\} = O(1/n)$). Therefore, $\bar \zeta_i = n \E \zeta_i = O(1)$,
and we can write:
$$
\E \sum_i w_i \zeta_i =  \frac{1}{n} \sum_i w_i \bar \zeta_i.
$$
Next,
$$
\E \sum_i \zeta_i^2 \le 2 \sum_i  \E \kappa_i^2  +
2n \frac{1}{n^2} \E (\sum_s \kappa_s)^2 \le C'_1,
$$
because $\sum_i  \E \kappa_i^2$ is uniformly upper bounded by the maximum (job among classes) second moment 
of the total size of all components. Assembling these bounds, we obtain
\beql{eq-n-times-drift}
n\Delta 
\le 2 \frac{1}{n} \sum_i w_i \bar \zeta_i + C'_1,
\eeql
where, recall, $\Delta$ is the expected increment  of $\Phi_2$ upon a job arrival into a fixed state $p^n$ such that $\Phi_2(p^n) \le C_2$.

Now consider the value of the generator $A^{(n)} \Phi_2^{(C_2)}(p^n)$ at point $p^n$. For that, consider 
the expected increment of $\Phi_2^{(C_2)}(p^n(t))$ over a small interval $[0,t/n]$, with $p^n(0)=p^n$.
First, note that, as $t\downarrow 0$, the contribution into this expected increment of the event that more than one job arrives,
is $o(t)$. (Because jobs arrive as Poisson process of rate $\lambda^n n$, and $\Phi_2^{(C_2)}$ is bounded.) With probability
$\lambda^n t + o(t)$ there will be exactly one job arrival in $[0,t/n]$. Moreover, the state into which this arrival occurs is close
to $p^n$, and recall (Lemma~\ref{lem-monotone-cont-basic}) the continuity in $p^n$
of the workload amounts added by a job. Finally, notice that, if a job arrival 
occurs into a state $p^n$ such that $\Phi_2(p^n) \le C_2$, the expected increment of $\Phi_2^{(C_2)}$ does not exceed that
of $\Phi_2$. Using these observations and the estimate \eqn{eq-n-times-drift}, we obtain \eqn{eq-generator-bound},
where $C_1$ is a constant such that $C_1 > \lambda^n C'_1$ for all $n$; we omit the straightforward $\epsilon/\delta$
formalities. 
$\Box$

\section{Mean-field limits}
\label{sec-ml}

\subsection{Mean-field limits and their fixed points}
\label{sec-fsp-general}

In this section we consider a system with infinite or finite frame size $c$,
and study the limits of the process $x^{n}(\cdot)$ as $n\to\infty$. Recall that for each $n$,  $\lambda^n n$ is the arrival rate 
and $\pi_j^n$ are job-class occurrence probabilities,
and we consider the asymptotic regime such that $\lambda^n \to \lambda$
and $\pi_j^n \to \pi_j, ~\forall j,$
 as $n\to\infty$. 
 
\begin{thm}
\label{lem-fsp-unique}
Assume that the initial conditions are such that $x^{n}(0) \to x(0)$
for some fixed $x(0) \in \cx$.  
Then, there exists a deterministic trajectory $x(t), ~t\ge 0,$ with values in $\cx$ and initial state $x(0)$, such that
\beql{eq-conv-to-ml}
(x^n(t), ~t\ge 0) ~ \stackrel{P}{\longrightarrow} (x(t), ~t\ge 0).
\eeql
(That is, the sequence of processes $x^{n}(\cdot)$, with trajectories in the Skorohod space of functions taking values in $\cx$, 
converges 
in probability to the deterministic trajectory $x(\cdot)$.)
As a function of $t$ (i.e., as an element of Skorohod space), $x(t)$ is continuous. 
\end{thm}

\begin{definition}[Mean-field limit (ML)]
The limiting deterministic trajectory $x(\cdot)$ in Theorem~\ref{lem-fsp-unique} will be called the {\em mean-field limit} (ML)
with initial state $x(0) \in \cx$.
\end{definition}

{\em Proof of Theorem~\ref{lem-fsp-unique}}. This proof closely follows the development in Section 8 of \cite{SnSt2020}. 
(The limiting trajectory $x(\cdot)$, which we call a mean-field limit $x(\cdot)$ in this paper, is called a fluid sample path (FSP) $x^c(\cdot)$
in  \cite{SnSt2020}, where $c$ explicitly indicates the frame size.)
There are some differences that we highlight here. 

Just like in \cite{SnSt2020} we can and do assume that, for each $n$, the server
indices $1,\ldots,n$ are assigned to the servers randomly, according to a permutation of $(1,\ldots,n)$ chosen uniformly at random. 
In our case, however, the initial states $x^n(0)$ and $x(0)$ are arbitrary (as opposed to having a special form in \cite{SnSt2020}), satisfying the convergence $x^n(0) \to x(0)$. 
Just like in \cite{SnSt2020}, we have the convergence of the arrival rates, $\lambda^n \to \lambda$. (In \cite{SnSt2020}, to avoid clogging notation, the construction is given under the assumption that $\lambda^n=\lambda$ and $\pi_j^n=\pi_j$ for all $n$. In this paper, 
the fact that $\lambda^n$ and $\pi_j^n$ may depend on $n$ is crucial, so we do not make this exposition-simplifying assumption to avoid confusion.) 

As in  \cite{SnSt2020}, let us denote: $\alpha_j=\lambda \pi_j d_j$ and $\alpha=\sum_j \alpha_j = \lambda \sum_j \pi_j d_j$.  
(Note that $\alpha \le \lambda \bar d$.)
Parameter $\alpha_j$ is naturally interpreted as the limiting rate at which 
a given server is selected by arriving jobs of class $j$, and  then $\alpha$ is the limiting total rate at which a server is selected. 
We will also use notation
$\hat \pi_j = \alpha_j/\alpha = \pi_j d_j/ [\sum_\ell \pi_\ell d_\ell]$ for the limiting probability that an arriving job selecting a server
is of class $j$. 

The definition of the formal single-server workload process $U_1(t)$ in the artificial ``infinite-server system,'' via the dependence set $\bar D_1(t)$, is same as in \cite{SnSt2020} (except, of course, the fact that the job workload placement rule is different). We denote
$x_w(t) \doteq \P\{U_1(t) > w\}, ~w\ge 0$. Let $W^n_i(t)$ be the workload of server $i$ at time $t$. (We drop the superscript $\lambda^n$,
because in the setting of this lemma $\lambda^n$ is fixed for each $n$.) Then, we can use exactly same coupling as in \cite{SnSt2020}, under which, w.p.1, for all sufficiently large $n$, the dependence sets $\bar D_1(t)$ and $D^n_1(t)$ are {\em equivalent} (as defined in \cite{SnSt2020}), and the component size vectors for each job are same. Furthermore, given that 
in the construction of $U_1(t)$ the initial workloads $W_i(0)$ of the servers in the dependence set $\bar D_1(t)$ are i.i.d. drawn 
from distribution $x(0)$, and the initial workloads $W^n_i(0)$ of the servers in the dependence set $D^n_1(t)$ are chosen 
uniformly at random (without replacement) from the set of workloads described by $x^n(0)$, the coupling can be such that, w.p.1,
$W^n_i(0) \to W_i(0)$ for each server in $\bar D_1(t)=D^n_1(t)$. 
Then, using Lemma~\ref{lem-monotone-cont-basic}, we conclude that, under the described coupling, w.p.1,
for every $t\ge 0$,
$$
\lim_{n\to\infty}  W_1^n(t) = U_1(t),
$$
implying, in particular,
$$
W_1^n(t) \Rightarrow U_1(t).
$$
Therefore, for any fixed $t$, at any point $w$ of continuity of $\P\{U_1(t) > w\}=x_w(t)$
(which is almost any $w$ w.r.t. Lebesgue measure),
$$
\lim_{n\to\infty} \P\{W^n_1(t) > w\} = \P\{U_1(t) > w\}.
$$
(This is slightly different from \cite{SnSt2020}, where the above convergence holds for any $w \ge 0$.)
Then, for any fixed $t$, for almost any $w$, we obtain (exactly as in Lemma 12 of \cite{SnSt2020})
$$
x^n_w(t) \stackrel{P}{\longrightarrow} x_w(t).
$$
Recalling that $x(t)$ and all $x^n(t)$ are elements of $\cx$, we obtain that, for any fixed $t$
\beql{eq-conv-t}
x^n(t) \stackrel{P}{\longrightarrow} x(t).
\eeql

From this point on, the proof is different from the development in \cite{SnSt2020}. (In \cite{SnSt2020} it sufficed to consider only
FSPs starting from special initial states $x(0)$. In this paper, we need FSPs starting from arbitrary initial state $x(0) \in \cx$.)

The next step is to show that 
\beql{eq-continuity-in-t}
\mbox{$x(t)$ is continuous in $t$.}
\eeql
This is equivalent to showing that $x_w(t)$ is continuous in $t$ at every point $w$ where $x_w(t)$ is continuous in $w$. 
This is true -- here we use the structure of process $U_1(\cdot)$, namely that, for a fixed $w$, increase or decrease of $x_w(t)$ in a small interval $[t,t+\delta]$
corresponds to the probabilities that the workload 
 crosses point $w$ in the interval $[t,t+\delta]$ from 
left to right or from right to left, respectively. 
(More precisely, crossing point $w$ from left to right means entering interval $(w,\infty]$.
Similarly, crossing from right to left means leaving interval $(w,\infty]$.)
As a result, $x_w(t)$, as as function of $t$, can never jump up,
simply because the jobs selecting the server arrive according to a finite rate Poisson process.
Function $x_w(t)$, of $t$, also cannot jump down 
at a point of continuity (in $w$) of $x_w(t)$. So, we obtain \eqn{eq-continuity-in-t}. 

Finally, from the structure of process $x^n(\cdot)$, it is easy to see that 
$$
\forall \epsilon>0,~~ \exists \delta>0, ~~\forall t, ~~ \lim_n \P\{\sup_{\tau\in [t,t+\delta]} L(x^{n}(\tau),x^{n}(t)) > \epsilon\} = 0.
$$
From this, \eqn{eq-continuity-in-t} and \eqn{eq-conv-t}, we obtain \eqn{eq-conv-to-ml}.
  $\Box$

It follows from the proof of Theorem~\ref{lem-fsp-unique}, that $U_1(t), ~t\ge 0,$ as a random process, 
can be viewed as the following evolution of a single server workload (or movement of a single particle). 
As $t$ increases, $U_1(t)$ and the dependence set $\bar D_1(t)$ evolve in the following fashion. 
At initial time $0$, $U_1(0)$ is chosen randomly according to the distribution $x(0)$, and $\bar D_1(0) = \{1\}$.
Let $\tau_1, \tau_2, \ldots$ be the points of the rate $\alpha$  Poisson process of the (artificial) job arrivals selecting server $1$.
In the interval $[0,\tau_1]$, the workload (particle location) $U_1(t)$ decreases at constant rate $-1$ unless and until it hits zero, and
$\bar D_1(t) = \{1\}$. We can and will adopt 
a (non-essential) convention that $\bar D_1(t)$ and $U_1(t)$ are LCRL. 
(After $U_1(t)$ is defined as a LCRL process, it can always be replaced by its RCLL version if necessary.)
The job arriving and selecting server $1$ at time $\tau_1$ is of some class $j$, and it selects $d_j - 1$ servers besides servers $1$; for illustration purposes, suppose those are servers $2$ and $3$. Then, 
$$
\bar D_1(\tau_1+) = \bar D_1(\tau_1) \cap \bar D_2(\tau_1) \cap \bar D_3(\tau_1) = \{1\} \cap \bar D_2(\tau_1) \cap \bar D_3(\tau_1),
$$
where $\bar D_2(\tau_1)$ and $\bar D_3(\tau_1)$ are the random dependence sets of servers $2$ and $3$ at time $\tau_1$. 
The dependence sets
$\bar D_2(\tau_1)$ and $\bar D_3(\tau_1)$ do not intersect with each other, do not contain server $1$, and constructed the same way 
as described in \cite{SnSt2020} (for $\bar D_1(t)$ with $t=\tau_1$). Moreover, {\em $\bar D_2(\tau_1)$ and $\bar D_3(\tau_1)$ are independent and identically distributed, with their distribution depending only on the time $\tau_1$.} (The independence and equality of distributions here are
in the sense that they would hold if the servers in both $\bar D_2(\tau_1)$ and $\bar D_3(\tau_1)$ would use, independently, the same set of  labels.) This means, in turn, that  $U_2(\tau_1)$ and $U_3(\tau_1)$ are i.i.d. with the distribution $x(\tau_1)$.
Therefore, the value of $U_1(\tau_1+)$ is determined by $U_1(\tau_1)$ and the realizations of $U_2(\tau_1)$, $U_3(\tau_1)$ and those of
the arriving job component sizes. Then,  in $(\tau_1,\tau_2]$, the workload (particle location) $U_1(t)$ again decreases at constant rate $-1$ unless and until it hits zero.
At time $\tau_2$ another job selecting server 1 arrives, which also selects some additional servers, say, $4,5,6$. The values of 
$U_4(\tau_2), U_5(\tau_2), U_6(\tau_2)$ are i.i.d. with the distribution $x(\tau_2)$.
The value of $U_1(\tau_2+)$ is determined by $U_1(\tau_2)$ and the realizations of $U_4(\tau_2), U_5(\tau_2), U_6(\tau_2)$ and those of
the arriving job relica sizes. And so on. 
We will use this interpretation of the process $U(\cdot)$ later. One may say that $U_1(\cdot)$ describes the random 
movement of a single particle in the ``mean-field environment'' defined
by the time-varying distribution $x(\cdot)$.

  The following Lemma~\ref{lem-fsp-props} easily follows from the construction of ML in Theorem~\ref{lem-fsp-unique}
  and the monotonicity of the pre-limit process.
  
  \begin{lem}
  \label{lem-fsp-props}
(i) Continuity. ML $x(\cdot)$ is continuous w.r.t. the initial state $x(0)$. 

(ii) Shift (``Markov'') property. If $x(\cdot)$ is a ML, then for any $\theta \ge 0$ the trajectory
$\hat x(t) = x(\theta+t), ~t\ge 0,$ is the ML with initial state $\hat x(0)=x(\theta)$.

(iii) Monotonicity. If $x(0) \le \hat x(0)$, then the ML $x(\cdot)$ and $\hat x(\cdot)$
with these initial conditions are such that $x(t) \le \hat x(t)$ for all $t\ge 0$.
  \end{lem}

 \begin{definition}
 Point $x\in \cx$ is called a {\em fixed point of the mean-field limit} (ML-FP), if $x(t) \equiv x$ is an ML.
 \end{definition}
 
   Note that, trivially, $x^{**}$ is an ML-FP for the infinite-frame system 
   and $x^\emptyset$ is the only ML-FP for the zero-frame system. ($x^{**}$ and $x^\emptyset$ are defined in \eqn{eq-inf-def}
   and \eqn{eq-empty-def}.)
  
\begin{lem}
\label{lem-fsp-props-3}
 Consider the ML $x(\cdot)$, starting with the empty initial state $x(0)=x^\emptyset$. Then $x(t)$ is monotone non-decreasing and $x(t) \to x^*$,
 where $x^*$ is an ML-FP, and moreover it is the {\em minimum ML-FP} (ML-MFP). (That is, $x^* \le x$ for any ML-FP $x$.)
If $x^* \ne x^{**}$, then $x^*$ is proper. 
  \end{lem}
  
  {\em Proof.} ML $x(t)$ is obviously monotone non-decreasing, and so the convergence $x(t) \to x^*$ for some element $x^* \in \cx$
  does hold. That $x^*$ must be an ML-FP  follows from Lemma~\ref{lem-fsp-props}. Indeed, the shifted versions of $x(\cdot)$, namely
  $x(\theta+t), ~t \ge 0,$ are ML. A sequence of such ML, with $\theta \uparrow \infty$, is such that
  $x(\theta+t) \to x^*$ at any fixed $t\ge 0$; in particular, the initial states $x(\theta) \to x^*$.
  We also have that the entire ML $x(\theta+\cdot)$ converges to the stationary trajectory at point $x^*$.
  The latter, by Lemma~\ref{lem-fsp-props}, is a ML. Thus, $x^*$ is an ML-FP.
  
 The ML-FP $x^*$ must be the minimum ML-FP, 
 which follows by comparing the ML $x(t)$ starting from the empty state with the stationary FL $\hat x(t)=x$ 
 for any ML-FP $x$. Since $x(t) \le x$ for all $t$, and $x(t) \to x^*$, we have $x^* \le x$.
  
  Suppose $x^*$ is not proper, but $x^* \ne x^{**}$, that is $x^*_\infty = \epsilon^* \in (0,1)$. This leads to a contradiction by the following argument, based on the monotonicity. For any $\epsilon < \epsilon^*$, arbitrarily close to $\epsilon^*$, the following property holds for the pre-limit processes starting from the empty state:  
  \beql{eq-prop-123}
  \forall C>0, ~\exists T>0, ~\mbox{for all sufficiently large $n$}, ~~ \P\{W^n_1(T) \ge C\} \ge \epsilon.
  \eeql
  Note that property \eqn{eq-prop-123} cannot hold when $\epsilon > \epsilon^*$. On the other hand, if we choose $T_1$
  large enough so that \eqn{eq-prop-123} holds with $C$ and $T$ replaced by $C+T$ and $T_1$, consider the time interval $[0,T_1+T]$,
  and use the monotonicity, we obtain:
  $$
  \mbox{for all sufficiently large $n$}, ~~ \P\{W^n_1(T_1+T) \ge C\} \ge \epsilon + (1-\epsilon) \epsilon.
  $$
  This means that \eqn{eq-prop-123} holds with $\epsilon$ replaced by $\epsilon + (1-\epsilon) \epsilon$,
  which is greater than $\epsilon^*$ when $\epsilon$ is close to $\epsilon^*$. A contradiction.
  $\Box$

 \subsection{ML-FP as a solution to a functional differential equation}
  
 In this subsection we will show that an element $x\in \cx$ is an ML-FP if and only if it is a solution to a certain functional differential equation (FDE); and we derive some properties of such solutions that we will need later. It is not hard to ``guess'' the basic form of the FDE. It is also 
 not hard to define the FDE formally and show that any ML-FP must satisfy the FDE. But, proving the converse (that any FDE solution is an ML-FP) and establishing the FDE solution properties require more involved analysis.
 
 We start with some definitions. For any element $x \in \cx$, let us use a slightly abusive notation $x_{[0,w]}$ for  
 $(x_u, ~u\in [0,w])$, i.e. for function $x_\cdot$ with the domain truncated to $[0,w]$. Suppose we randomly choose a job class
 $j$ according to the distribution $(\pi_j)$, and then place the class $j$ job on some abstract $d_j$ servers, whose workloads are chosen randomly and independently according to distribution $x$. Then, let
 $$
 h(x_{[0,w]})
 $$
 denote the expected number of those servers whose workload, as a result, will cross level $w$ from left to right
 or, more precisely, will enter the interval $(w,\infty]$. 
Clearly, the values of $x_u$ for $u>w$ have no effect on this expectation -- that is why the expectation depends only on $x_{[0,w]}$.

There is an alternative, equivalent way to define functional $h$, which we also be useful.
Recall our notation $\alpha = \lambda \sum_j \pi_j d_j$,  $\hat \pi_j = \pi_j d_j/ [\sum_\ell \pi_\ell d_\ell]$. 
For a fixed $x \in \cx$, $w\ge 0$, and a job class $j$, consider CDF $J_j^{(x,w)}(u), ~u\ge 0,$ defined as follows.
Consider an (abstract) ``tagged'' server, with workload $w$, selected by a class $j$ job; suppose the remaining $d_j -1$ servers, selected 
by this job, have random independent workloads chosen according to distribution $x$;
then $J_j^{(x,w)}(\cdot)$ is the CDF of the random amount of workload that will be added to the tagged server. 
Denote $\bar J_j^{(x,w)}(u) = 1-J_j^{(x,w)}(u)$.
Then, it is easy to check that
$$
\lambda h(x_{[0,w]}) = \alpha \sum_j \hat \pi_j \int_{0-}^w [-dx_u] \bar J_j^{(x,u)}(w-u),
$$
or
$$
h(x_{[0,w]}) = \sum_j \pi_j d_j \int_{0-}^w [-dx_u] \bar J_j^{(x,u)}(w-u).
$$
Denote further 
$$
J^{(x,w)}(u) = \sum_j \hat \pi_j J_j^{(x,w)}(u), ~~~ \bar J^{(x,w)}(u) = 1-J^{(x,w)}(u).
$$
Then,
$$
\lambda h(x_{[0,w]}) = \alpha \int_{0-}^w [-dx_u] \bar J^{(x,u)}(w-u).
$$

It very intuitive that any ML-FP should, speaking informally at this point, satisfy the following functional differential equation (FDE):
\beql{eq-diff-main}
- x'_w = \lambda h(x_{[0,w]}),
\eeql
with the LHS being the steady-state scaled rate at which particles cross point $w$ from right to left (leave the interval $(w,\infty]$), and
the RHS  being the steady-state scaled rate at which particles cross point $w$ from left to right (enter the interval $(w,\infty]$). 
The definitions and results that follow, give a precise meaning of the FDE \eqn{eq-diff-main},
and show that, indeed, $x$ is an ML-FP if and only if it is a solution to this equation.

Let us define the following operator $A$ which maps $\cx$ into itself. 
Let a distribution (``environment'') $x \in \cx$ be fixed; $x$ may or may not be proper. Consider a single particle moving in $\R_+$ as follows. Unless and until the particle jumps, it moves left at speed $1$; if/when the particle ``hits'' the left boundary $0$, it stays there until next jump; there is a Poisson process, rate $\alpha$, of time points when the the particle jumps forward; if the jump occurs when the particle location is $w$, the jumps size distribution is $J^{(x,w)}(\cdot)$, independent of the process history (besides current location $w$). Then $Ax \in \cx$ is the stationary distribution of the particle. Note that, if a particle is located at infinity, it stays there forever. Consequently, if $x$ is not proper, i.e. $x_\infty > 0$, then $(Ax)_\infty = x_\infty$.

\begin{definition}
A fixed point $x$ of the operator $A$ we will call an {\em OP-FP}.
\end{definition}
   
  \begin{lem}
  \label{lem-fsp-is-op-fp}
If $x$ is a ML-FP, then it is an OP-FP. (We will show later in Lemma~\ref{lem-fp-iff-de} that, in fact, $x$ is an ML-FP if and only if it is an OP-FP.)
  \end{lem}
  
  {\em Proof.}   Consider the ML construction given in Theorem~\ref{lem-fsp-unique}, and the interpretation of the 
  process $U_1(\cdot)$ given after the proof of Theorem~\ref{lem-fsp-unique}. 
  This construction defines the distribution of a single particle location
at any time $t$. Since $x$ is a ML-FP, $x(t)=x$ is a stationary ML, and therefore the distribution of the particle location is $x$ for any $t$.
This means that, if the particle is selected by a job at time $t$, then the locations of other selected particles are i.i.d. having distribution $x$. This, in turn, means exactly that $x$ is a fixed point of the operator A, i.e.
it is a stationary distribution of a particle moving within environment $x$.
$\Box$

\begin{definition}
An element $x=(x_w)\in \cx$
with infinite or finite frame size $c$ we call a DE-FP if it is an $\alpha$-Lipschitz [and then $(\lambda \bar d)$-Lipschitz] function,
satisfying the FDE \eqn{eq-diff-main} at every point $0<w<c$, where the derivative $x'_w$ exists (which is a.e.  w.r.t. Lebesgue measure). 
\end{definition}

  \begin{lem}
  \label{lem-op-fp-iff-de-fp}
  Any OP-FP $x=(x_w)$ is a DE-FP. Its frame size $c$ is equal to the truncation parameter $c$ of the model.
   Unless $c=\infty$ and OP-FP $x=x^{**}$, or $c=0$ and OP-FP $x=x^{\emptyset}$, we have $x_0 \in (0,1)$.
  \end{lem}

  {\em Proof.}  By definition, OP-FP $x$ is the stationary distribution of the location of a single particle 
  evolving in the (stationary) environment $x$. Obviously, when $c=0$, OP-FP $x=x^{\emptyset}$ is the unique OP-FP, and it trivially satisfies 
  the definition of DE-FP. 
   Similarly, when $c=\infty$, OP-FP $x=x^{**}$ trivially satisfies 
  the definition of DE-FP.
  So, in the rest of the proof it suffices to -- and we will -- consider an OP-FP $x$ such that either $c \in (0,\infty)$ or [$c=\infty$ and $x\ne x^{**}$]. 
  
  We can immediately observe that $x_0 \in (0,1)$ must hold, because in steady-state the particle must spend non-zero fraction of time both at the left boundary $0$ and away from it. Function $x_w$ is $\alpha$-Lipschitz, because in steady-state, for any interval 
  $[u,w] \subset [0,c)$ and any $\Delta$-long time interval, the average amount of time the particle spends in $[u,w]$ 
  is at most $\alpha (w-u) \Delta$, and therefore $x_u - x_w \le \alpha (w-u)$. It is also easy to see that $c' = \min\{w~|~x_w=0\} < c$
  is impossible, because otherwise in steady-state the particle must spend a non-zero fraction of time to the right of point $c'$ -- a contradiction.
    Therefore, $x_w >0$ for all $w < c$, which means that the frame size (of $x$) is equal to $c$.
  
  Finally, consider any point $w \in (0,c)$ where the derivative $x'_w$ exists (which is almost any point w.r.t. Lebesgue measure).
  Then \eqn{eq-diff-main} must holds, because its 
  LHS is the steady-state rate at which the particle crosses point $w$ from right to left (more precisely, leaves $(w,\infty]$),
  while the RHS is the steady-state rate at which the particle crosses point $w$ from left to right (more precisely, enters $(w,\infty]$).
$\Box$

  \begin{lem}
  \label{lem-h-Lipschitz}
 For any fixed $w \ge 0$, the functional $h(x_{[0,w]})$ on $\cx$ is $\bar d^2$-Lipschitz, w.r.t. the sup-norm on 
 $x_{[0,w]}$, namely
 $$
 | h(x_{[0,w]}) - h(y_{[0,w]})| \le \bar d^2 \sup_{0\le u \le w} | x_u - y_u |, ~~\forall x,y \in \cx.
 $$
 (Note that here $x_u$ and $y_u$ are not necessarily continuous.)
 \end{lem}
  
    {\em Proof.}  
It suffices to prove the following.
 Consider any fixed $x_{[0,w]}$ and any fixed $\delta>0$. 
 Consider $y_{[0,w]}$ defined by: $y_w=x_w$ and
 $y_u = (x_u+\delta) \wedge 1$ for $u < w$; so, it is $x_u$ shifted up by $\delta$, but not above $1$, in the interval $[0,w)$.
 Comparing $h(x_{[0,w]})$ and $h(y_{[0,w]})$, we see that, by monotonicity (Lemma~\ref{lem-monotone-cont-basic}), 
 the latter cannot be smaller, and recalling the definition of $h$ (as the level $w$ up-crossing rate) and using coupling, we see that the maximum possible up-crossing rate difference is $\bar d^2 \delta$.
Indeed, $\bar d \delta$ is the upper bound on the probability that a job arrival will select at least one server from the subset of scaled size $\delta$; then,
even if we assume that when it happens for $h(y_{[0,w]})$,  all components -- at most $\bar d$ -- of this job will cause their selected servers to 
 up-cross level $w$, the increase in the expected number of up-crossings compared to that for $h(x_{[0,w]})$ is at most $(\bar d \delta) \bar d$.
$\Box$

  \begin{lem}
  \label{lem-de-sol-unique}
 A DE-FP $x=(x_w)$, for a given initial condition $x_0$, exists, is unique and the dependence of this DE-FP on $x_0$ is continuous.
    \end{lem}
  
 {\em Proof.}  Suppose $x_0>0$. Note that \eqn{eq-diff-main} is equivalent to the
integral equation
\beql{eq-diff-main-integral}
 x_w = x_0 - \lambda \int_0^w h(x_{[0,u]}) du.
 \eeql
 A solution $x$ of \eqn{eq-diff-main-integral} is a fixed point of the operator
 \beql{eq-diff-main-operator}
 y_w  = x_0 - \lambda \int_0^w h(x_{[0,u]}) du.
 \eeql
 In a small interval $[0,\Delta]$, operator \eqn{eq-diff-main-operator} maps the space of $\alpha$-Lipschitz non-increasing functions
 with initial condition $x_0$,
 into itself. Using the Lipschitz property (Lemma~\ref{lem-h-Lipschitz}) of the functional $h(x_{[0,w]})$,
 we see that operator \eqn{eq-diff-main-operator} is a contraction, This, by the standard Picard iterations method
 establishes existence and uniqueness of solutions to \eqn{eq-diff-main-integral} in $[0,\Delta]$, and then
 in the entire interval $[0,c]$, where $c$ it the value of $w$ at which the the solution ``hits'' the $w$-axis.
 (In other words, $c$ is the frame of the solution.) This also establishes the continuity 
 of the solution in $x_0$ as long as $x_0>0$. For $x_0=0$, trivially, the solution to \eqn{eq-diff-main-integral} is unique, equal 
 to $x^\emptyset$. It is straightforward to check that the solution to \eqn{eq-diff-main-integral} is continuous w.r.t. $x_0$ at $x_0=0$ 
 as well.
$\Box$

It immediately follows from Lemma~\ref{lem-de-sol-unique} that $x^{**}$ is the unique DE-FP for $x_0=1$,
and $x^{\emptyset}$ is the unique DE-FP for $x_0=0$.

Consider any DE-FP $x$ with infinite frame size. From \eqn{eq-diff-main-integral}, we have
\beql{eq-diff-main-integral-full}
 \frac{x_0 - x_\infty}{\lambda} = \int_0^\infty h(x_{[0,u]}) du.
 \eeql
 Recall the definition of $h(x_{[0,u]})$ as the average number of particles up-crossing level $u$ upon a job arrival, 
 such that the locations of particles within its selection set are i.i.d. with distribution $x$. (Note that, if $x_\infty >0$, 
 with probability $x_\infty >0$ a selected particle is at infinity.)
 Then, the RHS of \eqn{eq-diff-main-integral-full} is equal the expected total amount of new workload added by such a job
 arrival to those selected servers whose workload is finite.
 
\begin{lem}
  \label{lem-de-fp-props}
  Consider any
  DE-FP $x=(x_w)\ne x^{**}$, with (finite or infinite) frame size $c>0$. (Note that, necessarily, $x_0 < 1$.)
  Then, $x$ satisfies 
   the following additional properties. Function
    $x_w$  is strictly decreasing in $[0,c)$.
  The right derivative $x'^{+}_w$ exists at 
every point $w\in [0,c)$. The left derivative $x'^{-}_w$ exists, and  $x'^{-}_w \le x'^{+}_w$, at 
every point $w\in (0,c)$. (These derivatives are non-positive.)
Equation \eqn{eq-diff-main} holds for the right derivative $x'^{+}_w$ at every point $w\in [0,c)$.
The right derivative $x'^{+}_w$ is RCLL, the left derivative $x'^{-}_w$ is LCRL,
and each is negative, bounded away from $0$ on any set $[0,a]$, where $a<c$. 
  \end{lem}

  {\em Proof.}  Consider a DE-FP $x$.
 It suffices to -- and we will -- consider the case $x_0 \in (0,1)$. (If $x_0=0$, then $x=x^\emptyset$ and the result trivially holds.)
It is easy to see that, when $x_w$ is continuous in $w$ (and it is),  $h(x_{[0,w]})$ is RCLL as a function of $w$.
Moreover, the left limit of $h(x_{[0,w]})$ at any point $u>0$ is greater than of equal to the right limit.
(The left limit at point $u$ may be strictly greater than the right limit when there is a non-zero probability that
the particle ``jumping forward" from point $0$  lands exactly at point $u$. Function $h(x_{[0,w]})$ is RCLL, as opposed to
LCRL, due to our convention that up-crossing level $w$ means entering $(w,c]$, as opposed to $[w,c]$.)
Given that $x_w$ is Lipschitz, and then $x'^{-}_w = x'^{+}_w$ a.e., we obtain that the left derivative is LCRL,
and that $x'^{-}_w \le x'^{+}_w$ must hold at every $w\in (0,c)$. 
Then $h(x_{[0,w]})$ must be positive, bounded away from zero on any interval $[0,a] \subset [0,c)$ -- otherwise,
by the right-continuity of $h(x_{[0,w]})$, there would exist the smallest point $u\ge 0$ such that $h(x_{[0,w]})=0$,
which is impossible, as easy follows from the $h(x_{[0,w]})$ definition.
$\Box$

Consider a DE-FP $x \ne x^{**}$. Its frame size $c$ can be finite or infinite. 
Denote by $\phi_u = x^{-1}_u, ~u\in [0,1]$ its inverse function.
The inverse $\phi_u$ has the following structure: $\phi_u = 0$ for $1 \ge u \ge x_0$; $\phi_u$ is strictly decreasing for
$x_\infty < u \le x_0$; $\phi_{x_\infty} = \lim_{u\downarrow x_\infty} \phi_u$; $\phi_u = \infty$ for $u < x_\infty$ (if $x_\infty >0$).  Moreover, on any interval $[\kappa,x_0] \subset (x_\infty,x_0]$, $\phi_u$ is Lipschitz. (This follows from $x_w$ having negative, bounded away from $0$ derivative in any interval $[0,a] \subset [0,c)$.) 

\begin{lem}
  \label{lem-hor-dist-no-local-max}
Consider any two DE-FP, $x$ and $\hat x$, neither equal to $x^{**}$. (They may have different frame sizes, finite or infinite, in general.)  Denote by $\phi_u, ~u\in [0,1],$ and $\hat \phi_u, ~u\in [0,1],$ their corresponding inverse functions. 
Denote $b = \max\{x_\infty, \hat x_\infty\} < 1$, and consider the difference 
$\psi_u \doteq \hat \phi_u - \phi_u$ for $u\in (b,1]$.
Then, if for some  $\chi \in (b,1]$
$$
\psi_\chi = \max_{u\in [\chi,1]} \psi_u > 0,
$$
the difference $\psi_u$ is non-increasing in $(b,\chi]$. 
  \end{lem}

 {\em Proof.}  The proof is by contradiction. If lemma does not hold, then the following holds for some $\zeta$ satisfying $\zeta \le \chi$ and
 $\zeta \in (b,1]$, 
 and some $\beta \in (b,\zeta)$:
 \beql{eq-local-max1}
 \psi_\zeta = \max_{u\in [\zeta,1]} \psi_u > 0
 \eeql
 and 
 \beql{eq-local-max2}
 \psi_u < \psi_\zeta, ~~\forall u\in (\beta,\zeta).
 \eeql
 Note that, in this case, we must have $\zeta \le \min\{\hat x_0, x_0\}$. Denote $w \doteq \phi_\zeta < \hat w \doteq \hat \phi_\zeta$. 
 From the definition of point $\zeta$, we have
 \beql{eq-local-max3}
 \hat \phi_\zeta - \hat \phi_u \le \phi_\zeta - \phi_u, ~~\forall \zeta \le u \le 1.
 \eeql
 Recall the definition of functional $h$ and compare $h(\hat x_{[0,\hat w]})$ and $h(x_{[0,w]})$. 
 By \eqn{eq-local-max3} and monotonicity (Lemma~\ref{lem-monotone-cont-basic}) we obtain
 \beql{eq-h-monotone}
 h(\hat x_{[0,\hat w]}) \le h(x_{[0,w]}). 
 \eeql
 Let us denote by $s_*$ the essential lower bound on the positive component sizes:
 \beql{eq-min-repl-size}
 s_* = \inf_{s\ge 0} \{ \exists j, ~\pr\{0< \xi_1 <s\}  > 0\}.
 \eeql
 Consider two cases: $s_* < \hat w$ and   $s_* \ge \hat w$. 
  
Note that left derivative of $\phi'^-_u=1/x'^+_z$, where $u=x_z$. So, the left derivative of $\phi_u$ is determined by the right derivative of $x_z$.
 
 If $s_* < \hat w$, it is easy to observe that the inequality \eqn{eq-h-monotone} must be strict, $h(\hat x_{[0,\hat w]}) < h(x_{[0,w]})$,
 implying $x'^+_{\hat w} > x'^+_{\hat w}$ and then $\hat \phi'^-_{\zeta} < \phi'^-_\zeta$.
Therefore $\psi'^-_\zeta<0$,
which contradicts \eqn{eq-local-max2}.

If $s_* \ge \hat w$, it is easy to see that, at points $u<\zeta$, sufficiently close to $\zeta$, inequality $\hat \phi'^-_u \le \phi'^-_u$ must hold,
and we, again, obtain contradiction to \eqn{eq-local-max2}.

We see that the assumptions \eqn{eq-local-max1} and \eqn{eq-local-max2} lead to contradictions, which completes the proof.
   $\Box$
   
As a corollary of Lemma~\ref{lem-hor-dist-no-local-max} we obtain the following 

 \begin{lem}
  \label{cor-de-fp-order}
Suppose we have two DE-FP $\hat x$ and $x$, neither equal to $x^{**}$, with $\hat x_0 > x_0$. Then, the difference of their inverses,
$\hat x^{-1}_u - x^{-1}_u$ is non-increasing in $u\in [0,1]$ (if we adopt convention that $+\infty - z = +\infty$ for any $z\le +\infty$).
In particular, these two DE-FP cannot intersect: 
if $x$ has a finite frame $[0,c]$, then $\hat x_w > x_w$ for all $w\in [0,c]$;
if $x$ has infinite frame, 
then $\hat x_w > x_u$ for all $w \in [0,\infty)$.
  \end{lem} 

The next Lemma~\ref{cor-de-fp-finite-unique} also easily follows from Lemma~\ref{lem-hor-dist-no-local-max}.

 \begin{lem}
  \label{cor-de-fp-finite-unique}
For any $0\le c < \infty$ (and a given $\lambda$) a DE-FP $x$ with finite frame $[0,c]$ is unique.
  \end{lem}
  
{\em Proof.}  Suppose not, that is two different DE-FP $\hat x$ and $x$ with the same finite frame $[0,c]$ exist.
WLOG, suppose $\hat x^{-1}_u > x^{-1}_u$ for some $u$.
  Then, by the continuity of the inverses, the difference $\hat x^{-1}_u - x^{-1}_u$ would attain a positive maximum at some point $u>0$,
  while $\hat x^{-1}_0 - x^{-1}_0 = 0$.
  This is impossible due to
  Lemma~\ref{lem-hor-dist-no-local-max}.
  $\Box$
   
  \begin{lem}
  \label{lem-fsp-finite-unique}
For any $0\le c < \infty$ (and a given $\lambda$) consider the system with finite frame $[0,c]$. Then, the following holds.\\
(i) The ML-FP $x^*=x^{*,\lambda}$ (with finite frame $[0,c]$) is unique, and is equal to the unique DE-FP and unique OP-FP, with this frame. \\
(ii) The convergence $x(t) \to x^*$ holds uniformly in all initial states $x(0)$.
  \end{lem}

 {\em Proof.} (i) Consider the $x^*$, which is the ML-MFP for the finite frame $[0,c]$. 
 It is an ML-FP, and then OP-FP and then DE-FP. The latter is unique by Lemma~\ref{cor-de-fp-finite-unique}.
 
 (ii) Any ML $x(\cdot)$ is ``sandwiched'' between the MLs starting from the minimum ($x^\emptyset$) and maximum ($x^{**,c}$) 
 initial states $x(0)$. The former is monotone non-decreasing, converging to the ML-MFP $x^*$. The latter is monotone
 non-increasing, converging to some $\hat x^* \ge x^*$. But, the limit $\hat x^*$ is also an ML-FP, and therefore is equal
 to the unique ML-FP $x^*$.
  $\Box$
  
      \begin{lem}
  \label{lem-fp-iff-de}
An element $x \in \cx$ is a ML-FP if and only if it is a DE-FP (and then a OP-FP as well).
  \end{lem}
  
 {\em Proof.}  The `only if' part has already been proved. Let us prove the 'if' part. 
 Consider the following truncated version $\hat x$ of this DE-FP. Namely, fix any $\phi > x_\infty$ and consider the unique finite $c$,
 such that $x_c = \phi$; then $\hat x = x_{[0,c]}$. This $\hat x$ is the unique ML-FP (= OP-FP) for the $\phi$-reduced
 system (as defined in Section~\ref{sec-equiv-view}) with finite frame $c$. Let us compare the ML $x(\cdot)$ for the original system with initial state $x(0)=x$ and the ML $\hat x(\cdot)$ for the $\phi$-reduced truncated system with initial state $\hat x(0) = \hat x= x_{[0,c]}$. 
From the ML construction we see that, for any $t \in [0, c/2]$ and any $w \in [0, c/2]$, we must have 
$x_w(t) = \hat x_w(t)$.
Since $\hat x(t) \equiv \hat x$, and $c$ can be made arbitrarily large by choosing $\delta$ sufficiently close to $x_\infty$,
we must have $x(t) \equiv x$, that is $x$ is an ML-FP.
 $\Box$

The results given so far in this section are for any fixed parameter $\lambda$ and any fixed distribution $\pi=(\pi_j)$ of job classes. 
We will also need the following extension of Lemma~\ref{lem-de-sol-unique}.

  \begin{lem}
  \label{lem-de-sol-unique2}
 The dependence of a DE-FP $x=(x_w)$ on $(x_0,\lambda, \pi)$ is continuous.
     \end{lem}

{\em Proof.} First, note that if for a fixed $x_{[0,w]}$ and $\lambda$ we consider the dependence of $h(x_{[0,w]})$ on $\pi$, 
this dependence is $\bar d$-Lipschitz. (This is easy to see using coupling.) From this and Lemma~\ref{lem-h-Lipschitz}, 
$h(x_{[0,w]})$ is Lipschitz in $(x_{[0,w]}, \pi)$.

Consider a converging sequence $(x_0^{(\ell)},\lambda^{(\ell)},\pi^{(\ell)}) \to (x_0,\lambda,\pi)$.
The corresponding sequence of (uniformly Lipschitz) DE-FP solutions $x^{(\ell)}$ is such that its any subsequence has a further subsequence
along which $x^{(\ell)} \to x$. Using the fact that \eqn{eq-diff-main-integral} holds for each  $x^{(\ell)}$ and the properties of $h(\cdot)$, it is straightforward to see that $x$ must satisfy \eqn{eq-diff-main-integral} for the triple $(x_0,\lambda,\pi)$. $\Box$

 \section{Proof of Theorem~\ref{thm-main-finite}}
 \label{sec-proof-main-th-finite}
 
 \sasha{This proof relies on the properties of the mean-field limits and their fixed points that we derived in Section~\ref{sec-ml}. 
 Most importantly, it uses Lemma~\ref{lem-fsp-finite-unique} -- the uniqueness of ML-FP $x^{*,\lambda}$ for a finite-frame system
 and the attraction to it, $x(t) \to x^{*,\lambda}$, for any ML $x(\cdot)$.}
 
 (i) Consider any subsequence of $n$, along which 
 the stationary distributions converge, i.e. $x^n(\infty) \Rightarrow x(\infty)$, where $x(\infty)$ is some random element in $\cx$.
 We know that the convergence of initial states $x^n(0) \to x(0)$ implies convergence $x^n(\cdot) \to x(\cdot)$ to the ML $x(\cdot)$
 with initial state $x(0)$, and the dependence of a ML on its initial state is continuous. Then, by Theorem 8.5.1 in \cite{Liptser_Shiryaev},
 the distribution of $x(\infty)$ is stationary for the deterministic process evolving along the ML trajectories. 
 But, by Lemma~\ref{lem-fsp-finite-unique}, all MLs converge to the ML-MFP $x^{*,\lambda}$.
 
 (ii) The continuity of $x^{*,\lambda}$ in $\lambda$ easily follows from the facts that the ML starting from the empty state 
 (and converging up to $x^{*,\lambda}$) on any finite time interval is continuous with respect to $\lambda$. 
 It is also easy to see that, as $\lambda \uparrow \infty$, the $x^{*,\lambda}$ cannot converge to anything except $x^{**,c}$.
 $\Box$

\section{A DE-FP/ML-FP uniqueness condition. Proof of Theorem~\ref{thm-main-D-monotonicity}}
\label{sec-cond-unique}

\sasha{We are now in position to formally state the following DE-FP/ML-FP uniqueness condition, which appears in the statements
of Theorems~\ref{thm-main} and \ref{thm-main-D-monotonicity}.}

\begin{condition}[DE-FP/ML-FP uniqueness]
\label{cond-unique}
We say that this condition holds for a fixed $\lambda \ge 0$, if 
for any sufficiently small $b \ge 0$, 
there is at most one DE-FP $x$ with infinite frame and $x_\infty = b$. 
\end{condition}

\subsection{Proof of Theorem~\ref{thm-main-D-monotonicity}(i).}  

\sasha{This proof relies on the properties of the mean-field limits and their fixed points that we derived in Section~\ref{sec-ml}.}

Fix any $\lambda \ge 0$. Let us prove that Condition~\ref{cond-unique} holds. In fact we will prove that it holds for 
 any (not necessarily small) $b \in [0,1)$. 
The proof is by contradiction. Suppose, there exist two different DE-FP, $x$ and $\hat x$, both with infinite frames, and such that 
$x_\infty = \hat x_\infty= b \in [0,1)$. Since these DE-FP are different, $x_0 \ne \hat x_0$. Let $x_0 < \hat x_0$ for concreteness.
Then, by
 Lemma~\ref{cor-de-fp-order}, the difference of the inverses $\hat x^{-1}_u - x^{-1}_u$ is non-increasing in $u\in (b,1]$, 
 where $b=x_\infty = \hat x_\infty$. Therefore, by D-monotonicity, the average amount of new workload brought by job arrival 
 selecting servers with the workloads being i.i.d. with distribution $\hat x$ does not exceed that corresponding to
 distribution $x$. On the other hand, by \eqn{eq-diff-main-integral-full}, those average new workload amounts
 corresponding to $\hat x$ and $x$ are $(\hat x_0 - \hat x_\infty)/\lambda$ and $(x_0 - x_\infty)/\lambda$, respectively,
 with the former being strictly greater than the latter. The contradiction completes the proof.
  $\Box$

\subsection{Proof of Theorem~\ref{thm-main-D-monotonicity}(ii)} 

\sasha{This proof is self-contained, in that it does not use any other results or proofs in this paper.}

It suffices to prove the theorem's statement (ii) for the class $j$ itself, because each sub-class also has i.i.d. component sizes and the same proof applies.

Consider the servers $1,\ldots,d_j$ in the selection set, in the order of their appearance in vector $D$.
Consider the non-trivial case when $k_j < d_j$. (Otherwise, $\E \eta_j(D) = \E \sum_{i=1}^{d_j} \xi^{(j)}_i$ is the expected total size -- it does not depend on $D$.) In the rest on the proof we will drop job index $j$, and write $d,k,H, \eta(D), \xi_i, \xi$ instead of $d_j,k_j,H_j, \eta_j(D), \xi_i^{(j)}, \xi^{(j)}$, to simplify notation.

Fix vector $D$ and $1 < m \le d$. Denote
$$
\tilde \eta(y) \doteq \eta(D_1, \ldots, D_{m-1}, D_m+y, D_{m+1}, \ldots, D_{d}), ~y\ge 0.
$$
It will suffice to show that, for any $D$, any fixed $1 < m \le d$, 
\beql{eq-eta-monotone-y}
\mbox{function $\E\tilde \eta(y)$ is non-increasing in $y\ge 0$.}
\eeql
Consider deterministic function $\tilde \eta(y), ~y\ge 0,$ for a fixed realization $\xi=(\xi_1,\ldots,\xi_{d})$ of component sizes; we will write 
$\tilde \eta(y,\xi)$ to indicate the dependence on $\xi$. It is elementary to see that 
$|\tilde \eta(y+\Delta,\xi) - \tilde \eta(y,\xi)| \le \Delta d$ for any $D,y,\xi$ and $\Delta \ge 0$, and therefore
function $\tilde \eta(y,\xi), ~y\ge 0,$ is  $d$-Lipschitz. This immediately implies that 
\beql{eq-eta-lipschitz}
\mbox{function $\E\tilde \eta(y)$ is $d$-Lipschitz, and then absolutely continuous.}
\eeql

We will prove \eqn{eq-eta-monotone-y}, first, in the special case when $H(\cdot)$ is absolutely continuous, and then in the general case.

{\em Case (a): $H(\cdot)$ is absolutely continuous.} This condition is equivalent to $H(\cdot)$ {\em not} having an atom at $y_{max}<\infty$,
where $y_{max}$ is the essential upper bound of the distribution support.
(When $y_{max}=\infty$ the condition automatically holds.)
We will show that, for any $D$, any fixed $1 < m \le d$, 
\beql{eq-eta-deriv22}
\frac{d \E \tilde \eta(y)}{dy} \le 0, ~~\forall y \ge 0.
\eeql

Fix $y\ge 0$. 
Denote: $\tilde D_i = D_i$ for $i\ne m$, and $\tilde D_i = D_i+y$ for $i = m$; $\tilde Z_i = \sum_1^i \tilde D_\ell$.
Since $H(\cdot)$ is absolutely continuous,
\beql{eq-no-ties}
\mbox{W.p.1, $\xi$ is such that
$\tilde Z_{i_1}+\xi_{i_1} \ne \tilde Z_{i_2}$ and 
$\tilde Z_{i_1}+\xi_{i_1} \ne \tilde Z_{i_2}+\xi_{i_2}$ for all $i_1 \ne i_2$.}
\eeql
It will suffice to consider only those realizations of $\xi$ satisfying \eqn{eq-no-ties}.
Denote by $E(i_1,i_2)$ the following event:  {\em $\{$server $i_1$ is such that the value 
$\tilde Z_{i_1}+\xi_{i_1}$ is the $k-th$ smallest among values $\tilde Z_{i}+\xi_i, ~i=1,\ldots,d$, and $ \tilde Z_{i_1}+\xi_{i_1} \in (\tilde Z_{i_2}, \tilde Z_{i_2}+\xi_{i_2})\}$}. We will refer to event $E(i_1,i_2)$ as {\em $\{i_1$ cancels $i_2\}$}. Then, it is elementary to see that function $\tilde \eta(y,\xi), ~y\ge 0,$ is such that
\beql{eq-eta-deriv}
\frac{d \tilde \eta(y,\xi)}{dy} = \sum_{i_1=1}^{m-1} \sum_{i_2=m}^{d} [-\bI(E(i_1,i_2) + \bI(E(i_2,i_1)].
\eeql
Then, by bounded convergence theorem, 
\beql{eq-eta-deriv2}
\frac{d \E\tilde \eta(y)}{dy} = \E \frac{d \tilde \eta(y,\xi)}{dy} = \sum_{i_1=1}^{m-1} \sum_{i_2=m}^{d} [-\pr(E(i_1,i_2)) + \pr(E(i_2,i_1))].
\eeql
Observe that in \eqn{eq-eta-deriv2} we alway have
\beql{eq-key-compare}
\pr(E(i_1,i_2)) \ge \pr(E(i_2,i_1)).
\eeql
Indeed, for a pair of indices $i_1 < i_2$ consider the event 
$$
E= E(i_1,i_2) \cup E(i_2,i_1) =
\{\mbox{either $i_1$ cancels $i_2$ or $i_2$ cancels $i_1$}\}.
$$
Then 
$$
\pr(E(i_1,i_2)~|~E) \ge \pr(E(i_2,i_1)~|~E),
$$
because when both components $i_1$ and $i_2$ are served simultaneously, 
and the elapsed service time of $i_2$ is $\tau_2$, 
the elapsed service time of $i_1$ is $\tau_1=\tau_2 + (\tilde Z_{i_1}-\tilde Z_{i_2}) \ge \tau_2$, and therefore, by IHR, conditioned on one of 
these two components completing service, the probability that it will be $i_1$ is $h(\tau_1)/[h(\tau_1)+h(\tau_2)] \ge 1/2$.
From \eqn{eq-eta-deriv2} and \eqn{eq-key-compare}, we obtain \eqn{eq-eta-deriv22}.

{\em Case (b): General.} It remains to consider the case when $H(\cdot)$ has an atom at $y_{max}<\infty$.
For any $D$, there exists its arbitrarily small perturbation such that condition \eqn{eq-no-ties} holds for all $y\ge 0$, except maybe a finite set of points $\{y_{\ell}\}$. Then, for such a perturbed $D$, repeating the argument given for the case (a), we can show that 
$\frac{d \E \tilde \eta(y)}{dy} \le 0$ for any $y \not\in \{y_{\ell}\}$, and then $\E \tilde \eta(y)$ is non-increasing. 
(Recall that $\E \tilde \eta(y)$ is absolutely continuous.) It remains to use the continuity \eqn{eq-eta-cont}, to obtain that
\eqn{eq-eta-monotone-y} holds for any $D$.
$\Box$

 \section{Proof of Theorem~\ref{thm-main}}
 \label{sec-proof-main-th}
 
\sasha{This proof relies -- in essential way -- on both the properties of systems with finite $n$ (Theorems~\ref{thm-crit-load-finite} and \ref{th-closeness}) and on the properties of the mean-field limits, derived in Section~\ref{sec-ml}. 
(The proof does not use Theorems~\ref{thm-main-finite} and \ref{thm-main-D-monotonicity}.)}

\sasha{We now informally describe the key ideas of the proof. In essence, it reduces to proving the following for any fixed $\rho<1$. Consider a converging sequence $\lambda^n \to \lambda$ 
such that, for all $n$, the system is stable, with the load exactly $\rho^n(\lambda^n) = \rho$. (Such a sequence exists by Theorem~\ref{thm-crit-load-finite}.) Moreover, we can choose this sequence so that the convergence of stationary 
distributions holds, $x^n(\infty) \Rightarrow x(\infty)$, where $x(\infty)$ is some random element. 
Recall that an ML trajectories $(x(t), t\ge 0)$, which are deterministic, 
closely approximate the random trajectories of process $(x^n(t), t\ge 0)$ when $n$ is large. 
In fact, the distribution of $x(\infty)$ is invariant for the deterministic process evolving along ML trajectories.
Using monotonicity 
it is not hard to see that, w.p.1, $x(\infty) \ge x^{*,\lambda}$, where $x^{*,\lambda}$ is the ML-MFP for the parameter $\lambda$,
and $x_0^{*,\lambda} \le \rho$. We need to show that, in fact, $x^n(\infty) \Rightarrow x^{*,\lambda}$ and $x_0^{*,\lambda}=\rho$. (The latter easily follows from the former.) The ML-FP/DE-FP uniqueness (Condition~\ref{cond-unique}) is employed to
show that any ML $x(\cdot)$ with proper initial state $x(0)$ dominating $x^{*,\lambda}$ must converge to $x^{*,\lambda}$ as $t\to\infty$. The intuition here is that, due to Condition~\ref{cond-unique}, we can always pick an improper DE-FP $\tilde x^*$ which  is just slightly larger than $x^{*,\lambda}$. Then any proper $x(0)$ is dominated by $\hat x(0)$ which is $\tilde x^*$ ``shifted sufficiently far to the right.'' If we consider the ML $\hat x(\cdot)$ starting from $\hat x(0)$, it will be non-increasing (by monotonicity) and it must converge down to $\tilde x^*$, because otherwise the limit $\lim_{t\to\infty} \hat x(t)$
would be a DE-FP different from $\tilde x^*$, thus violating Condition~\ref{cond-unique}. Now, if any ML $x(\cdot)$,
starting from any proper $x(0)$ dominating $x^{*,\lambda}$ is attracted to $x^{*,\lambda}$, this will imply that 
the probability of $x(\infty)$ being proper is equal to the probability of $x(\infty)=x^{*,\lambda}$.
Then we will show that the probability of $x(\infty)$ being improper is equal to the probability of $x(\infty)=x^{**}$.
The intuition here is that, by Theorem~\ref{th-closeness}, for any $n$, in steady-state, the ``relative distances between particles are of the order $O(1)$;'' therefore, if say $0.99$-th quantile of the distribution given by $x^n(t)$ is very large, then
with high probability any other fixed quantile will be very large as well. 
This will allow us to show that $x(\infty)$ is concentrated on at most two points: $x^{*,\lambda}$ and the infinite element $x^{**}$. 
It will still remain to show that $x^{**}$ cannot have a positive probability. 
This will be done by contradiction. Namely, assuming $x(\infty)=x^{**}$ has some positive probability, we obtain that,
uniformly in all large $n$, in steady-state, a fixed, say $0.99$-th, quantile of $x^n(0)$ is very large with some positive probability.
Using the mentioned above convergence of ML $\hat x(\cdot)$ down to DE-FP $\tilde x^*$, we will be able to conclude
that, with high probability, a lower, say $0.98$-th, quantile of $x^n(T_1)$ (at some large time $T_1$) would have to be 
``far to the left'' from the $0.99$-th quantile of $x^n(0)$. This would imply that, for large $n$, in steady-state,
the expected distance between the $0.98$-th and $0.99$-th quantiles of $x^n(\infty)$ could be arbitrarily large;
that would contradict Theorem~\ref{th-closeness}.}

 The formal proof starts here.
 
 Consider our asymptotic regime with $n\to\infty$. Fix $0< \rho < 1$. For each $n$ by Theorem~\ref{thm-crit-load-finite} we can and do choose $\lambda^n$ such that the steady-state load is equal to exactly $\rho$, i.e.
 $\E x^n_0(\infty) = \rho$. We can always choose a subsequence of $n$, along which 
 $\lambda^n \to \lambda$, where, necessarily, $\lambda \in (0,\infty)$; consider any such subsequence.
 
 Consider  any further subsequence of $n$, along which 
 the stationary distributions converge, i.e. $x^n(\infty) \Rightarrow x(\infty)$, where $x(\infty)$ is some random element in $\cx$.
 We know that the convergence of initial states $x^n(0) \to x(0)$ implies convergence $x^n(\cdot) \to x(\cdot)$ to the ML $x(\cdot)$
 with initial state $x(0)$, and the dependence of a ML on its initial state is continuous. Then, by Theorem 8.5.1 in \cite{Liptser_Shiryaev},
 the distribution of $x(\infty)$ is stationary for the deterministic process evolving along the ML trajectories. 
 Then we must have that, w.p.1, $x(\infty) \ge x^*$, where $x^*=x^{*,\lambda}$ is the ML-MFP,
 because for any ML, $x(t) \ge x^\emptyset (t)$ and $x^\emptyset (t) \uparrow x^*$. 
 Moreover, $x^*_0 \le \rho$. (Otherwise, $\liminf_n \E x^n_0(\infty)  > \rho$.) Let us prove that, in fact, 
 \beql{eq-limit-concentrates}
 \mbox{$x^*_0 = \rho$ and $x(\infty)=x^*$ w.p.1.}
 \eeql
 
First, we prove that
\beql{eq-proper-mfp}
\P\{x_\infty(\infty) =0\} = \P\{x(\infty)=x^{*}\},
\eeql
that is, almost surely, a proper state is equal to $x^*$.
To demonstrate \eqn{eq-proper-mfp}, 
we will prove that ML with any proper initial state $x(0)$ (dominating $x^*$) is such that to $x(t) \to x^*$;
moreover, we will prove uniform convergence, in the following sense. For any $\epsilon>0$ and $C>0$, there exists $\delta=\delta(\epsilon,C)>0$ such that for any $A>0$ there exists $T=T(\epsilon,C,\delta,A)>0$ such that, for any ML with $x_A(0) \le \delta$, 
\beql{eq-uniform-conv-down}
x_w(t) < x^*_w + \epsilon, ~ \forall w\le C, ~~ \forall t \ge T.
\eeql
Indeed, pick a sufficiently small $\delta>0$, and for it DE-FP 
$\tilde x^{*}$ with $\tilde x^{*}_\infty=\delta$ and $\tilde x^{*}_w < x^*_w+ \epsilon/2, ~w\le C$. Then, $x(0)$ is dominated by the state $\hat x(0)$, which is equal to $\tilde x^{*}$ shifted right by $A$. (That is,  $\hat x_w(0) = 1$ for $w<A$, and $\hat x_w(0) = \tilde x^{*}_{w-A}$ for $w \ge A$.)
The ML $\hat x(\cdot)$, starting from $\hat x(0)$, is monotone non-increasing, because, by construction, if we would have a left regulation boundary at $A$,
 $\hat x(0)$ would be a stationary ML.  Moreover, $\hat x(t) \downarrow \tilde x^{*}$. 
 (It must converge to something, which is a ML-FP = DE-FP, and this DE-FP must be equal to $\tilde x^{*}$,
 because otherwise we would have two DE-FP with the same limit $\delta$ at $w=\infty$ -- that is impossible by Condition~\ref{cond-unique}.)
 Therefore, $x_w(t) \le \hat x_w(t) < \tilde x^{*}_w  + \epsilon/2, ~w\le C,$ for all $t\ge T$ for a sufficiently large $T$,
 which proves \eqn{eq-uniform-conv-down}. In turn, \eqn{eq-uniform-conv-down} in particular implies that 
 any ML $\hat x(\cdot)$ with proper initial state $\hat x(0)\ge x^*$, converges to $x^*$. This, combined with $x(\infty) \ge x^*$,
proves \eqn{eq-proper-mfp}.
 
Next, we claim that
\beql{eq-improper-infinite}
\P\{x_\infty >0\} = \P\{x=x^{**}\}.
\eeql
The proof of \eqn{eq-improper-infinite} is given below in Section~\ref{sec-improper-infinite}.

The combination of \eqn{eq-proper-mfp} and \eqn{eq-improper-infinite} shows that the limit $x(\infty)$ is concentrated
on two points, $x^*$ and $x^{**}$. To prove \eqn{eq-limit-concentrates},
it remains to show that $\epsilon_1 \doteq \P\{x=x^{**}\} = 0$. Indeed, this will imply that $x(\infty)=x^*$ w.p.1,
and $x^*_0=\rho$ must hold (because otherwise we would have $\lim_n \E x^n_0 = x^*_0 < \rho$.)
We will prove that $\epsilon_1 = 0$ by contradiction. Suppose not: $\epsilon_1 = \P\{x=x^{**}\} > 0$.

Pick a sufficiently small $\delta>0$, and for it DE-FP 
$\tilde x^{*}$ with $\tilde x^{*}_\infty=\delta$; $\tilde x^{*}$ is close to and dominates $x^*$.
For a given $A>0$, consider, as we already did above, the ML $\hat x(\cdot)$ starting from
$\hat x(0)$, which is equal to $\tilde x^{*}$ shifted right by $A$. Recall that $\hat x(t)$ converges down to $\tilde x^{*}$.
Also, using monotonicity, it is easy to see that, if we increase $A$, it can only decrease $\hat x(t) - \hat x(0)$ for any $t$.
(Because increasing $A$ means moving the regulation boundary $0$ to the left w.r.t. $A$.)
These observations imply that,
for any $A_1>0$ we can find a sufficiently large $T_1>0$ and a sufficiently large $A$, such that we have to following property.
Denote 
$q^n_1(t) = [x^n(t)]^{-1}_{\delta}$ and $q^n_2(t) = [x^n(t)]^{-1}_{2\delta}$.
(These are, respectively, the $(1-\delta)$-th and $(1-2\delta)$-th quantile of the distribution 
$x^n(t)$.)
Then, as $n\to\infty$,
$$
\P\{q^n_2(T_1) - q^n_1(0) < - A_1 \} \to 1 ~~\mbox{uniformly in $x^n(0)$ with $q^n_1(0) \ge A$.}
$$
Consider the process $x^n(\cdot)$ in the stationary regime, in the interval $[0,T_1]$. 
When $n$ is large: (a) with probability close to $\epsilon_1$, $q^n_1(0)$ is large and
$q^n_2(T_1) - q^n_1(0) < - A_1$; 
(b) with probability close to $1-\epsilon_1$, $q^n_2(T_1) - q^n_1(0)$ is close to 
$[x^*]^{-1}_{2\delta} - [x^*]^{-1}_{\delta} \le 0$.
Then, 
\beql{eq-drift-quantile}
\limsup_{n\to\infty} \E [q^n_2(T_1) - q^n_1(0)] \le - \epsilon_1 A_1/2,
\eeql
assuming that
\beql{eq-st-upper-bound}
q^n_2(T_1) - q^n_1(0) \le_{st} V_n,
\eeql
where $V_n$ is a family of non-negative uniformly integrable random variables.
The proof of \eqn{eq-st-upper-bound} is given below in Section~\ref{sec-st-upper-bound}.
Thus, \eqn{eq-drift-quantile} indeed holds. But, since $A_1$ can be arbitrarily large,
this implies that $\E [q^n_1(\infty) - q^n_2(\infty)] \to \infty$ as $n\to\infty$, which is impossible.
(It follows from Theorem~\ref{th-closeness} that for any fixed $0 <\eta <1$,
$\E |[x^n(\infty)]^{-1}_\eta -\bar x^n(\infty)|$ is uniformly bounded in $n$.) The contradiction completes the proof 
of \eqn{eq-limit-concentrates}. 

Thus,  we have proved that $x^n(\infty) \Rightarrow x^*$
for the scenario where $\rho^n = \E x^n_0(\infty) = \rho$,
and for a subsequence of $n$ along which  $\lambda^n \to \lambda$ for some 
$\lambda$, and $x^*=x^{*,\lambda}$ being the ML-MFP for this $\lambda$, with $x^{*,\lambda}_0=\rho$.

It is easy to see that
the dependence of the ML-MFP $x^{*,\lambda}$ on $\lambda$ is strictly increasing continuous, in particular 
the dependence of $x^{*,\lambda}_0$ on $\lambda$ is strictly increasing continuous. From here we obtain
that the limiting $\lambda$ in the above scenario is uniquely determined by $\rho$,
and the dependence $\lambda=\lambda(\rho)$ is strictly increasing continuous in $[0,1)$.
Denote $\bar \lambda = \lim_{\rho\to 1} \lambda(\rho)$, and denote by $\rho=\rho(\lambda)$ the function inverse to $\lambda(\rho)$.
If we fix any $\lambda\in [0,\bar \lambda)$ and let $\lambda^n \to \lambda$, it is easy to see, 
using the properties of the functions $\rho(\lambda)$ and $\lambda(\rho)$, along with
 monotonicity, 
that the convergence $x^n(\infty) \Rightarrow x^{*,\lambda}$ is the only possibility.
$\Box$

\subsection{Proof of \eqn{eq-improper-infinite}.}
\label{sec-improper-infinite}

The intuition for this fact is rather simple. If $x_\infty(\infty)>0$ with positive probability $\epsilon_1$, then for a very small fixed $\eta>0$ and all large $n$ the inverse (quantile) $[x^n(\infty)]^{-1}_{\eta}$ is very large with probability close to $\epsilon_1$. However,
by Theorem~\ref{th-closeness}, the inverses $[x^n(\infty)]^{-1}_{\eta}$ and $[x^n(\infty)]^{-1}_{1-\eta}$ are within distance $O(1)$
with high probability, and therefore $[x^n(\infty)]^{-1}_{1-\eta}$ is also very large with probability close to $\epsilon_1$.
The formal proof is as follows.

By Theorem~\ref{th-closeness}, for any $\epsilon>0$ and $\eta>0$, there exists a sufficiently large $C>0$ such that, for any fixed 
$0 < \eta < 1$, 
for any $n$, $x^n(\infty)$ is such that 
\beql{eq-distances}
\P\{\mbox{any of the distances between $\bar x^n(\infty)$,
$[x^n(\infty)]^{-1}_{\eta}$ and $[x^n(\infty)]^{-1}_{1-\eta}$ exceeds $C$}\} < \epsilon.
\eeql

We only need to consider the case when 
$$
\P\{x_\infty(\infty) > 0\} = \epsilon_1>0.
$$
(Otherwise \eqn{eq-improper-infinite} is trivial.)
We need to show that this implies
\beql{eq-improper-infinite2}
\P\{x(\infty)= x^{**}\} = \epsilon_1.
\eeql
Indeed, pick any positive $\epsilon_2<\epsilon_1$, arbitrarily close to $\epsilon_1$. 
Pick $u>0$, sufficiently small so that $\P\{x_\infty(\infty) > u\} > \epsilon_2$. 
Pick a small $\epsilon>0$ (compared to 
$\epsilon_2$), pick $\eta \in (0,u)$, and then the corresponding $C>0$, so that property \eqn{eq-distances} holds. 
Observe that for an arbitrarily large $C_1>0$, for all sufficiently large $n$, $\P\{x^{(n)}_{C_1}(\infty) > u\} > \epsilon_2$
(where we use the fact that $\{x_a > u\}$ is an open set in $\cx$),
which means that $\P\{[x^n(\infty)]^{-1}_{\eta} \ge C_1\} > \epsilon_2$.
But this means that $\P\{[x^n(\infty)]^{-1}_{1-\eta} \ge C_1-C\} > \epsilon_2-\epsilon$.
From here, $\P\{x_w(\infty) \ge 1-\eta, ~0\le w < C_1-C\} \ge \epsilon_2-\epsilon$ 
(where we use the fact that $\{x_w \ge u, ~0 \le w < a\}$ is a closed set in $\cx$).
Since this is true for an arbitrarily large $C_1$,
$\P\{x_\infty(\infty) \ge 1-\eta\} \ge \epsilon_2-\epsilon$. Since this is true for an arbitrarily small positive $\eta$, 
arbitrarily small positive $\epsilon$,
and $\epsilon_2$ arbitrarily close to $\epsilon_1$, we obtain \eqn{eq-improper-infinite2}. $\Box$

\subsection{Proof of \eqn{eq-st-upper-bound}.}
\label{sec-st-upper-bound}

The intuition for the proof is that we upper bound the value of $q_2^n(T_1)$ at $T_1$ by considering the ``worst case'' scenario when $x^n$
is ``maximum possible'' subject to $q_1^n(0)=A$ and neither services nor cancellations occur in $[0,T_1]$. 
For this worst case we show that the increase of the appropriately defined ``mean'' of $x^n(t)$ in $[0,T_1]$
is uniformly integrable; then $q_2^n(T_1)-A$ must be uniformly integrable as well. The formal argument is as follows.

We construct the stochastic upper bound $V_n$ as follows. Suppose $q_1^n(0)=A$. Then $q_2^n(T_1)-q_1^n(0)
= q_2^n(T_1)-A$ is stochastically dominated by the
random variable $V_n = [x^n(T_1)]^{-1}_{2\delta}- [x^n(0)]^{-1}_{\delta} = [x^n(T_1)]^{-1}_{2\delta}- A$  in the ``worst case'' when: 
(a) $x^n(0)$ is such that there are $\delta n$ servers with infinite workload and the remaining $(1-\delta) n$ servers
have workload $A$, i.e. $x^n_w(0)=1$ for $w< A$, and $x^n_w(0)=\delta$ for $w\ge A$; (b) no service occurs in $[0,T_1]$;
(c) there is no cancellation in $[0,T_1]$ -- the component sizes of each arriving job are the actual amounts of workload added by the job.
Denote by ${\bar {\bar x}}^n(t) \doteq \int_0^\infty [x^n_w(t)-\delta] dw$ the ``mean'' of $x^n(t)$, which only takes into account the servers 
with finite workloads and ``ignores'' the infinite-workload servers. Then, 
using the fact that the arrival process is Poisson and the component sizes have finite second moment, 
we see that, as $n\to\infty$, the increase ${\bar {\bar x}}^n(T_1) - {\bar {\bar x}}^n(0) = {\bar {\bar x}}^n(T_1) -A$ (for the worst case process) 
has converging mean and vanishing variance. Therefore, the family ${\bar {\bar x}}^n(T_1) - A$ is uniformly integrable. From here we easily obtain that (for the worst case process) the family $V_n = [x^n(T_1)]^{-1}_{2\delta}- A$ 
is uniformly 
integrable, because otherwise ${\bar {\bar x}}^n(T_1) - A$ could not be uniformly integrable.
$\Box$

\section{Proof of Theorem~\ref{thm-ihr}}
\label{sec-proof-ihr} 

\sasha{The proof relies on Theorems~\ref{thm-main} and \ref{thm-main-D-monotonicity}.
Under the theorem assumptions, Condition~\ref{cond-unique} holds for all $\lambda\ge 0$ by 
Theorem~\ref{thm-main-D-monotonicity}, and 
Assumption~\ref{cond-second-moment} holds as well. Assumption ~\ref{cond-ntriv-add}(ii) also holds automatically.}
Note that  Assumption~\ref{cond-ntriv-add}(i) does {\em not} necessarily hold. (For example, in the case when for all job classes
w.p.1 all component sizes are equal.) So, we cannot apply Theorem~\ref{thm-main} immediately.
However, it is easy to ``get around'' a situation when Assumption~\ref{cond-ntriv-add}(i) does not hold,
and still use Theorem~\ref{thm-main}, using the following argument based on monotonicity and continuity.

In addition to the original system, consider the following modified system. We add an artificial class $\tilde j$, which has a small probability 
$\pi^{(\epsilon)}_{\tilde j}=\epsilon>0$, 
while all other classes $j$ probabilities are $\pi^{(\epsilon)}_j = (1-\epsilon) \pi_j$;
class $\tilde j$
has $k_{\tilde j} < d_{\tilde j}$, and has, say, i.i.d. exponential component sizes. With this artificial class added,
the modified system does satisfies Assumption~\ref{cond-ntriv-add}(i), while still satisfying 
Condition~\ref{cond-unique} (for all $\lambda\ge 0$) and 
Assumptions~\ref{cond-second-moment} and \ref{cond-ntriv-add}(ii). 
For the modified system we do have Theorem~\ref{thm-main}. 

We consider a sequence of modified systems,
parameterized by $\epsilon \downarrow 0$. For each fixed $\epsilon$, let 
$\tilde \rho^{(\epsilon)}(\lambda), ~\lambda \in [0, \tilde \lambda^{(\epsilon)})$ be the 
corresponding function $\rho(\lambda)$ from Theorem~\ref{thm-main}, 
where we deliberately use notations $\tilde \rho^{(\epsilon)}$ and $\tilde \lambda^{(\epsilon)}$
instead of $\rho$ and $\bar \lambda$. If we use the convention $\tilde \rho^{(\epsilon)}(\lambda) = 1$ for $\lambda \ge \tilde \lambda^{(\epsilon)}$, the function $\tilde \rho^{(\epsilon)}(\cdot)$ is continuous in $[0,\infty)$. 
Now, for each $\epsilon$, we will use the change of variable $\beta = (1-\epsilon) \lambda$,
and consider the function 
$$
\rho^{(\epsilon)}(\beta) = \tilde \rho^{(\epsilon)}(\beta/(1-\epsilon)), ~~\beta \in [0, \bar \lambda^{(\epsilon)}),
$$
where $\bar \lambda^{(\epsilon)} = (1-\epsilon)\tilde \lambda^{(\epsilon)}$.
(Function $\rho^{(\epsilon)}(\cdot)$ is continuous in $[0,\infty)$.) 
The reason for this change of variable is that, for a fixed $\beta$,
$$
\frac{\beta}{1-\epsilon} \pi^{(\epsilon)}
$$
is decreasing as $\epsilon\to 0$. Then
it follows from the monotonicity properties in Lemma~\ref{lem-monotone1} that the sequence of functions $\rho^{(\epsilon)}(\cdot)$ is non-increasing as $\epsilon\to 0$,
with $\bar \lambda^{(\epsilon)}$ non-decreasing and such that $\bar \lambda^{(\epsilon)} \uparrow \bar \lambda^{(0)} < \infty$.
The functions $\rho^{(\epsilon)}(\cdot)$ converge down to some non-decreasing RCLL function $\rho^{(0)}(\cdot)$,
in the sense of convergence at every point of continuity. 

The next step is to observe that function $\rho^{(0)}(\cdot)$ is in fact continuous, including at point $\bar \lambda^{(0)}$,
which means $\rho^{(0)}(\beta) \uparrow 1$ as $\beta \uparrow \bar \lambda^{(0)}$. Indeed,
recall that for each $\epsilon>0$, each $\beta \in [0, \bar \lambda^{(\epsilon)})$, and corresponding
$\lambda = \beta/(1-\epsilon) \in [0,\tilde \lambda^{(\epsilon)})$, we have a well defined DE-FP $x$ for the modified system,
with $x_0<1$. If $\rho^{(0)}(\cdot)$ would have a positive jump at some point $a \in [0, \bar \lambda^{(\epsilon)}]$,
then using Lemma~\ref{lem-de-sol-unique2} we would be able to obtain two different proper infinite-frame DE-FP, $x$ and $y$, for the original system.
This would contradict Condition~\ref{cond-unique}. Thus, $\rho^{(0)}(\cdot)$ is indeed continuous.
Moreover, for each $\lambda \in [0, \bar \lambda^{(0)})$, again using Lemma~\ref{lem-de-sol-unique2},
for the original system with parameter $\lambda$, 
we obtain a DE-FP $x$ with $x_0 = \rho^{(0)}(\lambda)$, which must be the unique ML-MFP due to Condition~\ref{cond-unique}.

Given the continuity of $\rho^{(0)}(\cdot)$, we have the uniform convergence 
$\rho^{(\epsilon)}(\cdot) \stackrel{\|\cdot\|}{\longrightarrow}\rho^{(0)}(\cdot)$. 
Fix $\lambda \in [0, \bar \lambda^{(0)})$ and a small $\epsilon >0$.
Compare the original system with parameter $\lambda$ and the modified system with parameter $\lambda/(1-\epsilon)$.
Any limit of stationary distributions of the original system must be ``sandwiched'' between its ML-MFP (below)
and the ML-MFP for the modified system (above). But those two ML-MFP can be made arbitrarily close to each other
by making $\epsilon$ small. We conclude that the limit of stationary distributions for the original system
is concentrated on its ML-MFP. Then conclusion of Theorem~\ref{thm-main} for the original system easily follows,
with the function $\rho(\cdot)$ being $\rho^{(\epsilon)}(\cdot)$.
$\Box$

\section{Further discussion of the DE-FP/ML-FP uniqueness condition (Condition~\ref{cond-unique}) and additional steady-state asymptotic independence results for infinite-frame systems}
\label{sec-further-conditions}

\sasha{The following Theorem~\ref{thm-majorization} gives a ``majorization'' sufficient condition for the Condition~\ref{cond-unique}
to hold for a given $\lambda$. Informally speaking, it states the following.
Suppose, a non-trivial (i.e. not equal to $x^{**}$) infinite-frame
DE-FP exists for a system with strictly large arrival rates, for all classes and subclasses, then those in the original system with given $\lambda$. Then Condition~\ref{cond-unique}
holds for the original system with the given $\lambda$.} 

\begin{thm}
\label{thm-majorization}
Consider our model with the arrival rate parameter $\lambda$, and the job classes' probability distribution $\pi=(\pi_j)$ as defined for the model.
WLOG assume that the set of classes $j$ includes all the classes of the original system and all their sub-classes
(where some or all sub-classes may have zero probability of occurrence). 
Suppose there exist other parameters, $\hat \lambda$ and $\hat \pi=(\hat \pi_j)$, such that: (a) $\hat \pi_j >0$ for all classes and sub-classes;
(b) $\lambda \pi < \hat \lambda \hat \pi$; (c) for parameters $\hat \lambda$ and $\hat \pi$ there exists an infinite-frame DE-FP $\hat x \ne x^{**}$.
Then, Condition~\ref{cond-unique} holds for this $\lambda$ (and $\pi$).
\end{thm}

{\em Proof.} It suffices to consider the non-vacuous case when the ML-MFP $x^*$ for $\lambda$ (and $\pi$) is proper; this
 $x^*$ is a DE-FP with infinite frame and $x^*_\infty = 0$. Let us now fix a small $b \ge 0$, and consider a $b$-reduced system
 (see Section~\ref{sec-equiv-view}). Then, for the $b$-reduced system, $\tilde \lambda^{(b)} = \lambda/(1-b)$ is the arrival rate parameter
 let $\tilde \pi^{(b)} = (\tilde \pi_j^{(b)})$ denote the job class distribution. Clearly, when $b$ is small, we have 
 $$
 \tilde \lambda^{(b)} \tilde \pi^{(b)} < \hat \lambda \hat \pi.
 $$
 Since $\hat x$ is a DE-FP for $\hat \lambda$ and $\hat \pi$, by monotonicity we obtain that ML-MFP $\tilde x^{*,(b)}$
  is proper (because it must be dominated by $\hat x$). Moreover, the dependence of $\tilde x^{*,(b)}$, which is also a DE-FP,
  on $b$ is continuous, with $\tilde x^{*,(0)} = x^*$. But then 
$$
x^{*,(b)} = b + (1-b) \tilde x^{*,(b)}
$$
is an infinite-frame DE-FP for parameters $\lambda$ and $\pi$, and the dependence of $x^{*,(b)}$
on $b$ is continuous, with $x^{*,(0)} = x^*$. Therefore, as $b$ is increasing from $0$ while remaining small, 
$x^{*,(b)}$ is {\em the} family of DE-FP for $\lambda$ and $\pi$, as the initial condition $x^{*,(b)}_0$ is increasing
from $x^*_0$. Therefore, for all sufficiently small $b\ge 0$, $x^{*,(b)}$ is the only infinite-frame DE-FP $x$ for $\lambda$ and $\pi$,
with $x_\infty=b$.
$\Box$

Using Theorem~\ref{thm-majorization}, along with adjusted versions of the proof of Theorem~\ref{thm-main}, 
we can obtain some additional steady-state asymptotic independence results. 

First, consider the case when the system has inherently sub-critical load. Suppose, for each class $j$ there is some upper bound $\bar s_j$ on the expected amount of work $\E \eta_j(D)$ brought by a class $j$ job, which is uniform in the workload-differential vector $D$ of the selection set. (The definition of $\eta_j(D)$ is in Section~\ref{sec-cond-sufficient}.) \sasha{Note that 
this does {\em not} assume some uniform stochastic upper bound on $\eta_j(D)$ -- only a uniform upper bound 
on the {\em expectation}.}
We say that the system is
{\em inherently sub-critical} if
\beql{eq-cubcrit}
\bar \rho =\lambda \sum_j \pi_j \bar s_j < 1.
\eeql
Trivially, the expected total size of all components, $d_j \E \xi_1^{(j)}$, can always be chosen as $\bar s_j$ for any class $j$.
\sasha{One example of a non-trivial bound $\bar s_j$ is for a class $j$ which is $D$-monotone decreasing. 
Then
$$
\bar s_j = \sup_D \E \eta_j(D) = \E \eta_j(0,\ldots,0)
$$
is the tight upper bound on $\E \eta_j(D)$, which can be found, at least numerically. This includes the case
when the distribution 
of components sizes is such that they are i.i.d. with IHR. (Or, more generally, it may be a mixture of such distributions.)}

Now, analogously to Definition~\ref{def-d-monotone}, a class $j$ 
is called $D$-monotone increasing, if its distribution $F_j$ of component sizes is such that $\E \eta_j(D)$ is non-decreasing in $D$. For such a class $j$,
\beql{eq-555}
\bar s_j = \sup_D \E \eta_j(D) = \E \eta_j(\infty,\ldots,\infty) = k_j \E \xi_1^{(j)}
\eeql
is the tight upper bound on $\E \eta_j(D)$. 
A distribution $H(\cdot)$ on $\R_+$ is said to have {\em decreasing hazard rate} (DHR), if
condition \eqn{eq-ihr-cond} holds with `$\le$' replaced by `$\ge$'. 
Then, a proof analogous to that of Theorem~\ref{thm-main-D-monotonicity}(ii) applies to show that {\em a class $j$
is $D$-monotone increasing when its 
distribution $F_j$ is such that the component sizes are i.i.d. DHR.}
(Or, more generally, $F_j$ may be a mixture of such distributions.) Thus, the bound \eqn{eq-555} applies for a class $j$
with i.i.d. DHR components sizes.

\begin{thm}
\label{thm-inh-subcrit}
Suppose that the infinite-frame system with parameters $\lambda$ and $\pi$ is inherently sub-critical
and Assumption~\ref{cond-second-moment} holds.
Then 
the ML-MFP $x^*$ is the unique proper infinite-frame DE-FP and $x^n(\infty) \Rightarrow x^*$.
\end{thm}

{\em Proof.} Given \eqn{eq-cubcrit}, it is straightforward to observe that ML-MFP $x^*$ is proper (with $x^*_0 \le \bar \rho$).
Let us first suppose that Assumption~\ref{cond-ntriv-add} holds. (We later show how to remove this assumption.)
Consider the sequence $x^n(\infty)$. For all large $n$, $\rho^n = \E x_0^n(\infty) \le \bar \rho + \epsilon < 1$.
Consider any subsequence along which $x^n(\infty) \Rightarrow x(\infty)$ for some $x(\infty)$;
denote $\rho = \E x_0(\infty) = \lim \E x^n_0(\infty) < 1$. Clearly, $x(\infty) \ge x^*$. 
Now, given \eqn{eq-cubcrit}, clearly $\hat \lambda$ and $\hat \pi$ exist, for which the system is also inherently subcritical
with the corresponding proper ML-MFP $\hat x$, and conditions of Theorem~\ref{thm-majorization} hold.
We conclude that Condition~\ref{cond-unique} holds for $\lambda$ (and $\pi$).
From this point on we can repeat the proof of Theorem~\ref{thm-main} to show that $x(\infty) = x^*$ must hold.

If Assumption~\ref{cond-ntriv-add} does not necessarily hold, we can consider a modified system where we
augment the set of job classes by adding a class $\tilde j$
(say with i.i.d. exponential component sizes), and pick parameters $\tilde \lambda$ and $\tilde \pi$ so that:
$\lambda \pi \le \tilde \lambda \tilde \pi$; $|\lambda \pi - \tilde \lambda \tilde \pi |$ is small, so that the inherent subcriticality holds for 
the modified system as well; the ML-MFP $\tilde x^*$ for the modified system is close to (and dominates) $x^*$.
For the modified system the theorem statement does hold. 
We see that, for our original system, 
any subsequential distributional limit $x(\infty)$ of $x^n(\infty)$ must be ``sandwiched'' between $x^*$ and $\tilde x^*$.
Since $\tilde x^*$ can be arbitrarily close to $x^*$, we must have $x^n(\infty) \Rightarrow x^*$.
$\Box$

\begin{thm}
\label{thm-main-majorize}
Consider an infinite-frame (non-truncated)  system. Assume that the set of job classes is such that along with each class it contains all
its subclasses, and moreover $\pi_j >0$ for all classes (and subclasses). 
Suppose the additional Assumptions~\ref{cond-second-moment} and \ref{cond-ntriv-add} hold.
Denote 
$$
\bar \lambda = \sup \{\lambda \ge 0 ~|~ \mbox{corresponding ML-MFP $x^{*,\lambda}$ is proper}\},
$$
$$
\bar \rho = \lim_{\lambda \uparrow \bar \lambda} x_0^{*,\lambda}.
$$
(Note that the case $\bar \rho < 1$ is {\em not} excluded here.) 
Then, as $n\to\infty$,
\beql{eq-conv-main-majorize}
x^n(\infty) \Rightarrow x^{*,\lambda} ~~\mbox{for $0 \le \lambda < \bar \lambda$}
\eeql
and 
\beql{eq-conv-main-majorize-inf}
x^n(\infty) \Rightarrow x^{**} ~~\mbox{for $\lambda > \bar \lambda$}.
\eeql
\end{thm}

{\em Proof.} The proof is an adjusted version of that of Theorem~\ref{thm-main}. We only provide a sketch here.
We pick a $\rho \in [0,\bar \rho)$ and, as in the proof of Theorem~\ref{thm-main}, for each $n$ consider the unique
$\lambda^n$ under which $\rho^n = \E x^n_0(\infty) = \rho$. Then we consider a subsequence along which $\lambda^n \to \lambda$
for some $\lambda$. For this $\lambda$, as in the proof of Theorem~\ref{thm-main}, we must have that any subsequential distributional limit
$x(\infty)$ of $x^n(\infty)$ is such that $x^{*,\lambda} \le x(\infty)$ and $x^{*,\lambda}_0 \le \E x_0(\infty) = \rho$.
By Theorem~\ref{thm-majorization}, Condition~\ref{cond-unique} holds for $x^{*,\lambda}$.
(Indeed, for a $\tilde \lambda$ slightly greater than $\lambda$, proper $x^{*,\tilde \lambda}$ exists, because 
$x^{*,\lambda}_0 \le \rho < \bar \rho$.) 
From this point on we can repeat the proof of Theorem~\ref{thm-main} to show that $x(\infty) = x^{*,\lambda}$ must hold.
We then easily see that $\lambda$ is determined by $\rho$ uniquely, and 
the corresponding function $\lambda(\rho)$ is strictly increasing, mapping $[0,\bar \rho)$ onto $[0,\bar \lambda)$.
From here we obtain claim \eqn{eq-conv-main-majorize}. The claim \eqn{eq-conv-main-majorize-inf} 
follows from the fact that, for $\lambda > \bar \lambda$, the ML-MFP is the infinite element $x^{**}$,
and the ML-MFP is always a lower bound on any subsequential  limit $x(\infty)$ of $x^n(\infty)$.
$\Box$

\section{A ``free'' particle system: Tightness of stationary distributions of centered states}
\label{sec-free}

For a fixed $n$, consider an artificial ``free'' system, where particle locations (server workloads)
are {\em not} regulated at $0$; namely, they evolve in $\R$ rather than in $\R_+$, with each particle keeping moving left at constant rate $-1$
unless/until it jumps right due to a job arrival. Clearly, for such system, the centered process $\mathring x^n(\cdot)$, 
taking values in $\mathring{\cx}^{(n)}$, is well defined, and is itself a Markov process 
(not just a projection of Markov process $x^{(n)}(\cdot)$). 
This process only changes its state (jumps) upon job arrivals. Therefore, the arrival rate 
$\lambda^n$ -- as long as it is positive -- only determines the rate at which changes take place, and  has no impact on the existence/non-existence/form of the process stationary distributions. A version of Theorem~\ref{th-closeness} holds for the free system, 
under an additional technical condition that needs to be introduced to be able to demonstrate
stability (positive Harris recurrence).
Recall that the stability for a given $\lambda^n$ was assumed in Theorem~\ref{th-closeness}. 
For the free system, $\lambda^n$
-- as long as it is positive -- 
is irrelevant for stability/instability and the form of a stationary distribution, but 
the stability does need to be shown.

State $\mathring x^n(\cdot)$ can be equivalently described as $p^n(t) = (w_1(t), \ldots, w_n(t)) \in \R^n$, the same way as in Section~\ref{sec-closeness-proof}, except here $p^n(t)$ does {\em not} contain component $z^n(t)$, which is irrelevant for the free process. 
A subset $B \subset R^n$ is called {\em small} (see, e.g. \cite{Bramson-book}) for the process $p^n(\cdot)$,
 if there exists $\tau >0$ and a finite measure 
$\mathcal M$ on $R^n$, such that, uniformly on the initial state $p^n(0)\in B$, the distribution of $p^n(\tau)$ dominates $\mathcal M$.
(Existence of a small set allows one to use Nummelin splitting to view the process as having a renewal atom.) Note that the property of a set being small for $p^n(\cdot)$
does not depend on $\lambda^n$, as long as $\lambda^n>0$.

\begin{thm}
\label{th-closeness-free}
Suppose the additional Assumptions~\ref{cond-second-moment} and \ref{cond-ntriv-add}  hold.
Suppose further that, for all large $n$, 
\beql{eq-small-set}
\mbox{for any $C>0$ the (closed) set $\{|p^n| \le C\}$ is small for $p^n(\cdot)$.} 
\eeql
Then, there exist $\bar C>0$ and $\bar n$ such that for all $n \ge \bar n$,  $\mathring x^n(\cdot)$ (with any $\lambda^n>0$)
is positive recurrent and, moreover,
\beql{eq-tight-uniform-free}
\E \Phi_1(\mathring x^n(\infty)) \le \bar C.
\eeql
\end{thm}

Condition \eqn{eq-small-set} is not very restrictive -- it can be verified in many cases of interest.
For example, it holds in the case 
when there exists a class $j$ and $\delta>0$, such that the component size distribution $F_j$ dominates a scaled 
(by a small constant) 
Lebesgue measure on $[0,\delta]^{d_j}$.
It is easy to check that the set $\{|p^n| \le C\}$ is small, with any $\tau>0$ and measure $\mathcal M$ defined as follows.
Fix $\epsilon > 0$, and consider set 
$$
S^n \doteq \{\Delta w_i = w_{i+1} - w_i \in [0,\epsilon], ~i=1,2,\ldots, n-(d_j-k_j)-1; ~~w_{i} = w_{n-(d_j-k_j)}, ~i \ge n-(d_j-k_j)\}.
$$
Then, $\mathcal M$ is restricted to $S^n$ and is
defined by a scaled (by a small constant) Lebesgue measure 
on $(\Delta w_1,\ldots, \Delta w_{n-(d_j-k_j)-1}) \in [0,\epsilon]^{n-(d_j-k_j)}$.

{\em Proof of Theorem~\ref{th-closeness-free}.} The proof is a version of that of Theorem~\ref{th-closeness}. We only provide
a sketch, highlighting the differences. 

Let $n$ (and any $\lambda^n>0$) be fixed. Consider a sequence of processes $p^{n,(r)}(\cdot)$, indexed by by $r\uparrow \infty$, with initial states such that 
$|p^{n,(r)}(0)|=r$ and $(1/r) p^{n,(r)}(0) \to q(0)$, where $|q(0)| = 1$. 
(This is a sequence of processes that defines fluid limits of $p^{n}(\cdot)$.)
Then, there exists $T>0$ such that
\beql{eq-fluid-stabil-weak-p}
\frac{1}{r} |p^{n,(r)}(rT)| \Rightarrow 0.
\eeql
Property \eqn{eq-fluid-stabil-weak-p} follows from the following properties of the
fluid limits $q(\cdot)$ of $p^{n}(\cdot)$: $q(\cdot)$ is Lipschitz;
as long as $q_i(t) > q_\ell(t)$, the difference $q_i(t) - q_\ell(t)$ cannot increase;
as long as $|q(t)| >0$, or equivalently, $\min_i q_i(t) < 0$ and $\max_i q_i(t) > 0$, 
\beql{eq-max-min-deriv}
\frac{d}{dt} \min_i q_i(t) - \frac{d}{dt} \max_i q_i(t) \le -\epsilon,
\eeql
for some fixed $\epsilon>0$ (which may depend on $n$); and, therefore,
the fluid limits are such that
 $|q(t)| = 0$ for all $t \ge 1/\epsilon$.
 In turn, the fluid limit properties, except \eqn{eq-max-min-deriv}, are straightforward to obtain.
Property \eqn{eq-max-min-deriv} is obtained using Assumption~\ref{cond-ntriv-add}, for example as follows.
Denote the sets $I_{min}(t)=\argmin_i q_i(t)$ and $I_{max}(t)=\argmax_i q_i(t)$. 
At any time within a sufficiently small fixed (scaled) time interval $[t,t+\Delta t]$, 
 the pre-limit process $(1/r) p^{n,(r)}(r s)$
 is such that if we randomly uniformly pick a pair of particles $i_{min} \in I_{min}$
and $i_{max} \in I_{max}$, the distribution of the jump size of the former (upon a job arrival) 
always dominates that of the latter, with the expectation of the former exceeding that of the latter
at least by some positive constant; this, in turn, follows from the fact that, 
with non-zero probability, the picked particle $i_{min}$ is one of the left-most particles (in the pre-limit process)
and $i_{max}$ is one of the right-most particles (in the pre-limit process) -- and we 
can use Assumption~\ref{cond-ntriv-add} to show that in this case the expected jump size of
$i_{min}$ must be greater than that of $i_{max}$ by at least some positive constant.

Property \eqn{eq-fluid-stabil-weak-p}, along with condition \eqn{eq-small-set}, implies (cf. \cite{Dai95,Bramson-book}) that $p^n(\cdot)$
is positive Harris recurrent (stable), and therefore has unique stationary distribution. 

Finally, given 
the stability of $p^n(\cdot)$
the proof of the uniform bound \eqn{eq-tight-uniform-free}
is essentially same as that of \eqn{eq-tight-uniform} in Section~\ref{sec-closeness-proof}. In fact, it is simpler,
because the dynamics of $p^n(\cdot)$ for the ``free'' system is simpler, -- there is no regulation boundary $z(t)$ on the left,
therefore the particles never hit it, and the only process transitions are those occurring upon job arrivals.
$\Box$

The free system, in addition to being of independent interest, has the following natural connection to the stability 
(Theorem~\ref{thm-crit-load-finite}) 
of our original infinite-frame system. Suppose the conditions of Theorem~\ref{th-closeness-free} hold, including condition \eqn{eq-small-set}.
Denote by $\bar \lambda^n$ the unique value of $\lambda^n$, 
for which the steady-state average drift of the mean $\bar x^n(\cdot)$ of the (non-centered) free process $x^n(\cdot)$ is $0$. 
(Such value is indeed unique, because, clearly, the drift has the form $-1 + a \lambda^n$, where constant $a>0$ is determined by the
stationary distribution of $\mathring x^n(\cdot)$; so, $\bar \lambda^n = 1/a$.) Then it is not hard to show that this $\bar \lambda^n$
is exactly equal to the $\bar \lambda^n$ in Theorem~\ref{thm-crit-load-finite}. (This can be done by combining 
arguments in the proofs of Theorems~\ref{thm-crit-load-finite} and \ref{th-closeness-free}.) 
We note, however, that this conclusion requires the additional technical condition \eqn{eq-small-set}, 
while Theorem~\ref{thm-crit-load-finite} does {\em not} require it.

\bibliographystyle{abbrv}

\end{document}